\documentclass[a4paper,10pt]
{amsart}
\usepackage[foot]{amsaddr}
\newtheorem{theorem}{Theorem}[section]

\theoremstyle{definition}

\theoremstyle{remark}

\usepackage[margin=2.7cm]{geometry}
\usepackage{systeme}
\usepackage{graphicx}
\usepackage[]{units}
\usepackage{color}
\usepackage{mathrsfs}
\usepackage{bm}
\usepackage[ansinew]{inputenc}
\usepackage{float}
\usepackage{hyperref}
\numberwithin{equation}{section}

\usepackage{subcaption}
\usepackage{cases}
\usepackage{enumitem}
\usepackage{etoolbox}
\usepackage{algorithm, algorithmicx, algpseudocode}
\usepackage[superscript,biblabel]{cite}

\AfterEndEnvironment{problem}{\noindent\ignorespaces}

\newcommand{\R}{{\rm I\!R}}  
\newcommand*\interior[1]{\mathring{#1}}
\newcommand\norm[1]{\left\lVert#1\right\rVert}

\DeclareMathOperator*{\argmin}{argmin}

\newtheoremstyle{noIndent}
  {}
  {}
  {\itshape}
  {}
  {\bfseries}
  {.}
  {.5em}
  {}
\theoremstyle{noIndent}
\newtheorem{problem}{Problem}[section]

\begin{document}

\title[UMBERTO EMIL MORELLI et al]{A numerical approach for heat flux estimation in thin slabs continuous casting molds using data assimilation}
\thanks{Funded by the European Union's Horizon 2020 research and innovation programme under the Marie
Skaodowska-Curie Grant Agreement No. 765374.
It also was partially supported by the Ministry of Economy, Industry and Competitiveness through the Plan Nacional de I+D+i (MTM2015-68275-R), by the Agencia Estatal de Investigacion through project [PID2019-105615RB-I00/ AEI / 10.13039/501100011033], by the European Union Funding for Research and Innovation - Horizon 2020 Program - in the framework of European Research Council Executive Agency:  Consolidator Grant H2020 ERC CoG 2015 AROMA-CFD project 681447 "Advanced Reduced Order Methods with Applications in Computational Fluid Dynamics" and  INDAM-GNCS project "Advanced intrusive and non-intrusive model order reduction techniques and applications", 2019.}

\author{Umberto Emil Morelli\textsuperscript{1,2,3*}}
\email{umbertoemil.morelli@usc.es}

\author{Patricia Barral\textsuperscript{1,2}}
\author{Peregrina Quintela\textsuperscript{1,2}}
\author{Gianluigi Rozza\textsuperscript{3}}
\author{Giovanni Stabile\textsuperscript{3}}

\address{\textsuperscript{1}Universidade de Santiago de Compostela, Santiago de Compostela, Spain}
\address{\textsuperscript{2} Technological Institute for Industrial Mathematics (ITMATI), Santiago de Compostela, Spain}
\address{\textsuperscript{3}Scuola Internazionale Superiore di Studi Avanzati (SISSA), Trieste, Italy}

\date{}

\dedicatory{}

\begin{abstract}
    In the present work, we consider the industrial problem of estimating in real-time the mold-steel heat flux in continuous casting mold.
    We approach this problem by first considering the mold modeling problem (direct problem).
    Then, we plant the heat flux estimation problem as the inverse problem of estimating a Neumann boundary condition having as data pointwise temperature measurements in the interior of the mold domain.
    We also consider the case of having a total heat flux measurement together with the temperature measurements.
    We develop two methodologies for solving this inverse problem.
    The first one is the traditional Alifanov's regularization, the second one exploits the parameterization of the heat flux.
    We develop the latter method to have an offline-online decomposition with a computationally efficient online part to be performed in real-time.
    In the last part of this work, we test these methods on academic and industrial benchmarks. 
    The results show that the parameterization method outclasses Alifanov's regularization both in performance and computational cost.
    Moreover, it proves to be robust with respect to the measurements noise.
    Finally, the tests confirm that the computational cost is suitable for real-time estimation of the heat flux.
\end{abstract}
\keywords{Inverse Problem, Heat Transfer, Continuous Casting, Real-time, Data Assimilation, Boundary Condition Estimation}

\maketitle

\section{Introduction}
\label{section:MoldControl}

\begin{figure}[!htb]
    \begin{subfigure}[htb]{0.5\linewidth}
        \centering
        \includegraphics[width=0.6\textwidth]{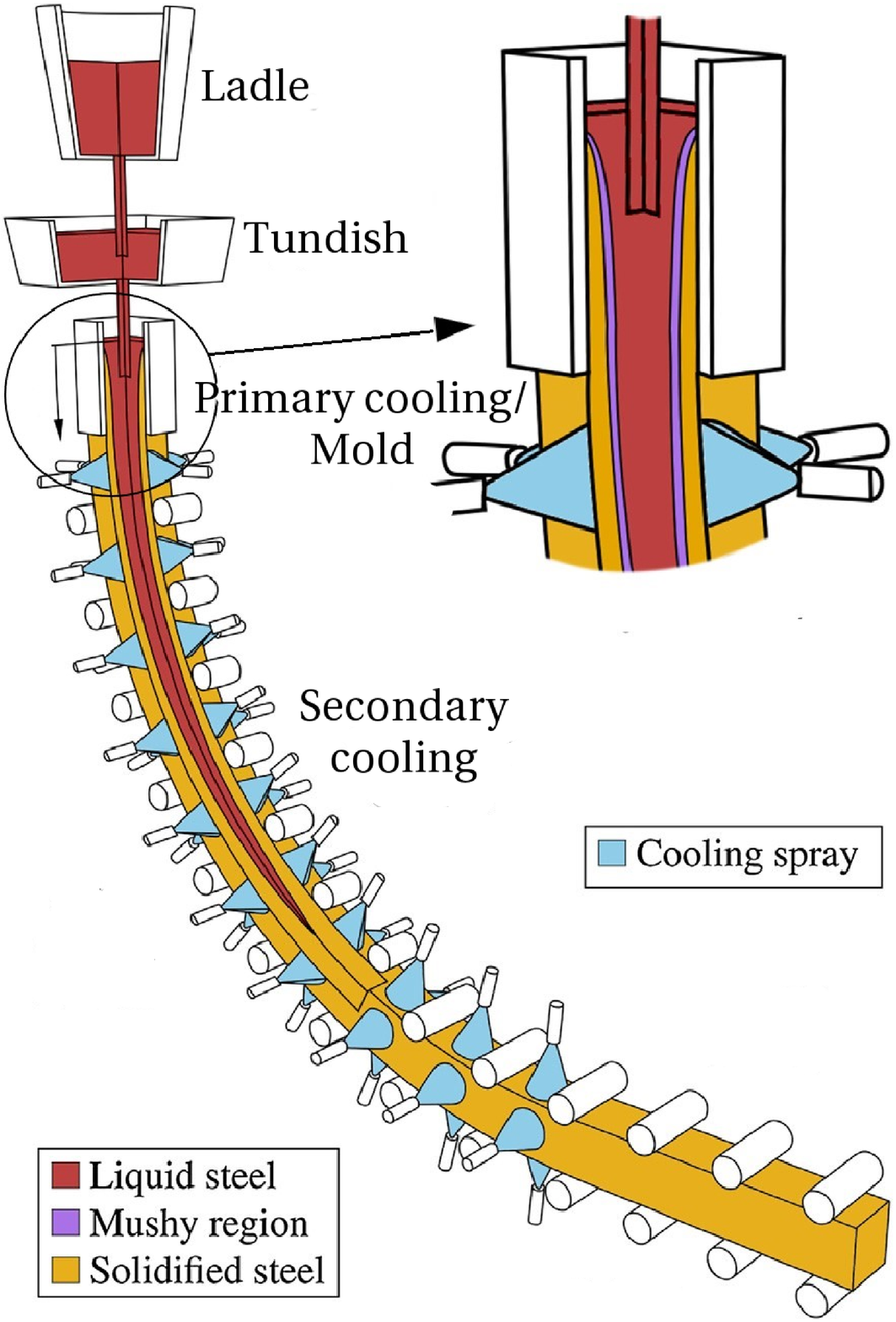}
        \caption{CC process.}
    \end{subfigure}%
    \begin{subfigure}[htb]{0.5\linewidth}
        \centering
        \includegraphics[width=0.9\textwidth]{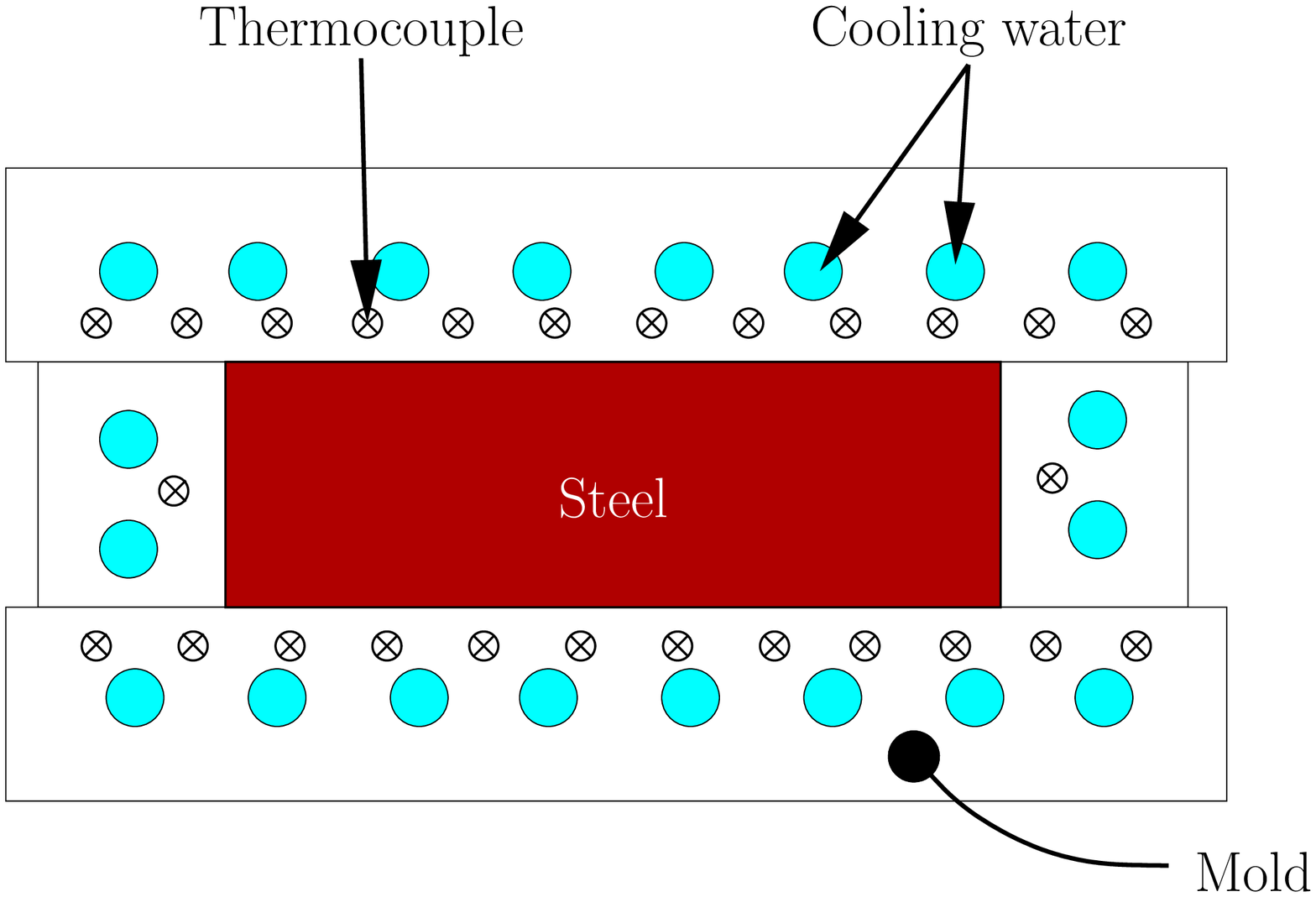}
        \caption{Mold section.}
    \end{subfigure}%
    \caption{(A) Schematic overview of the continuous casting process (adapted from Klimes et al\cite{Klimes2015}).
    (B) Schematic of a horizontal section of the mold (the casting direction is perpendicular to the image).}
\label{fig:castingSchematic}
\end{figure}

Continuous Casting (CC) of steel is presently the most used process to produce steel worldwide.
For example, in 2017, 96\% of the steel was produced by CC.\cite{WorldSteel2018}
This industrial process is not new at all.
In fact, continuous casters as in Figure~\ref{fig:castingSchematic} have been used for many decades now.
Consequently, the process has undergone a long sequence of improvements driven by experience of the commercial operators and, more recently, numerical simulations.\cite{Thomas2018}

We can summarize the CC process as follows.
The metal is liquefied and then tapped into the ladle. 
When it is at the correct temperature, the metal goes into the tundish.
In the tundish, the metal flow is regulated and smoothed.
Through the Submerged Entry Nozzle~(SEN), the metal is drained into a mold.
The role of the mold in the CC process is to cool down the steel until it has a solid skin which is thick and cool enough to be supported by rollers.

At the outlet of the mold, the metal is still molten in its inner region, thus a secondary cooling section follows the mold.
Supported by rollers, it is cooled until complete solidification by directly spraying water over it.
At the end of this secondary cooling region, the casting is completed.
To be ready for its final application, the strand generally continues through additional mechanisms which may flatten, roll or extrude the metal into its final shape.
This is just a brief overview on the CC process.
We refer to Irving's monograph for a detailed description.\cite{Irving1993}

In this work, we focus on CC of thin slabs.
Slabs are cataloged thin when their thickness is smaller than 70~mm, while their width is between 1 and 1.5~m, in general.
Thanks to the small thickness, the solidification in the slab is relatively fast, consequently the casting speed is generally high, between 7 and 14 meters per minute.

Thin slab molds are made of four different plates: two wide plates and two lateral plates, all made of copper (see Figure~\ref{fig:castingSchematic} (b)).
In general, lateral plates can be moved or changed to modify the slab section dimensions.
The geometry of these plates is more complex than one can expect: they have drilled channels where the cooling water flows, slots in the outside face for thermal expansion, thermocouples, and fastening bolts.
To compensate the shrinkage of the slab with the cooling and minimize the gap, the molds are tapered.
Moreover, the upper portion of the mold forms a funnel to accommodate the SEN.

Several phenomena related to steel flow, solidification, mechanics and heat transfer appear in the mold region.
This complexity makes the mold the most critical part of the CC process.
Thus, when running a continuous caster, productivity and safety issues must be addressed at the mold.

Regarding quality, the presence of imperfections on the external surface of the casted piece (cracks, inclusions, etc.) must be avoided.
In fact, since casted pieces are generally laminated in later productions, surface defects would become evident affecting also the mechanical properties of the final products.

However, quality control is not the only issue arising at the mold.
Due to the creation of the solid skin, a frequent problem arising during CC is the sticking of the steel to the mold.
After the detection of a sticking, the casting speed is reduced to reestablish the desired metal flow before restoring the nominal casting speed.
This affects the product quality and the productivity of the caster.
Also, if not detected on time, it can lead to dangerous events forcing the shutdown of the caster.

Less frequent but more catastrophic events are the liquid break-out and the excessive increase of the mold temperature.
The former is due to a non-uniform cooling of the metal with the skin being so thin to break.
The latter is generally considered as the most dangerous event in a casting plant.
In fact, if the mold temperature is high enough to cause the boiling of the cooling water, we have a dramatic decrease in the heat extraction.
Then, the temperature in the mold quickly rise, that could cause the melting of the mold itself.
Both these incidents are very dangerous and costly.
In fact, they generally require the shutdown of the caster, the substitution of expensive components and an extended turnaround.

For all these reasons, the early detection of problems in the mold is crucial for a safe and productive operation of continuous casters.
Their detection becoming more difficult as casting speed (thus productivity) of the casters increases.

Since, continuous caster has been running for decades, operators already faced all these problems.
To have insight of the scenario in the mold, they provided the molds with measuring equipment.
In particular, they measure the pointwise temperature of the mold by thermocouples (see Figure~\ref{fig:castingSchematic} (b)) and the cooling water temperature as well as its flow at the inlet and outlet of the cooling system.

The way CC operators use the data coming from the measurement equipment is the following.
The thermocouples' temperatures are used to have insight of the mold temperature field.
On the other hand, the water temperature rise is used to approximate the heat extracted from the steel.

This approach allowed to run continuous casters for decades.
Nevertheless, it has several drawbacks: it relies on the experience of operators, gives very limited information about the heat flux at the mold-slab interface, and is customized for each geometry so it requires new effort to be applied to new designs.
So, a new tool for analyzing the mold behavior is necessary.

We begin by reporting that CC operators consider that knowing the local heat flux between mold and slab is the most important information in analyzing the casting in this region.
By considering the mold itself to be our domain and focusing our interest in its thermal behavior, the mold-slab heat flux can be seen as a Neumann Boundary Condition (BC) in the model.
To compute its value, we pose the following  inverse problem: from the temperature measurements provided by the thermocouples, estimate the boundary heat flux at the mold-slab interface.

In general, this is a complex problem which can be divided into three different ones but with related ingredients:
\begin{itemize}
  \item Accurate modeling of the thermal problem in the physical mold.
  \item Solution of the inverse problem of estimating the heat flux.
  \item Reduction of the computational cost of the inverse problem solution to achieve real-time computation.
\end{itemize}
In this work, we address in detail the above three problems giving an overview on the state of the art and presenting our approach and contribution for their solution.

\section{Mold Thermal Model}

In this section, we give a description of the physical phenomena that occur in the mold region of a caster, giving an overview on previous efforts in modeling them.
Then, we present the physical problem that we will consider in the present work.
Finally, we present the mathematical model we use in the rest of the paper.

\subsection{Physical Problem}
Going from the inside to the outside of the mold, we encounter several physical phenomena.
In the inner part of the mold, we have the liquid pool of steel.
There, we have a molten metal flow with dispersed argon bubbles and inclusion particles.
All around the liquid pool, we encounter the solid skin and, in between, the mushy region.
Here, the steel changes phase undergoing solidification.
Between the steel and the mold, there is a thin layer of flux powder which is liquid close to the steel and solid where in contact with the mold.
Finally, we encounter the mold which is surrounding the flux powder (in case of not perfect casting, we can also have an air gap between the mold and the slab).

The mold is composed of a solid (copper) region and a liquid region (water) representing its cooling system.
In the copper, we have heat conduction due to the temperature gradients.
In the water, we have an incompressible flow inside tubes.
To prevent the water from boiling, it is pumped at a very high pressure and flow rate.
Therefore, a turbulent flow with high Nusselt number is expected.

According to the previous description, the main physical phenomena for CC include\cite{AISE2003}:
\begin{itemize}
  \item Fully-turbulent, transient fluid motion in a complex geometry (SEN, strand liquid pool), affected by dispersed particles and thermal buoyancy.
  \item Thermodynamic reactions.
  \item Multiphase flow and heat transport.
  \item Dynamic motion of free liquid surfaces and interfaces.
  \item Thermal, fluid and mechanical interactions between solids, liquids and gases.
  \item Heat transfer in solids.
  \item Distortion and wear of the mold.
  \item Solidification of steel shell.
  \item Shrinkage of the solidifying steel shell.
  \item Crack formation.
\end{itemize}

Due to its complexity, the literature on CC mold modeling is extensive.
For each physical phenomenon in the previous list, we have at least a dedicated model and several investigations.
Just to name the most relevant works, Meng and Thomas\cite{Thomas2003} investigated the modeling of transient slag-layer phenomena in the
steel-mold gap.
Thomas et al\cite{Thomas2011} studied the steel flow and solidification in the mold including the transport of superheat with the turbulent transient flow of molten steel, surface level fluctuations, and the transport and entrapment of inclusion particles.
On the other hand, a detailed description of the solidification in casting is available in the Stefanescu's monograph,\cite{Stefanescu2015} while a review on the initial solidification models was done by Wang.\cite{Wang2017} 

As mentioned the literature on the subject is extensive.\cite{Goldschmit1999, Williams1979, Koric2006}
Thus, we redirect the interested reader to the review articles by Thomas et al\cite{Thomas2002,Thomas2018}.

\subsection{Specific Physical Problem}
As discussed in Section~\ref{section:MoldControl}, our goal is to study the real-time behavior of CC molds and, in particular, to compute the mold-slab heat flux.
One way to compute it could be to simulate all the phenomena discussed above from the SEN to the secondary cooling region.
However, the resulting model would be quite complex and computationally expensive to deal with, especially for real-time applications.
Then, we discard this option.

The fully experimental approach is also not feasible since it is not possible to make direct measurements in the solidification region. 
The only measurements available are temperature measurements by thermocouples that are buried inside the mold plates and cooling water temperature measurements.
Then, our approach to study the real-time behavior of continuous casting molds and be able to compute the mold-slab heat flux is to solve an inverse problem having these data as control data.
In the rest of this section, we describe the mold thermal model that we use in the present investigation and the related assumptions.

In modeling the thermal behavior of the mold, we consider the following well established assumptions:
\begin{itemize}
  \item The copper mold is assumed a homogeneous and isotropic solid material.
  \item The cooling water temperature is known.
  \item The thermal expansion of the mold and its mechanical distortion are negligible.
  \item The material properties are assumed constant.
  \item The boundaries in contact with air are assumed adiabatic.
  \item No boiling in the water is assumed.
  \item The heat transmitted by radiation is neglected.
\end{itemize}
Since we want to have solution in real-time (e.g., at each second) and the casting speed is of few meters per minute, we consider steady-state models.
Moreover, we only consider 3D mold models because we are interested in the heat flux in all the mold-slab interface.

As a final remark, the running parameters of the cooling system and its geometry ensure a fully developed turbulent flow.
In fact, these molds are equipped with a closed loop cooling system where the water is pumped at a high pressure.
The average velocity in each cooling channel is approximately 10 m/s, the diameter being approx. 10 mm.
Thus, the Reynolds number in the cooling system is around 10\textsuperscript{5}, which ensures a turbulent flow.

Thanks to the high Reynolds number of the flow, we can further assume that the cooler and hotter water molecules are well mixed.
Consequently, the temperature in each section of the cooling channel is approximately constant.
Moreover the water is pumped in a closed circuit, so we can assume that the water flow rate is constant.
In turn, since the channels have constant section, the velocity of the fluid is also uniform and constant (plug flow).

Then, we focus our attention on the following model:
\begin{enumerate}
 \setcounter{enumi}{0}
 \item\label{model:steadyHeatConduction} The computational domain is only composed of the (solid) copper mold.
    We consider a steady-state three-dimensional heat conduction model with a convective BC in the portion of the boundary in contact with the cooling water.
    The water temperature is known at the inlet and outlet of the cooling system.
    The water temperature is assumed to be linear.
\end{enumerate}

\subsection{Computational Domain and Notation}

\begin{figure}[htb]
    \centering
    \includegraphics[width=0.6\textwidth]{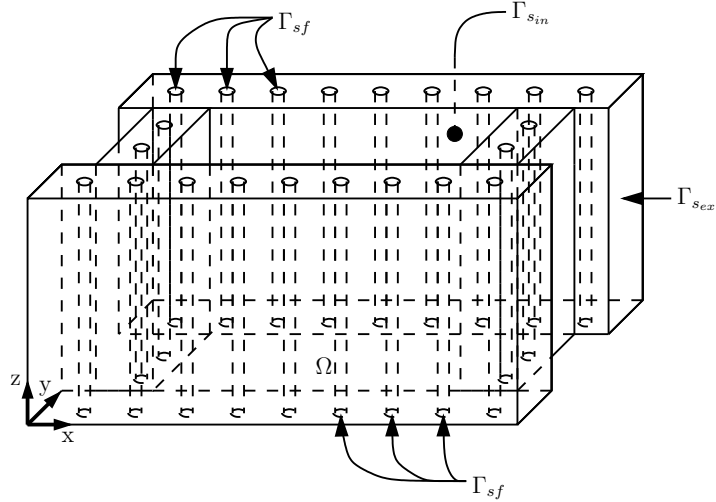}	
    \caption{Schematic of the mold domain, $\Omega$, and its boundaries.}
\label{fig:directProblem_schematicDomain}
\end{figure}

Consider a solid domain, $\Omega$, which is an open Lipschitz bounded subset of $\R^3$,  with smooth boundary $\Gamma$ (see Figure~\ref{fig:directProblem_schematicDomain}).
Let $\Gamma = \Gamma_{s_{in}}\cup\Gamma_{s_{ex}} \cup \Gamma_{sf}$ where $\interior{\Gamma}_{s_{in}}$, $\interior{\Gamma}_{s_{ex}}$ and $\interior{\Gamma}_{sf}$ are disjoint sets.
The Eulerian Cartesian coordinate vector is denoted by $\mathbf{x}\in \Omega$ and $\mathbf{n}(\mathbf{x})$ the unit normal vector that is directed outwards from $\Omega$.

In this setting, $\Omega$ corresponds to the region of the space occupied by the mold.
The interface between the mold and the cooling system is denoted by $\Gamma_{sf}$.
While, $\Gamma_{s_{in}}$ is the portion of the mold boundary in contact with the solidifying steel.
Finally, we denote the remaining part of the mold boundary with $\Gamma_{s_{ex}}$.

\subsection{Mathematical Model}\label{section:directProblem_assumptions}

We shall assume all along the following assumptions on the data:
\begin{enumerate}[start=1,label={(H\arabic*)}]
        \setcounter{enumi}{0}
  \item\label{ass:steadyness} The process is assumed to be steady-state. 
  \item\label{ass:thermalConductivityPositivness} The thermal conductivity is constant and strictly positive: $k\in\R^+$.
  \item\label{ass:heatFluxInL2} The mold-steel heat flux, $g$, belongs to $L^2 (\Gamma_{s_{in}})$.
  \item There is no heat source inside the mold domain.
  \item\label{ass:steadyHeatConduction_heatTransferCoeff} The heat transfer coefficient is known, constant and strictly positive: $h\in \R^+$.
  \item\label{ass:steadyHeatConduction_knownCoolingTemperature} The cooling water temperature, $T_f$, is known and belongs to $L^2(\Gamma_{sf})$. 
\end{enumerate}

Under the assumptions \ref{ass:steadyness}-\ref{ass:steadyHeatConduction_knownCoolingTemperature}, we propose the following three-dimensional, steady-state, heat conduction model
\begin{problem}{}
Find $T$ such that
\begin{equation}
   - k \Delta T = 0, \text{ in }\Omega,
\label{eq:3DhcModelSteady_laplacian}
\end{equation}
with BCs
\begin{numcases}{}
  -k \nabla T \cdot \mathbf{n} = g		&       on $\Gamma_{s_{in}}$, \label{eq:3DhcModelSteady_laplacian_BC1}\\
  -k \nabla T \cdot \mathbf{n} = 0		&       on $\Gamma_{s_{ex}}$, \label{eq:3DhcModelSteady_laplacian_BC2}\\
  -k \nabla T \cdot \mathbf{n} = h(T - T_f)	&       on $\Gamma_{sf}$\label{eq:3DhcModelSteady_laplacian_BC3}.
\end{numcases}
\label{prob:3DhcModelSteady}
\end{problem}

We recall that for this problem the following result is well established (see Nittka, Theorem 3.14)\cite{Nittka2011}:
\begin{theorem}
  Under the assumptions \ref{ass:heatFluxInL2},\ref{ass:steadyHeatConduction_heatTransferCoeff} and \ref{ass:steadyHeatConduction_knownCoolingTemperature}, the solution to Problem~\ref{prob:3DhcModelSteady} exists and is unique in $H^1(\Omega)$.
  Moreover, there exists a $\gamma > 0$ such that the solution to Problem~\ref{prob:3DhcModelSteady} belongs to $C^{0,\gamma}({\Omega})$.
  \label{theo:Nittka2011_3.14}
\end{theorem}

As a final remark, we recall (see Raymond, Theorem 3.3.6)\cite{Raymond2013},
\begin{theorem}
    If $g$ and $T_f$ belong to $L^s(\Gamma_{s_{in}})$ and $L^s(\Gamma_{sf})$ respectively, with $s > 2$, then the solution $T$ to Problem~\ref{prob:3DhcModelSteady} belongs to $C(\bar{\Omega})$ and
    \begin{equation}
        \norm{T}_{C(\bar{\Omega})} \leq C \left( \norm{g}_{L^s(\Gamma_{s_{in}})} + \norm{T_f}_{L^s(\Gamma_{sf})} \right),
    \end{equation}
    where the constant $C$ is independent of $h$.
\end{theorem}

Regarding the numerical solution of Problem~\ref{prob:3DhcModelSteady}, we use the finite volume method for its discretization.
Given a tessellation $\mathcal{T}$ of the domain, $\Omega$, we write the discrete unknown $(T_C)_{C\in \mathcal{T}}$ as the real vector $\mathbf{T}$, belonging to $\R^{N_h}$ with $N_h = \text{size}(\mathcal{T})$.
Then, we write the discretized problem as the linear system.
\begin{equation}
    \mathbf{A} \mathbf{T} = \mathbf{b}, 
\label{eq:benchmark3DhcModelSteady_discreteDirectProblemLinSys}
\end{equation}
where $\mathbf{A}$ is the stiffness matrix and $\mathbf{b}$ the source term.
The value of each element of $\mathbf{A}$ and $\mathbf{b}$ depends on the particular finite volume scheme for the discretization and the mesh used.
Since our problem is a classic diffusion problem, we refer for further details regarding the finite volume discretization to the Eymard's monograph.\cite{Eymard2000}

\section{Inverse Problem}\label{sec:inverseProblem}

This section is devoted to the formulation and the study of an inverse problem related to Problem~\ref{prob:3DhcModelSteady}.
We consider two different inverse problems.
In the first one, the available data are the thermocouples' measurements only.
In the second one, the total heat flux measurement is available together with the thermocouples' measurements.

\subsection{State of the Art}

The literature on inverse heat transfer problems is vast.\cite{Ling2003, Loulou2006, Jin2007, Huang1996}
We refer to Alifanov's\cite{Alifanov1988}, Orlande's\cite{Orlande2010}, Beck and Clair's\cite{Beck1985} and Chang's\cite{Chang2017} works for a detailed review.
This literature also includes the particular problem of computing the mold-slab heat flux from temperature measurements in the mold.\cite{Ranut2012, Udayraj2017, Mahapatra1991}
From a mathematical point of view, the present problem is the estimation of a Neumann BC (the heat flux) having as data pointwise measurements of the state inside the domain.
Such problems were also addressed in investigations not related to heat transfer.\cite{Raymond2013, Vitale2012, Huang1998}
Due to the vastness of the literature on the subject, we merely report the most relevant works on this subject.

Research in inverse heat transfer problems started in the 50s.
It was driven by the interest in knowing thermal properties of heat shields and heat fluxes on the surface of space vehicles during re-entry.
From a heuristic approach in the 50s, researchers moved to a more mathematically formal approach.
In fact, in the 60s and 70s, most of the regularization theory that we use nowadays to treat ill-posed problems was developed.\cite{Tikhonov1963, Alifanov1988, Beck1968, Beck1970, Chen1976}
Here, we discuss in general the most popular methodologies used for the solution of inverse heat transfer problems.

Traditionally in estimating the boundary heat flux in CC molds, a heat flux profile is selected, then by trial and error it is adapted to match the measured temperatures.\cite{Mahapatra1991}
Pinhero et al\cite{Pinheiro2000} were the first to use an optimal control framework and regularization methods.
They used a steady-state version of the 2D mold model proposed by Samarasekera and Brimacombe\cite{Samarasekera1982} and parameterized the heat flux with a piecewise constant function.
Finally, they used Tikhonov's regularization for solving the inverse problem and validated the results with experimental measurements.
A similar approach was used by Rauter et al.\cite{Rauter2008}
Ranut et al\cite{Ranut2011, Ranut2012} estimated the heat flux transferred from the solidifying steel to the mold wall both in a 2D and 3D domain.
They used a steady-state heat conduction model for the mold and parameterized the heat flux with a piecewise linear profile in 2D and symmetric cosine profile in 3D.
For the solution of the inverse problem, they used the Conjugate Gradient Method (CGM) and a mixed GA-SIMPLEX algorithm\cite{Nelder1965} in 2D while in 3D they only used the GA-SIMPLEX algorithm.
Their results were also tested with experimental data.

Hebi et al\cite{Man2004, Hebi2006} attempted to estimate the solidification in CC round billets by using a 3D transient heat conduction model in the strand and the mold with a Robin condition at the mold-strand interface.
Then, they posed the following inverse problem: find the inverse heat transfer coefficient between mold and strand such that the distance between measured and computed temperatures at the thermocouples is minimal.
They assumed the inverse heat transfer coefficient to be piecewise constant.
Then, by using sensitivity coefficients, each piece was iteratively adapted to match the measured temperature.
To allow convergence, a relaxation factor was introduced in between the iterations.
They validated the results with plant measurements without obtaining a good agreement.
A similar approach was used by Gonzalez et al\cite{Gonzalez2003} and Wang et al,\cite{Wang2016, ZhangWang2017, Hu2018, Tang2012} the latter using a Neumann condition at the mold-strand interface.

Wang and Yao\cite{WangYao2011} used the aforementioned inverse problem solution technique to estimate the inverse heat transfer coefficient for CC round steel billets.
Then, they used the results obtained to train a Neural Network (NN) for on-line computation.
Similarly, Chen et al\cite{Chen2014} used the fuzzy inference method for estimating the mold heat flux.
They modeled the mold with a 2D steady-state heat conduction model in the solid and parameterized the boundary heat flux.
They tested the results on a numerical benchmark obtaining a good agreement.

Yu and Luo\cite{Yu2015} considered a 2D vertical section of a slab and the solidification problem therein.
They developed a modified Levenberg-Marquardt method to estimate the inverse heat transfer coefficient in the secondary cooling region from temperature measurements on the surface of the slab.

Udayraj et al\cite{Udayraj2017} applied CGM with adjoint problem for the solution of the inverse steady-state 2D heat conduction problem, this methodology was first proposed by Alifanov\cite{Alifanov1988} for the regularization of boundary inverse heat transfer problems.
By using this method there is no need of parameterizing the heat flux.
However, the method underestimates the heat flux away from the measurements.
To overcome this issue, the authors proposed to average the computed heat flux at each step and use the uniform averaged value as initial estimation for the following step.
Similarly, Chen et al\cite{Chen2009} tackled the problem of estimating the steady boundary heat flux in 2D circular CC supporting rollers based on temperature measurements inside the domain.
For its solution, they used the CGM proposed by Alifanov\cite{Alifanov1988}.

We conclude this section by describing previous works that are related to the present research but not to CC.
Ambrosi et al\cite{Ambrosi2006, Vitale2012} studied the mathematical formulation of the force traction microscopy problem.
This inverse problem consists in obtaining the boundary stress field on a portion of the boundary (Neumann BC) based on the pointwise measurement of the displacement (state variable) inside the domain.
The similarity with the present research is the presence of pointwise observations and a boundary inverse problem with (linear) elliptic direct problem.
In Vitale et al\cite{Ambrosi2006}, they stated the 2D direct problem and the related inverse problem in the standard optimal control framework due to Lions\cite{Lions1971} for which the unknown BC is the distributed boundary control.
Then in a following work\cite{Vitale2012}, they extended the formulation to the 3D linear elasticity model proving existence and uniqueness of the optimal control solution.

Our contribution to the literature is to develop a novel method for solving the 3D inverse heat transfer problem in CC molds that can achieve real-time performances.
Moreover, we design two benchmark cases for this application and we use them to compare the performances of the proposed method with the classical Alifanov's regularization.

\subsection{Inverse Problem with Thermocouples' Measurements}\label{sec:inverseProblem_soloThermocouples}
The inverse problem we want to solve is to estimate the heat flux $g$ capable of reproducing the measured temperatures at the thermocouples' points.
This can be stated as an optimal control problem with pointwise observations.

We introduce the following notation.
Let $\Psi:=\{\mathbf{x}_1, \mathbf{x}_2 , \dots, \mathbf{x}_M \}$ be a collection of points in $\Omega$.
We define the application $\mathbf{x}_i \in \Psi\rightarrow \hat{T}(\mathbf{x}_i)\in \R^+$, $\hat{T}(\mathbf{x}_i)$ being the experimentally measured temperature at $\mathbf{x}_i \in \Psi$.
Moreover, let $G_{ad}$ be a bounded set in $L^2(\Gamma_{s_{in}})$.

Using a least square, deterministic approach, we state the inverse problem as
\begin{problem}{(\textbf{Inverse})}
  Given $\{ \hat{T}(\mathbf{x}_i) \}_{i=1}^M$, find the heat flux $g \in G_{ad}$ that minimizes the functional $J_1:L^2(\Gamma_{s_{in}}) \rightarrow \R^+$,
  \begin{equation}
    J_1[g]:=\frac{1}{2}\sum^M_{i=1} [T[g](\mathbf{x}_i) - \hat{T}(\mathbf{x}_i)]^2,
    \label{eq:inverseProblem_costFunction}
  \end{equation}
  where $T[g](\mathbf{x}_i)$ is the solution of Problem~\ref{prob:3DhcModelSteady} at points $\mathbf{x}_i$, for all $i=1,2,\dots,M$.
  \label{inverseProblem_3DhcModelSteady}
\end{problem}
Notice that, thanks to Theorem~\ref{theo:Nittka2011_3.14} the state variable $T$ is continuous in $\Omega$, then its value at pointwise observations is well-defined.

We now introduce the sensitivity problem related to Problem~\ref{prob:3DhcModelSteady}.
We derive it by perturbing in Problem~\ref{prob:3DhcModelSteady} the heat flux $g\rightarrow g+\delta g$, causing a variation of the temperature field, $T[g]\rightarrow T[g]+\delta T[\delta g]$.
Subtracting Problem~\ref{prob:3DhcModelSteady} from the obtained problem, we have
\begin{problem}{(\textbf{Sensitivity})}
Find $\delta T$ such that
\begin{equation}
   -k\Delta \delta T [\delta g]= 0, \quad      \text{in }\Omega,
\label{eq:benchmark3DhcModelSteady_sensitivity}
\end{equation}
with BCs
\begin{numcases}{}
  -k \nabla \delta T[\delta g] \cdot \mathbf{n} = \delta g       	&       on $\Gamma_{s_{in}}$,\\
  -k \nabla \delta T[\delta g]\cdot \mathbf{n} = 0                   	&       on $\Gamma_{s_{ex}}$,\\
  -k \nabla \delta T[\delta g] \cdot \mathbf{n} = h(\delta T[\delta g]) 	&       on $\Gamma_{sf}$.
\end{numcases}
\label{prob:benchmark3DhcModelSteady_sensitivity}
\end{problem}
Then, it is verified that $T[g + \delta g] = T[g] + \delta T[\delta g]$.
Besides, $\delta T$ is linear: $\delta T[\delta g_1 + \delta g_2] = \delta T[\delta g_1] + \delta T[\delta g_2]$.

We now derive in a formal way the adjoint of Problem~\ref{inverseProblem_3DhcModelSteady}.
Firstly, we multiply (\ref{eq:3DhcModelSteady_laplacian}) by a Lagrange multiplier $\lambda$.
Then, we integrate over $\Omega$ and add it to (\ref{eq:inverseProblem_costFunction}) obtaining
\begin{equation}
   \mathcal{L}[g, \lambda] = \frac{1}{2} \sum_{i=1}^M (T[g](\mathbf{x}_i) - \hat{T}(\mathbf{x}_i))^2 + \int_{\Omega} k\Delta T[g](\mathbf{x}) \lambda(\mathbf{x}) d\mathbf{x}.
\label{eq:benchmark3DhcModelSteady_adjointDerivation1}
\end{equation}
To compute the Fr{\'e}chet derivative with respect to $g$ of $\mathcal{L}[g, \lambda]$, denoted by $d\mathcal{L}_g[\delta g, \lambda]$, we first write
\begin{equation}
   \mathcal{L}[g + \delta g, \lambda]-  \mathcal{L}[g, \lambda] = \sum_{i=1}^M  \delta T[\delta g](\mathbf{x}_i)(T[g](\mathbf{x}_i) + \frac{1}{2}\delta T[\delta g](\mathbf{x}_i) - \hat{T}(\mathbf{x}_i)) + \int_{\Omega} k \lambda(\mathbf{x}) \Delta \delta T[\delta g](\mathbf{x}) d\mathbf{x}.
\end{equation}
The Fr{\'e}chet derivative of $\mathcal{L}$ is then obtained by neglecting the second order terms
\begin{equation}
    d\mathcal{L}[\delta g, \lambda] = \sum_{i=1}^M  \delta T[\delta g](\mathbf{x}_i)(T[g](\mathbf{x}_i) - \hat{T}(\mathbf{x}_i)) + \int_{\Omega} k \lambda(\mathbf{x}) \Delta \delta T[\delta g](\mathbf{x}) d\mathbf{x}.
\end{equation}
Finally, integrating the second addition of the previous equality twice by parts and applying the BCs of Problem~\ref{prob:benchmark3DhcModelSteady_sensitivity}, we can write 
\begin{equation}
\begin{aligned}
    d\mathcal{L}[\delta g, \lambda] =& \sum_{i=1}^M  \delta T[\delta g](\mathbf{x}_i)(T[g](\mathbf{x}_i) - \hat{T}(\mathbf{x}_i)) + \int_{\Omega} k \Delta\lambda(\mathbf{x})  \delta T[\delta g](\mathbf{x}) d\mathbf{x} - \int_{\Gamma_{s_{in}}\cup\Gamma_{s_{ex}}\cup\Gamma_{sf}} k \delta T[\delta g](\mathbf{x}) \nabla \lambda(\mathbf{x}) \cdot \mathbf{n}(\mathbf{x}) d\Gamma \\
    &- \int_{\Gamma_{s_{in}}}  \lambda(\mathbf{x}) \delta g(\mathbf{x}) d\Gamma  - \int_{\Gamma_{sf}}  h \lambda(\mathbf{x}) \delta T[\delta g](\mathbf{x}) d\Gamma.
\end{aligned}
\label{eq:benchmark3DhcModelSteady_adjointDerivationEnd}
\end{equation}

We can now state the adjoint problem as
\begin{problem}{(\textbf{Adjoint})}
    Find $\lambda[g]$ such that
    \begin{equation}
       k \Delta \lambda[g] +  \sum_{i=1}^M  (T[g](\mathbf{x}) - \hat{T}(\mathbf{x}))\delta(\mathbf{x}-\mathbf{x}_i)= 0, \quad      \text{in }\Omega,
    \label{eq:benchmark3DhcModelSteady_adjoint}
    \end{equation}
    with BCs
    \begin{numcases}{}
      k \nabla \lambda[g] \cdot \mathbf{n} = 0		&       on $\Gamma_{s_{in}}\cup \Gamma_{s_{ex}}$,\\
      k \nabla \lambda[g] \cdot \mathbf{n} +  h\lambda[g] = 0     &       on $\Gamma_{sf}$,
    \end{numcases}
    \label{prob:benchmark3DhcModelSteady_adjoint}
\end{problem}
$\delta(\mathbf{x}-\mathbf{x}_i)$ being the Dirac function centered at $\mathbf{x}_i$.

We notice that if $\lambda[g]$ is solution of Problem~\ref{prob:benchmark3DhcModelSteady_adjoint}, $-\lambda[g]$ represents the Fr{\'e}chet derivative of the Lagrange function with respect to the inner product in $L^2(\Gamma_{s_{in}})$.
Then, we have
\begin{equation}
    d\mathcal{L}[\delta g, \lambda[g]] =- \int_{\Gamma_{s_{in}}}  \lambda[g](\mathbf{x}) \delta g(\mathbf{x}) d\Gamma = \langle -\lambda[g], \delta g \rangle_{L^2(\Gamma_{s_{in}})}.
\end{equation}
Considering that $\mathcal{L}[g, \lambda[g]] = J_1[g]$,  the  G{\^a}teaux derivative of the functional $J_1[g]$ is
\begin{equation}
    J'_{1}[g] = - \lambda[g] \text{ in } L^2(\Gamma_{s_{in}}).
  \label{eq:GateauxDerivativeOfJ1}
\end{equation}

Different methods can be used for the solution of this minimization problem.
Here, we discuss its solution by Alifanov's regularization method\cite{Alifanov1988} and by parameterization of the heat flux, $g$.

\subsubsection{Alifanov's Regularization}\label{sec:inverse_3DhcModelSteady_alifanov}
The Alifanov's regularization method is a CGM applied on the adjoint equation.\cite{Neto2012}

We consider the following iterative procedure for the estimation of the function $g$ that minimizes the functional (\ref{eq:inverseProblem_costFunction}). 
Given an initial estimation $g^0 \in L^2(\Gamma_{s_{in}})$, for $n>0$ a new iterant is computed as: 
\begin{equation}
	g^{n+1} = g^n - \beta^n P^n ,\quad n=0,1,2,\dots\,
\label{eq:inverseProblem_heatFluxEstimation}
\end{equation}
where $n$ is the iteration counter, $\beta^n$ is the stepsize, also called correction factor, in the conjugate direction $P^n$ given by
\begin{equation}
	P^0 = J'_1 [g^0],\quad P^{n+1} = J'_1 [g^{n+1}] + \gamma^{n+1} P^n \text{ for } n \geq 1,
\label{eq:inverseProblem_searchDirection}
\end{equation}
$\gamma^{n+1}$ being the conjugate coefficient, and $J'_1 [g]$ the G{\^a}teaux derivative of $J_1$  given by (\ref{eq:GateauxDerivativeOfJ1}).

The stepsize $\beta^n$ in (\ref{eq:inverseProblem_heatFluxEstimation}) is obtained by minimizing the functional $J_1[g^n - \beta P^n]$ with respect to $\beta$.
Therefore, $\beta^n$ is the solution of the critical point equation of the functional $J_1$, restricted to a line passing through $g^n$ in the direction defined by $P^n$, i.e. $\beta^n$ is the critical point of  $J_1[g^n -\beta P^n]$ which then satisfies
\begin{equation}
    J_1[g^n-\beta^n P^n] = \min_\beta\left\{\frac{1}{2} \sum_{i=1}^M \left[ T[g^n-\beta P^n](\mathbf{x}_i) - \hat{T}(\mathbf{x}_i)\right]^2 \right\}.
\end{equation}
Recalling Problem~\ref{prob:benchmark3DhcModelSteady_sensitivity},
\begin{equation}
  J_1[g^n - \beta P^n] = \frac{1}{2} \sum_{i=1}^M \left[ T[g^n - \beta P^n](\mathbf{x}_i) - \hat{T}(\mathbf{x}_i) \right]^2 = \frac{1}{2} \sum_{i=1}^M \left[ (T[g^n] - \beta \delta T [P^n])(\mathbf{x}_i) - \hat{T}(\mathbf{x}_i) \right]^2.
  \label{eq:inverseProblem_searchStep1}
\end{equation}
Differentiating with respect to $\beta$, we obtain the critical point equation
\begin{equation}
  \frac{d J_1[g^n - \beta^n P^n]}{d\beta} = \sum_{i=1}^M [(T[g^n] - \beta^n \delta T [P^n])(\mathbf{x}_i) - \hat{T}(\mathbf{x}_i)] (-\delta T [P^n](\mathbf{x}_i)) = 0.
\end{equation}
Finally, we have
\begin{equation}
  	\beta^n=\frac{\sum_{i=1}^M \left[ T[g^n](\mathbf{x}_i) - \hat{T}(\mathbf{x}_i) \right] \delta T[P^n](\mathbf{x}_i)}{\sum_{i=1}^M (\delta T[P^n](\mathbf{x}_i))^2}.
\label{eq:inverseProblem_searchStep}
\end{equation}

With respect to the conjugate coefficient, $\gamma$, its value is zero for the first iteration and for other iterations it can be calculated using Fletcher-Reeves expression as follows\cite{Fletcher1964}
\begin{equation}
	\gamma^n = \frac{\left\lVert J'_1 [g^n]\right\rVert_{L^2(\Gamma_{s_{in}})}^2} {\left\lVert {J'_1 [g^{n-1}]}\right\rVert_{L^2(\Gamma_{s_{in}})}^2}. 
\label{eq:inverseProblem_conjugateCoeff}
\end{equation}
Notice that, to use this iterative procedure, we have to compute at each iteration the G{\^a}teaux derivative $J'_{1}[g](\mathbf{x})$ which is given by (\ref{eq:GateauxDerivativeOfJ1}).
Thus, we must solve the adjoint problem to compute it.

Alifanov's regularization algorithm is summarized in Algorithm~\ref{alg:AlifanovReg}.

\begin{algorithm}[htb]
 \caption{Alifanov's regularization.}
 \label{alg:AlifanovReg}
    \begin{algorithmic}
        \State Set $g^0$ and $n = 0$
        \While{$n < n_{max}$}
            \State Compute $T[g^n]$ by solving Problem \ref{prob:3DhcModelSteady}
            \State Compute $J_1[g^n]$ by (\ref{eq:inverseProblem_costFunction})
            \If{$J_1[g^n] < J_{1_{tol}}$}
                \State Stop
            \EndIf
            \State Compute $\lambda[g^n]$ by solving Problem \ref{prob:benchmark3DhcModelSteady_adjoint}
            \State Compute $J'_1 [g^n]$ by (\ref{eq:GateauxDerivativeOfJ1})
            \If{$n \geq 1$}
                \State Compute the conjugate coefficient, $\gamma^n$, by (\ref{eq:inverseProblem_conjugateCoeff})
                \State Compute the search direction, $P^n$, by (\ref{eq:inverseProblem_searchDirection})
            \Else
                \State $P^0 = J'_1 [g^0]$
            \EndIf
            \State Compute $\delta T[P^n]$ by solving Problem~\ref{prob:benchmark3DhcModelSteady_sensitivity} with $\delta g =  P^n$ 
            \State Compute the stepsize in the search direction, $\beta^n$, by (\ref{eq:inverseProblem_searchStep})
            \State Update heat flux $g^{n}$ by (\ref{eq:inverseProblem_heatFluxEstimation})
            \State $n = n + 1$
        \EndWhile
        \State \textbf{return} $g^{n}$
    \end{algorithmic}
\end{algorithm}

\subsubsection{Parameterization of the Boundary Conditions}\label{sec:parameterizedBC_method}

In this section, we consider the parameterization of the boundary heat flux $g$.
In the literature, the parameterization of $g$ has already been proposed.\cite{Ranut2012}
However, we propose a novel approach both for the parameterization and for the solution of the resulting inverse problem.

For the parameterization, we start considering that we want to parameterize an unknown function in $L^2(\Gamma_{s_{in}})$.
We then notice that in thin slab casting molds, the thermocouples are all located few millimeters inward from $\Gamma_{s_{in}}$.
All together they form a uniform 2D grid.
Then, to parameterize $g$, we use Radial Basis Functions (RBFs) centered at the projections of the thermocouples' points on $\Gamma_{s_{in}}$.\cite{Buhmann2003} 
Due to this choice we have as many basis functions as thermocouples.
In particular, we use Gaussian RBFs.
However, the following discussion can be applied to other basis functions.

The parameterization of the boundary heat flux reads (see Prando's appendix\cite{Prando2016})
\begin{equation}
  g(\mathbf{x}) \approx \sum_{j=1}^M w_j \phi_j(\mathbf{x}),
  \label{eq:parametrizedHeatFlux}
\end{equation}
where the $\phi_j(\mathbf{x})$ are $M$ known base functions, and the $w_j$ are the respective unknown weights.
Notice that by doing the parameterization, we change the problem from estimating a function in an infinite dimensional space to estimating a vector $\mathbf{w} =( w_1, w_2 ,\dots, w_M)^T$ in $\R^M$.

Let $\pmb{\xi}_i, 1\leq i\leq M,$ be the projection of the measurement point $\mathbf{x}_i \in \Psi$ on $\Gamma_{s_{in}}$ such that 
\begin{equation}
  \pmb{\xi}_i = \argmin_{\pmb{\xi}\in\Gamma_{s_{in}}} \left\lVert \mathbf{x}_i - \pmb{\xi}  \right\rVert_2, \quad \mathbf{x}_i\in \Psi.
\end{equation}
By centering the RBFs in these points, their expressions are
\begin{equation}
    \phi_j(\mathbf{x}) = e^{- \left( \eta \left\lVert \mathbf{x} - \pmb{\xi}_j  \right\rVert_2 \right)^2},\quad \text{ for } j=1,2,\dots,M, 
    \label{eq:GaussianRBF}
\end{equation}
where $\eta$ is the shape parameter of the Gaussian basis, increasing its values the radial decay of the basis slows down.

Suppose to have the solutions of Problem~\ref{prob:3DhcModelSteady}, $T[\phi_j]$, for $j=1,2,\dots,M$.
Denote by $T_{ad}$ the solution of 
\begin{problem}{}
Find $T_{ad}$ such that
\begin{equation}
   -k\Delta T_{ad} = 0, \quad      \text{in }\Omega,
\end{equation}
with BCs
\begin{numcases}{}
    -k \nabla T_{ad} \cdot \mathbf{n} = 0              &       on $\Gamma_{s_{in}} \cup \Gamma_{s_{ex}}$,\\
    -k \nabla T_{ad} \cdot \mathbf{n} = h(T_{ad} +  T_f) 	&       on $\Gamma_{sf}$.
\end{numcases}
\label{prob:additive}
\end{problem}
Then, we see that 
\begin{equation}
    T[\mathbf{w}] = \sum_{j=1}^M w_j (T[\phi_j] + T_{ad}) - T_{ad},
  \label{eq:directProbDecomposition}
\end{equation}
is solution of Problem~\ref{prob:3DhcModelSteady} since
\begin{equation}
   - k \Delta (\sum_{j=1}^M w_j (T[\phi_j] + T_{ad}) - T_{ad})= 0, \quad      \text{in }\Omega,
\end{equation}
and verifies the BCs associated to that problem
\begin{numcases}{}
  -k \nabla (\sum_{j=1}^M w_j (T[\phi_j] + T_{ad}) - T_{ad}) \cdot \mathbf{n} = \sum_{j=1}^M w_j \phi_j=g	& on $\Gamma_{s_{in}}$,\\
  -k \nabla (\sum_{j=1}^M w_j (T[\phi_j] + T_{ad}) - T_{ad}) \cdot \mathbf{n} = 0				& on $\Gamma_{s_{ex}}$,\\
  -k \nabla (\sum_{j=1}^M w_j (T[\phi_j] + T_{ad}) - T_{ad}) \cdot \mathbf{n} =  h \Bigg(\sum_{j=1}^M w_j (T[\phi_j] + T_{ad}) - T_{ad} - T_f \Bigg)    &  on $\Gamma_{sf}$.
\end{numcases}

Now, the objective of the inverse problem is to determine $\mathbf{w}$ which identifies $g$ once the elements of the base $\phi_j$, $j=1,2,\dots,M$ are fixed.
Notice that we consider all vectors as column vectors.

We rewrite the inverse Problem \ref{inverseProblem_3DhcModelSteady} as
\begin{problem}{(\textbf{Inverse})}
  Given the temperature measurements $\hat{T}(\Psi)\in \R^M$, find $\mathbf{w}\in\R^M$ which minimizes the functional
  \begin{equation}
    J_1[\mathbf{w}]=\frac{1}{2} \sum_{i=1}^M [T[\mathbf{w}](\mathbf{x}_i)-\hat{T}(\mathbf{x}_i)]^2,
  \end{equation}
\label{prob:benchmark3DhcModelSteady_inverseProblem_parametrizedBC}
\end{problem}
where to simplify notation, and if there is no room for error, $T[\mathbf{w}]$ represents the solution $T[g]$ of Problem \ref{prob:3DhcModelSteady} with $g$ as in (\ref{eq:parametrizedHeatFlux}).  

Given $\mathbf{w}$, we define the residual $\mathbf{R}[\mathbf{w}] \in \R^M$ as the vector whose components are
\begin{equation}
  (\mathbf{R}[\mathbf{w}])_i := (\mathbf{T}[\mathbf{w}])_i - (\hat{\mathbf{T}})_i,
  \label{eq:benchmark3DhcModelSteady_residual}
\end{equation}
where $\mathbf{T}[\mathbf{w}]$ and $\hat{\mathbf{T}}$ denote the vectors of $\R^M$ whose i-component is  $(\mathbf{T}[\mathbf{w}])_i = T[\mathbf{w}](\mathbf{x}_i)$ and $(\hat{\mathbf{T}})_i = \hat{T}(\mathbf{x}_i)$, respectively.
So, we rewrite the cost functional
\begin{equation}
  J_1[\mathbf{w}]=\frac{1}{2} \mathbf{R}[\mathbf{w}]^T \mathbf{R}[\mathbf{w}].
\end{equation}
To minimize the functional $J_1[\mathbf{w}]$, we solve the critical point equation
\begin{equation}
  \frac{\partial J_1[\mathbf{w}]}{\partial w_j} = \sum_{i=1}^M R[\mathbf{w}]_i \frac{\partial (T[\mathbf{w}])_i}{\partial w_j} = 0,\text{ for } j=1,2,\dots,M.
  \label{eq:benchmark3DhcModelSteady_criticalPoint}
\end{equation}
Thanks to (\ref{eq:directProbDecomposition}), equation (\ref{eq:benchmark3DhcModelSteady_criticalPoint}) can be written as
\begin{equation}
    \mathbf{R}[\mathbf{w}]^T (\mathbf{T}[\phi_j] +  \mathbf{T}_{ad}) = 0,\text{ for } j=1,2,\dots,M,
  \label{eq:benchmark3DhcModelSteady_criticalPoint2}
\end{equation}
being $\mathbf{T}_{ad}$ the vector of $\R^M$ whose i-component is $T_{ad}(\mathbf{x}_i)$.
Then, the vector associated to the solution of the direct problem in the measurement points, $\mathbf{T}[\mathbf{w}] \in \R^M$,  is given by
\begin{equation}
  \mathbf{T}[\mathbf{w}] = \sum_{j=1}^M w_j \mathbf{T}[\phi_j] + (\sum_{j=1}^M w_j - 1) \mathbf{T}_{ad}.
  \label{eq:solAtMeasurements}
\end{equation}

We denote by $\Theta$ the matrix of $\R^{M \times M}$ such that
\begin{equation}
  \Theta_{ij} := T[\phi_j](\mathbf{x}_i) + T_{ad}(\mathbf{x}_i).
\end{equation}
Therefore, equation (\ref{eq:benchmark3DhcModelSteady_criticalPoint2}) can now be written as
\begin{equation}
    \Theta^T \mathbf{R}[\mathbf{w}] = \mathbf{0}.
\end{equation}

Recalling the definition of $\mathbf{R}$ and (\ref{eq:solAtMeasurements}), we have
\begin{equation}
  \Theta^T \mathbf{R}[\mathbf{w}] = \Theta^T(\Theta \mathbf{w} - \mathbf{T}_{ad} - \hat{\mathbf{T}}) = \mathbf{0}.
\end{equation}
The solution of the inverse problem is then obtained by solving the linear system
\begin{equation}
  \Theta^T \Theta \mathbf{w} = \Theta^T(\hat{\mathbf{T}} + \mathbf{T}_{ad}).
  \label{eq:linSys_parametrizedBC}
\end{equation}
This is generally called the normal equation.

In this setting, (\ref{eq:linSys_parametrizedBC}) is a linear map from the observations to the heat flux weights.
Consequently, we have that the existence and uniqueness of the solution of the inverse problem depends on the invertibility of the matrix $\Theta^T \Theta$.

We can easily see that the matrix $\Theta^T \Theta $ is symmetric and positive semi-definite.
In general, however, we cannot ensure that it is invertible.
In fact, the invertibility depends on the choice of the basis function, the computational domain and the BCs.

In the numerical tests, we will see that this matrix tends to be ill-conditioned.
This is a reflect of the ill-posedness of the inverse problem.
Different regularization techniques for linear systems are available to overcome this issue.\cite{Bardsley2018}
Here, we consider the Truncated Singular Value Decomposition (TSVD) regularization.
We denote the Singular Value Decomposition (SVD) of $\Theta^T \Theta$ by
\begin{equation}
    \Theta^T \Theta = U \Sigma V^T = \sum_{i=1}^r \mathbf{u}_i \sigma_i \mathbf{v}_i^T,
\end{equation}
where $\sigma_i$ denotes the i-th singular value of $\Theta^T \Theta$ (numbered according to their decreasing value), $r$ denotes the first no null singular value, i.e. the rank of $\Theta^T \Theta$, $\mathbf{u}_i$ and $\mathbf{v}_i$ are the i-th columns of the semi-unitary matrices $U$ and $V$ respectively (both belonging to $\R^{M \times  r}$), and $\Sigma$ is the square diagonal matrix of $\R^{r \times  r}$ such that $\Sigma_{ii} = \sigma_i$ and $\Sigma_{ij} = 0$ if $i\neq j$.
Then, the TSVD regularized solution of (\ref{eq:linSys_parametrizedBC}) is
\begin{equation}
    \mathbf{w} = \sum_{i=1}^{\alpha_{TSVD}} \left(\frac{\mathbf{u}_i^T \Theta^T (\hat{\mathbf{T}} + \mathbf{T}_{ad}) }{\sigma_i}\right) \mathbf{v}_i.
\end{equation}
This solution differs from the least square solution only in that the sum is truncated at $i = \alpha_{TSVD}$ instead of $i=r$.

We conclude our discussion of this method by noticing its most interesting feature for our investigation.
In fact, it is already suitable for real-time computation since we can divide it into an offline (expensive) phase and an online (cheap) phase.
In the offline phase, we compute $T[\phi_j]$ for $j=1,2,\dots,M$ and $T_{ad}$ by solving Problem~\ref{prob:3DhcModelSteady} with each base as boundary heat flux and Problem~\ref{prob:additive}.
Then in the online phase, we input the measurements $\hat{\mathbf{T}}$ and solve the linear system (\ref{eq:linSys_parametrizedBC}).
For the choice made when selecting the basis functions, the linear system has the dimensions of the number of thermocouples.
Then, its solution can be easily done in real-time even with limited computational power.
This makes this method very promising for our real-time application.

\subsection{Inverse Problem with Thermocouples and Total Heat Measurement}\label{sec:inverseProblem_totalHeat}

In CC molds, we can have together with thermocouples' pointwise measurements also total heat flux measurements.
Assuming all boundaries but $\Gamma_{s_{in}}$  and  $\Gamma_{sf}$ to be adiabatic, all heat is extracted by the mold by the cooling water at $\Gamma_{sf}$.
Further, assuming the water heat capacity, $C_{p_f}$, to be constant and the the water mass flow rate $\dot{m}$ to be known, the total heat flux is given by
\begin{equation}
    \hat{G} = \int_{\Gamma_{s_{in}}} g d\Gamma = \dot{m} C_{p_f} (T_{f_{out}} - T_{f_{in}}),
    \label{eq:totalHeatMeasurement}
\end{equation}
where $T_{f_{in}}$ and $T_{f_{out}}$ are the cooling water temperatures at the inlet and outlet of the cooling system, respectively.
Then, the total heat flux measurements is obtained by (\ref{eq:totalHeatMeasurement}) after experimentally  measuring $T_{f_{in}}$, $T_{f_{out}}$ and the  water mass flow rate $\dot{m}$.

In this section, we discuss the formulation and solution of the inverse problem of estimating the boundary heat flux, $g$, by considering both the thermocouples' and total heat flux measurements.

Using again a least square, deterministic approach, we state the inverse problem as
\begin{problem}{(\textbf{Inverse})}
  Given $\{ \hat{T}(\mathbf{x}_i) \}_{i=1}^M$ and $\hat{G}$, find the heat flux $g \in G_{ad}$ that minimizes the functional $J_2:L^2(\Gamma_{s_{in}}) \rightarrow \R^+$,
  \begin{equation}
      J_2[g]:=\frac{1}{2}\sum^M_{i=1} [T[g](\mathbf{x}_i) - \hat{T}(\mathbf{x}_i)]^2 + \frac{1}{2} p_g \bigg( \int_{\Gamma_{s_{in}}} g d\Gamma - \hat{G} \bigg)^2,
    \label{eq:inverseProblemTotalHeat_costFunction}
  \end{equation}
    where $T[g](\mathbf{x}_i)$ is the solution of Problem~\ref{prob:3DhcModelSteady} at points $\mathbf{x}_i$, for all $i=1,2,\dots,M$, and $p_g[\frac{K^2}{W^2}]$ is a weight applied to the total heat measurement.
  \label{inverseProblemTotalHeat_3DhcModelSteady}
\end{problem}
Notice that, thanks to Theorem~\ref{theo:Nittka2011_3.14} the state variable $T$ is continuous in $\Omega_s$, then its value at pointwise observations is well-defined.

To derive the adjoint of Problem~\ref{inverseProblemTotalHeat_3DhcModelSteady}, we redo computations (\ref{eq:benchmark3DhcModelSteady_adjointDerivation1})-(\ref{eq:benchmark3DhcModelSteady_adjointDerivationEnd}).
It turns out that the adjoint of Problem~\ref{inverseProblemTotalHeat_3DhcModelSteady} is again Problem~\ref{prob:benchmark3DhcModelSteady_adjoint}.
However, the Fr{\'e}chet derivative with respect to the inner product in $L^2(\Gamma_{s_{in}})$ of $J_2$ is
\begin{equation}
        d\mathcal{L}_g[\delta g, \lambda] =- \int_{\Gamma_{s_{in}}}  \bigg[ \lambda[g](\mathbf{x}) - p_g \bigg( \int_{\Gamma_{s_{in}}} g d\Gamma - \hat{G} \bigg) \bigg]  \delta g(\mathbf{x}) d\Gamma = \langle -\lambda[g] + p_g \bigg( \int_{\Gamma_{s_{in}}} g d\Gamma - \hat{G} \bigg) , \delta g \rangle_{L^2(\Gamma_{s_{in}})}.
\end{equation}
Considering that $\mathcal{L}[g, \lambda[g]] = J_2[g]$,  the  G{\^a}teaux derivative of the functional $J_2[g]$ is
\begin{equation}
    J'_{2}[g] = -\lambda[g] + p_g \bigg( \int_{\Gamma_{s_{in}}} g d\Gamma - \hat{G} \bigg)  \text{ in } L^2(\Gamma_{s_{in}}).
  \label{eq:GateauxDerivativeOfJ2}
\end{equation}

Different methods can be used for the solution of this minimization problem.
As for the minimization of $J_1$, we discuss its solution by Alifanov's regularization method and by parameterization of the heat flux, $g$.

\subsubsection{Alifanov's Regularization}\label{sec:inverse_3DhcModelSteady_totalHeatFlux_alifanov}
In this section, we expand the discussion in Section~\ref{sec:inverse_3DhcModelSteady_alifanov} on Alifanov's regularization to the inverse Problem~\ref{inverseProblemTotalHeat_3DhcModelSteady}.

We consider the following iterative procedure for the estimation of the function $g$ that minimizes functional (\ref{eq:inverseProblemTotalHeat_costFunction}). 
Given an initial estimation $g^0 \in L^2(\Gamma_{s_{in}})$, for $n>0$ a new iterant is computed by (\ref{eq:inverseProblem_heatFluxEstimation}) with the conjugate direction given by
\begin{equation}
	P^0 = J'_2 [g^0],\quad P^{n+1} = J'_2 [g^{n+1}] + \gamma^{n+1} P^n \text{ for } n \geq 1,
\label{eq:inverseProblemTotalHeat_searchDirection}
\end{equation}
and the search step computed by 
\begin{equation}
    \beta^n=\frac{\sum_{i=1}^M \left[ T[g^n](\mathbf{x}_i) - \hat{T}(\mathbf{x}_i)\right] \delta T[P^n](\mathbf{x}_i) + p_g \left(  \int_{\Gamma_{s_{in}}} P^n  d\Gamma \right) \left(\int_{\Gamma_{s_{in}}}  g^n d\Gamma - \hat{G} \right) }{\sum_{i=1}^M (\delta T[P^n](\mathbf{x}_i))^2 + p_g \left( \int_{\Gamma_{s_{in}}} P^n  d\Gamma \right)^2}.
\label{eq:inverseProblemTotalHeat_searchStep}
\end{equation}

Alifanov's regularization algorithm is then as in Algorithm~\ref{alg:AlifanovReg} where the functional $J_1$ is substituted by $J_2$ and the search step and conjugate direction are computed by (\ref{eq:inverseProblemTotalHeat_searchStep}) and (\ref{eq:inverseProblemTotalHeat_searchDirection}), respectively.

\subsubsection{Parameterization of the Boundary Conditions}\label{sec:inverse_3DhcModelSteady_totalHeatFlux_parameterizedBC_method}
In this section, we apply the discussion made in Section~\ref{sec:parameterizedBC_method} to Problem~\ref{inverseProblemTotalHeat_3DhcModelSteady}.

Considering the parameterization (\ref{eq:parametrizedHeatFlux}), due to (\ref{eq:benchmark3DhcModelSteady_residual}), we rewrite (\ref{eq:inverseProblemTotalHeat_costFunction}) as
\begin{equation}
        J_2[\mathbf{w}] = \frac{1}{2} \mathbf{R}[\mathbf{w}]^T  \mathbf{R}[\mathbf{w}] + \frac{1}{2} p_g \bigg( \int_{\Gamma_{s_{in}}} \sum_{j=1}^M w_j \phi_j(\mathbf{x}) d\Gamma - \hat{G} \bigg)^2 = \frac{1}{2} \mathbf{R}[\mathbf{w}]^T  \mathbf{R}[\mathbf{w}] + \frac{1}{2} p_g \bigg( \sum_{j=1}^M w_j \int_{\Gamma_{s_{in}}} \phi_j(\mathbf{x}) d\Gamma - \hat{G} \bigg)^2.
\end{equation}
Defining the vector in $\R^M$ such that
\begin{equation}
    (\pmb{\phi})_i := \int_{\Gamma_{s_{in}}} \phi_i(\mathbf{x}) d\Gamma,
\end{equation}
we write
\begin{equation}
    J_2[\mathbf{w}] = \frac{1}{2} \mathbf{R}[\mathbf{w}]^T  \mathbf{R}[\mathbf{w}] + \frac{1}{2} p_g (\mathbf{w}^T \pmb{\phi} - \hat{G})^2.  
\end{equation}

As in Section~\ref{sec:parameterizedBC_method}, we now write the critical point equation for $J_2[\mathbf{w}]$
\begin{equation}
    \frac{\partial J_2[\mathbf{w}]}{\partial w_j} = \sum^M_{i=1} (R[\mathbf{w}])_i \frac{\partial (T[\mathbf{w}])_i}{\partial w_j} + p_g (\mathbf{w}^T \pmb{\phi} - \hat{G})(\pmb{\phi})_j = 0,\quad \text{ for } j=1,2,\dots,M.
\end{equation}
Then, introducing the matrix in $\R^{M \times M}$ such that
\begin{equation}
    \Phi_{ij} = (\pmb{\phi})_i (\pmb{\phi})_j,
\end{equation}
we can write the critical point equation as
\begin{equation}
    (\Theta^T \Theta + p_g \Phi) \mathbf{w} = p_g \hat{G} \pmb{\phi} + \Theta^T(\mathbf{T}_{ad} + \hat{\mathbf{T}}).
    \label{eq:inverse_3DhcModelSteady_totalHeatFlux_paramBClinSys}
\end{equation}
By solving the linear system (\ref{eq:inverse_3DhcModelSteady_totalHeatFlux_paramBClinSys}) we obtain the weights $\mathbf{w}$ of the parameterization.
Then, by (\ref{eq:parametrizedHeatFlux}) we compute the estimated heat flux $g$.
Also in this setting, the discussion at the end of Section~\ref{sec:parameterizedBC_method} on the regularization of the linear system and the offline-online decomposition holds.

\section{Analytical Benchmark}\label{sec:analyticalBenchmark}

In this section, we propose an academic benchmark case.
It is a steady-state heat conduction problem in a homogeneous isotropic solid occupying a rectangular parallelepiped domain.
By carefully selecting the BCs on the faces of the parallelepiped, we are able to compute the analytical solution of the heat conduction problem.
Then, we use this academic test to validate the numerical solution of the direct problem.
Moreover, by arbitrarily selecting some temperature measurements points, we test the different inverse problem solution methodologies discussed in Section~\ref{sec:inverseProblem}.

\begin{figure}[!htb]
\centering
\includegraphics[width=0.5\textwidth]{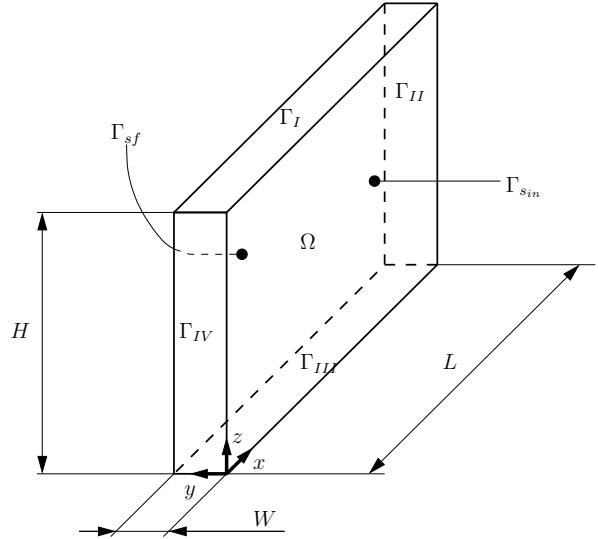}
    \caption{Schematic of the solid rectangular parallelepiped domain.}
\label{fig:analyticalBenchmarkDomainSchematic}
\end{figure}

Let the domain be  $\Omega = (0,L) \times (0,W) \times (0,H)$ as in Figure~\ref{fig:analyticalBenchmarkDomainSchematic} with positive real constants $L, W$ and $H$.
Let $\Gamma$ be boundary of $\Omega$.
Then, the different boundaries of the domain to be considered are
\begin{equation}
\begin{aligned}
    &\Gamma_{sf} :=\{\mathbf{x}\in\Gamma |\ \mathbf{x}=(x,W,z)\}, \quad && \Gamma_{s_{in}} :=\{\mathbf{x}\in\Gamma |\ \mathbf{x}=(x,0,z)\},\\
    &\Gamma_{I}  :=\{\mathbf{x}\in\Gamma |\ \mathbf{x}=(x,y,H)\}, \quad && \Gamma_{III}:=\{\mathbf{x}\in\Gamma |\ \mathbf{x}=(x,y,0)\},\\
    &\Gamma_{II} :=\{\mathbf{x}\in\Gamma |\ \mathbf{x}=(L,y,z)\}, \quad && \Gamma_{IV} :=\{\mathbf{x}\in\Gamma |\ \mathbf{x}=(0,y,z)\}.
\end{aligned}
\end{equation}

To have an analytical solution, $T_{an}$, in $\Omega$, we consider a slight modification of Problem~\ref{prob:3DhcModelSteady} that does not change its essential aspects.
\begin{problem}
  Find $T$ such that
  \begin{equation}
     - k \Delta T = 0, \text{ in }\Omega,
  \end{equation}
  with BCs
  \begin{numcases}{}
      -k \nabla T \cdot \mathbf{n} = g_{an}       &   on $\Gamma_{s_{in}}$,\label{eq:analyticalBenchmark_GammaInBC}\\
      -k \nabla T \cdot \mathbf{n} = q_{L}        &   on $\Gamma_{L}, L \in \{I, II, II ,IV\}$,\\
      -k \nabla T \cdot \mathbf{n} = h(T - T_f) &   on $\Gamma_{sf}$.
  \end{numcases}
  \label{prob:analyticalBenchmark}
\end{problem}

Let $a$, $b$, $c$ be real constants.
To have an analytical solution in $\Omega$, we consider the following data as BCs for Problem~\ref{prob:analyticalBenchmark},
\begin{equation}
  \begin{aligned}
      &q_{I}  (\mathbf{x})     = 2 k a H,           \quad && q_{III}(\mathbf{x}) = 0,\\
      &q_{II} (\mathbf{x})     = -k ( 2 a L + b y), \quad && q_{IV} (\mathbf{x}) = k b y, \\
      &T_f    (\mathbf{x})     =\frac{k (b x + c)}{h} + a x^2 + c y - a z^2 + b x W + c, 
  \end{aligned}
\end{equation}
with
\begin{equation}
    g_{an} (\mathbf{x}) = k (b x + c),
    \label{eq:analyticalHeatFlux}
\end{equation}
$k$ being the thermal considered conductivity, that is assumed constant.
Then,
\begin{equation}
    T_{an}(\mathbf{x})= a x^2+ b x y + c y - a z^2 + c,
  \label{eq:analyticalBenchmark_analyticalSolution}
\end{equation}
is the solution to Problem~\ref{prob:analyticalBenchmark}.

\subsection{Direct Problem}

We now discuss the numerical solution of Problem~\ref{prob:analyticalBenchmark}.
Due to its simplicity, the domain $\Omega$ is discretized by uniform, structured, orthogonal, hexahedral  meshes.
To study the convergence of the numerical solution to the analytical one, we consider grids with different degree of refinement.
In all tests, we use the same number of edges for the three axes.

With respect to the used finite volume scheme, since we have a structured orthogonal grid, no correction is needed when computing the gradient normal to the cells faces.
Moreover, we use linear interpolation to interpolate the values from cell centers to face centers. 
The resulting scheme is second order accurate.

From the discretization of Problem~\ref{prob:analyticalBenchmark}, we obtain a linear system.
We solve it by using the preconditioned conjugate gradient solver with diagonal incomplete Cholesky preconditioning.
The tolerance used for the linear system solver is $10^{-12}$.
All the computations are performed in ITHACA-FV\cite{Stabile2017,ithaca} which is a C++ library based on OpenFOAM\cite{Moukalled2015} developed at the SISSA Mathlab.

Finally, Table~\ref{tab:analyticalBenchmark_parameters} summarizes the parameters used for the computations. 

\begin{table}[!htb]
    \centering
    \caption{Parameters used for the simulation of the analytical benchmark.}
    \label{tab:analyticalBenchmark_parameters}
    \begin{tabular}{ ll }
        \hline
        \textbf{Parameter}    &   \textbf{Value}\\
        \hline
        Thermal conductivity, $k$ 	& $3.0~ W/(m K)$ \\
        Heat transfer coefficient, $h$ 	& $5.0~W/(m^2 K)$\\
        $a$	                        & $5 ~ K/m^2$ \\
        $b$	                        & $10~ K/m^2$ \\
        $c$ 	                        & $15~ K/m^2$ \\
        $L$ 	                        & $1~m$\\
        $W$	                        & $1~m$\\
        $H$	                        & $1~m$\\
        \hline
    \end{tabular}
\end{table}

\begin{figure}[!htb]
    \centering
    \includegraphics[width=0.6\textwidth]{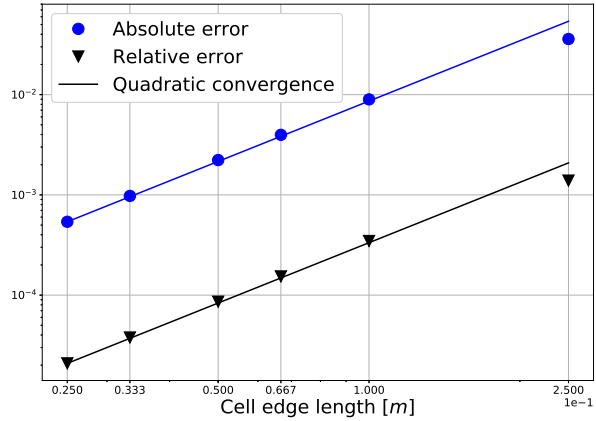}
    \caption{Decay of the absolute and relative difference in $L^2$-norm between the computed and true temperature field with the mesh refinement.}
    \label{fig:analyticalBenchmark_gridRefinement}
\end{figure}

To evaluate the accuracy of the numerical solutions, we show in Figure~\ref{fig:analyticalBenchmark_gridRefinement} the decay of the absolute and relative difference in the $L^2$-norm between the computed and true temperature field.
The test confirms the second order accuracy of the used finite volume scheme.
We conclude that Problem~\ref{prob:analyticalBenchmark} is numerically well solved.

\subsection{Inverse Problem with Temperature Measurements}

To numerically analyze the performances of the inverse solvers, we design the following test: we select a surface inside $\Omega_s$ which is parallel to $\Gamma_{s_{in}}$, and on this surface we locate $M$ measurement points  which correspond to the location of M virtual thermocouples.
The temperature in these points is given by $\hat{T}(\mathbf{x}_i) = T_{an}(\mathbf{x}_i), i=1,\dots,M$, being $T_{an}$ the solution of Problem~\ref{prob:analyticalBenchmark}, given by (\ref{eq:analyticalBenchmark_analyticalSolution}).
Using these temperatures as measurements, we apply the methods described in Section~\ref{sec:inverseProblem} to solve the inverse Problem~\ref{inverseProblem_3DhcModelSteady}, considering $T[g](\mathbf{x}_i)$ as the solution of Problem~\ref{prob:analyticalBenchmark} replacing $g_{an}$ by $g$.

\begin{figure}[!htb]
\centering
\includegraphics[width=0.4\textwidth]{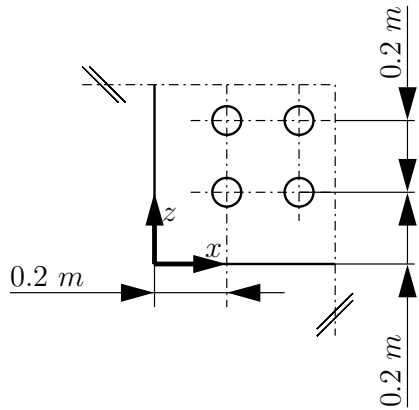}
    \caption{Positions of the virtual thermocouples for the analytical test case.}
\label{fig:analyticalBenchmark_thermocouplesPosition}
\end{figure}

The virtual thermocouples are located in the plane $y=0.2~m$.
Their $(x,z)$ coordinates are shown in Figure~\ref{fig:analyticalBenchmark_thermocouplesPosition}.
Then, we have 16 thermocouples located on the nodes of a uniform lattice at the plane $y=0.2~m$, unless otherwise stated.

The parameters used for the computations are summarized in Table~\ref{tab:analyticalBenchmark_inverseProblem_parameters}.
In this section, we test the inverse methodologies of Section~\ref{sec:inverseProblem} analyzing the effect of different parameters such as grid refinement, CG stopping criterion, RBF shape parameter, measurement noise, etc.
To analyze the numerical results, we will often use the following error norms 
\begin{equation}
  \left\lVert \varepsilon \right\rVert_{L^2(\Gamma_{s_{in}})} = \left\lVert  \frac{g - g_{an} }{g_{an}} \right\rVert_{L^2(\Gamma_{s_{in}})}, \quad 
  \left\lVert \varepsilon \right\rVert_{L^\infty(\Gamma_{s_{in}})} = \left\lVert  \frac{g - g_{an} }{g_{an}} \right\rVert_{L^\infty(\Gamma_{s_{in}})}. \quad 
    \label{eq:analyticalBenchmark_errorNorms}
\end{equation}
Notice that from (\ref{eq:analyticalHeatFlux}), $g_{an} > 0$.

\begin{table}[!htb]
\centering
\caption{Parameters used in testing the inverse problem solvers for the analytical benchmark case.}
\label{tab:analyticalBenchmark_inverseProblem_parameters}
\begin{tabular}{ ll }
    \hline
    \textbf{Parameter}    &   \textbf{Value}\\
    \hline
    N. of thermocouples	& 16 \\
    Thermocouples plane	& $y = 0.2~m$ \\
    $g^0$			& $0~W/m^2$\\
    RBF kernel		& Gaussian\\
    N. of RBF		& 16 \\
    Shape parameter, $\eta$	& 0.7 \\
    \hline
\end{tabular}
\end{table}

\subsubsection{Alifanov's Regularization}

In this section, we analyze the effect that the grid refinement and the stop criterion have on the results obtained by the Alifanov's regularization.

\begin{figure}[!htb]
    \begin{subfigure}[htb]{0.5\linewidth}
        \centering
        \includegraphics[width=\textwidth]{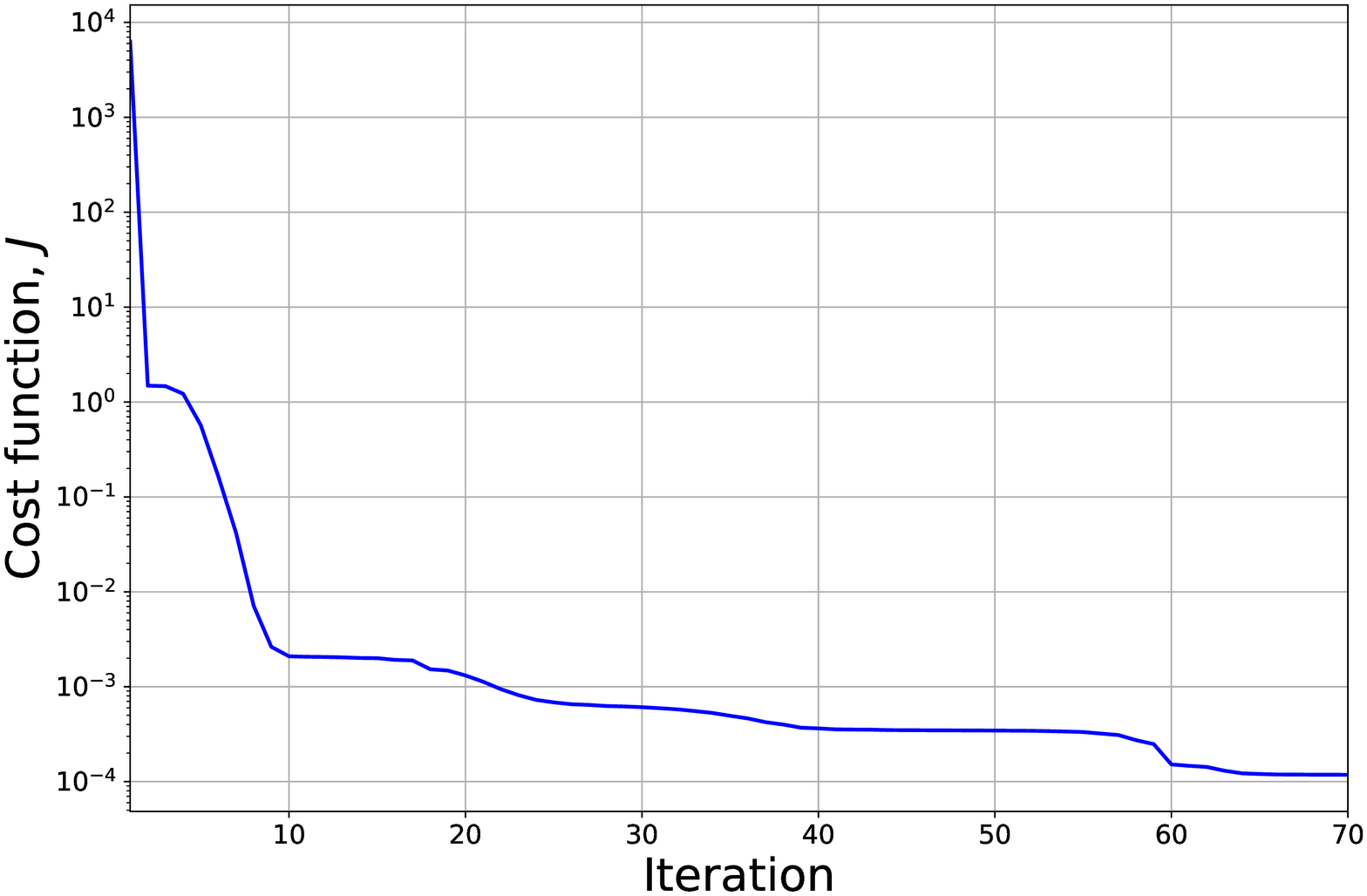}
        \caption{Cost functional, $J_1$.}
    \end{subfigure}%
    \begin{subfigure}[htb]{0.5\linewidth}
        \centering
        \includegraphics[width=0.95\textwidth]{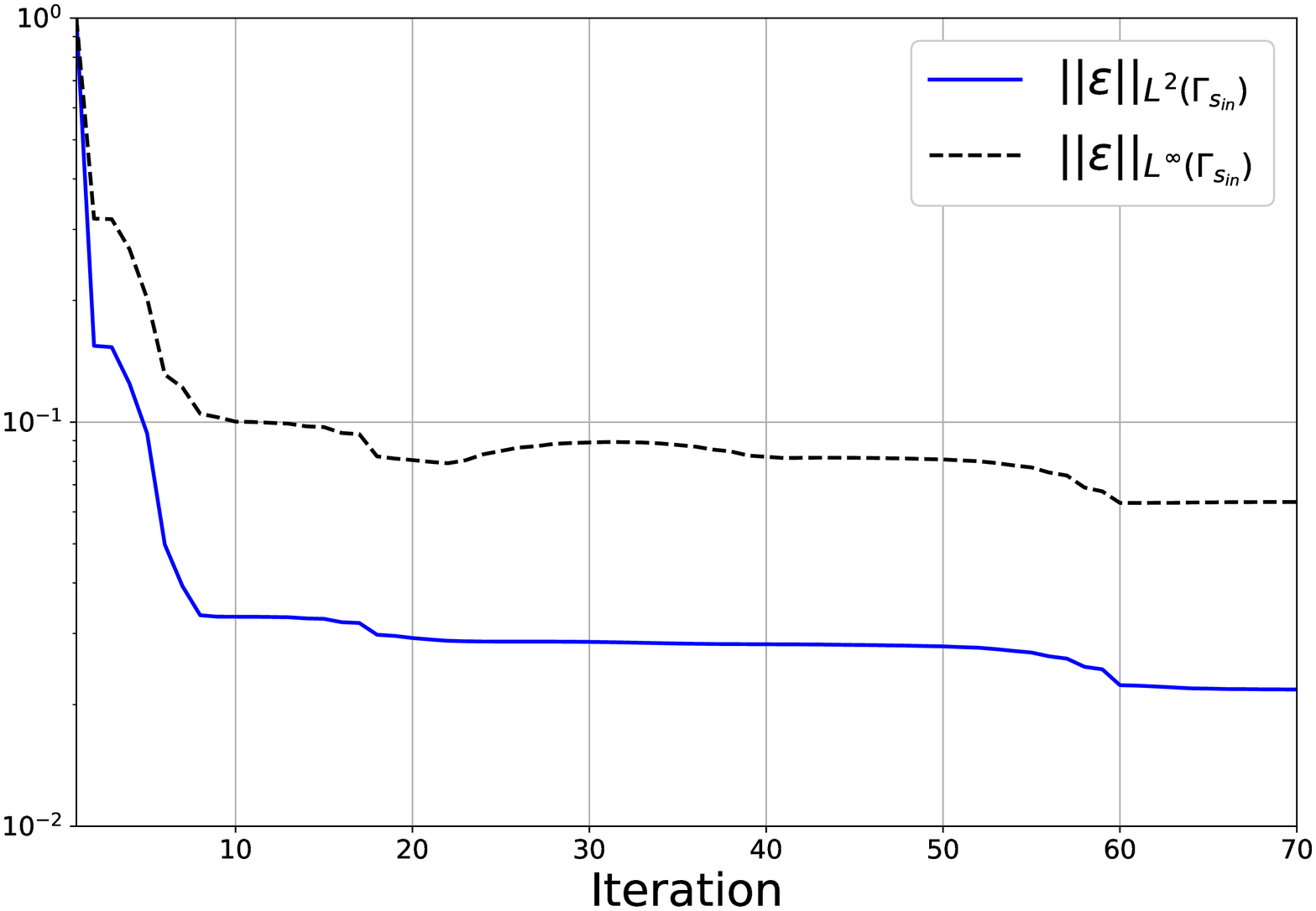}
        \caption{Relative error.}
    \end{subfigure}%
    \caption{Behavior of the cost functional $J_1$ (A) and of the heat flux  relative error $L^2$- and $L^\infty$-norms (B) with respect to the Alifanov's regularization iterations for the analytical benchmark case.}
\label{fig:analyticalBenchmark_y02_CGnoInt_convergence}
\end{figure}

We begin by comparing in Figure~\ref{fig:analyticalBenchmark_y02_CGnoInt_convergence} (a) and (b) the behavior of the functional  $J_1$ together with the $L^2$- and $L^\infty$-norm of the relative error  defined in (\ref{eq:analyticalBenchmark_errorNorms}) as functions of the number of iterations of the algorithm.
Both the cost function and the relative error have a sharp decay in the first 10 iterations.
Then, the convergence rate has a dramatic decrease reaching a plateau after 60 iterations.

\begin{figure}[!htb]
    \begin{subfigure}[htb]{\linewidth}
        \centering
        \begin{subfigure}[c]{.33\linewidth}
            \includegraphics[width=.95\textwidth]{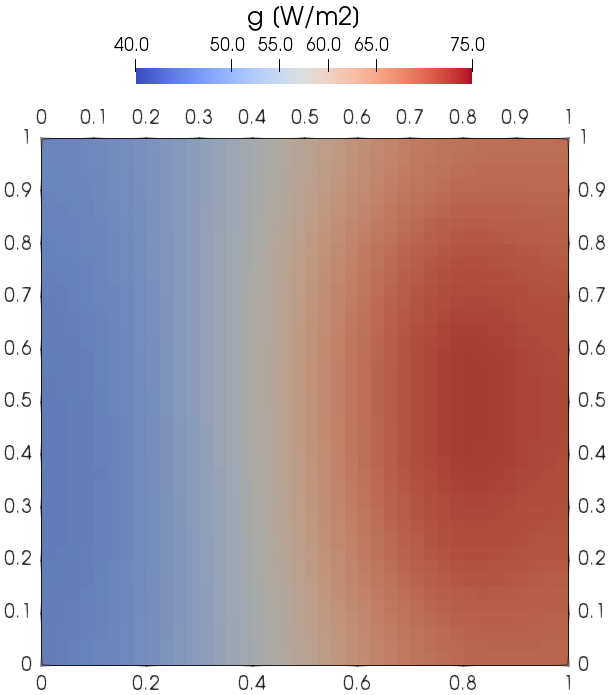}
            \caption{Iteration 10}
        \end{subfigure}%
        \begin{subfigure}[c]{.33\linewidth}
            \includegraphics[width=.95\textwidth]{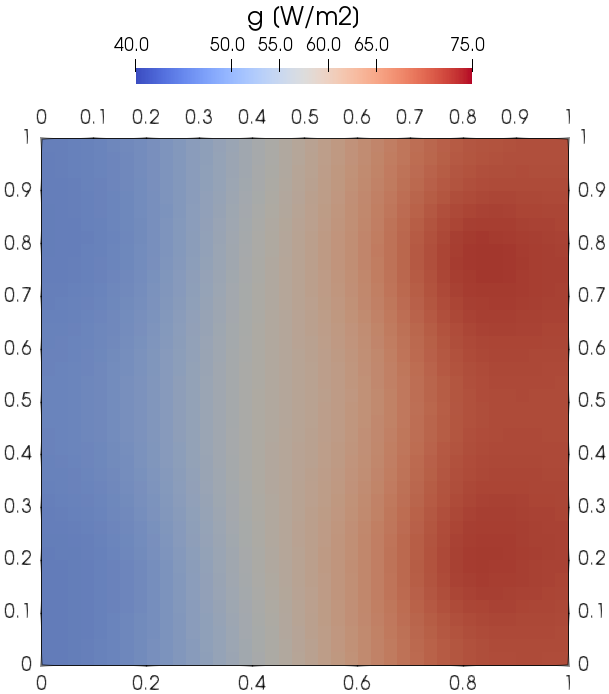}
            \caption{Iteration 70}
        \end{subfigure}%
        \begin{subfigure}[c]{.33\linewidth}
            \includegraphics[width=.95\textwidth]{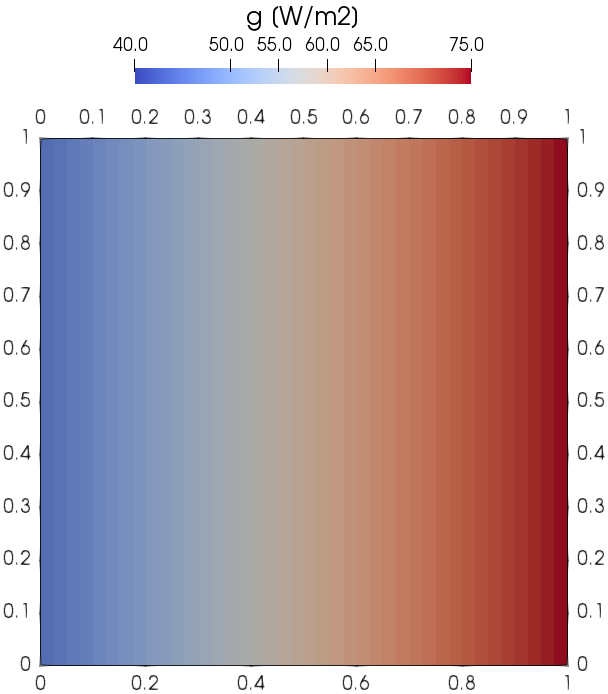}
            \caption{$g_{an}$}
        \end{subfigure}%
    \end{subfigure}%
    \caption{The estimated heat flux by Alifanov's regularization at different iterations (A,B) is compared to the analytical value (C) in the analytical benchmark case.} 
\label{fig:analyticalBenchmark_y02_CGnoInt_heatFlux}
\end{figure}

To have qualitatively insight on the results, we compare the computed heat flux at different iterations in Figure~\ref{fig:analyticalBenchmark_y02_CGnoInt_heatFlux}.
In few iterations, the estimated heat flux is already in good agreement with the analytical BC.
Then, the last iterations improve slightly the estimation.

\begin{figure}[!htb]
    \begin{subfigure}[htb]{0.5\linewidth}
        \centering
        \includegraphics[width=\textwidth]{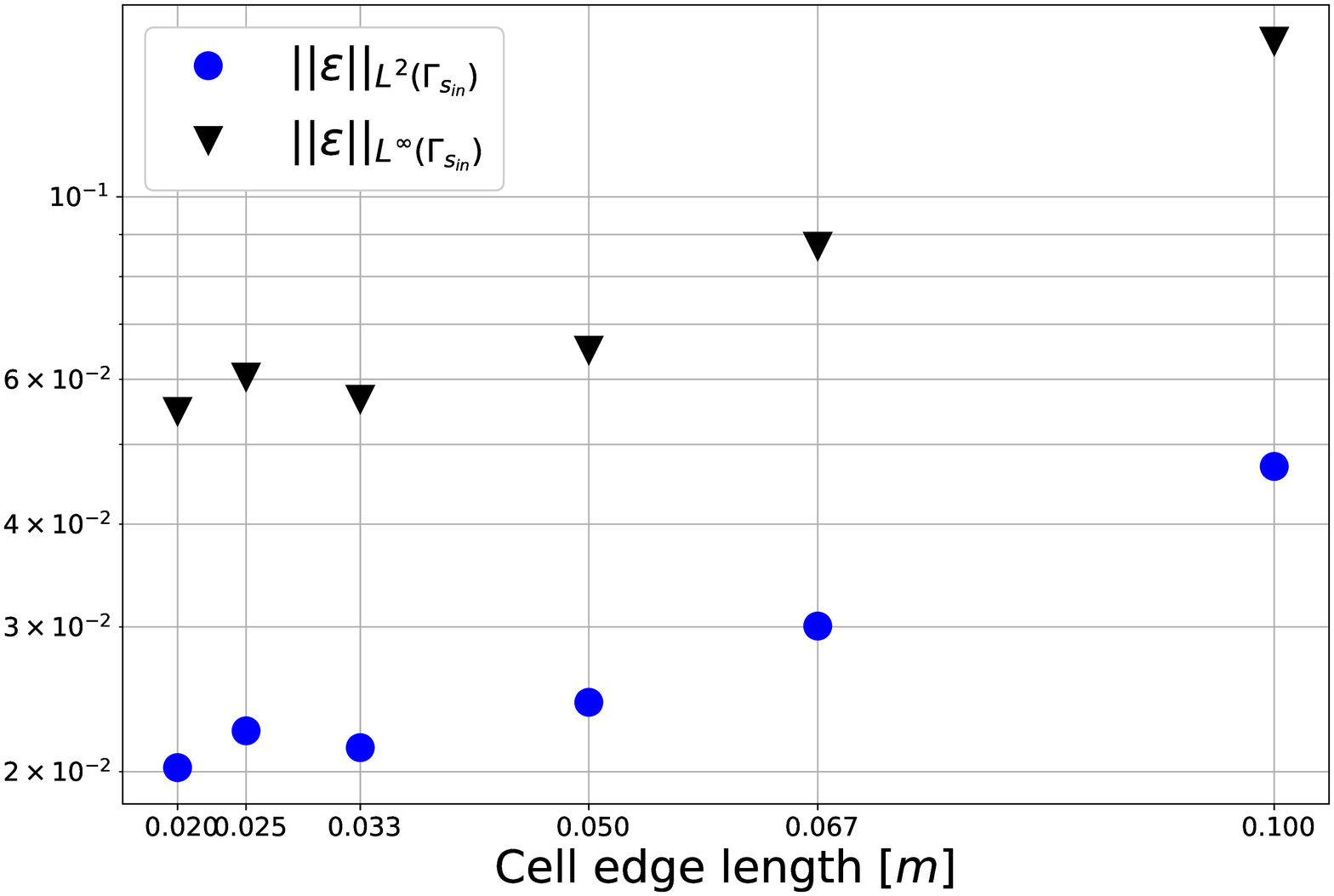}
        \caption{Mesh refinement.}
    \end{subfigure}%
    \begin{subfigure}[htb]{0.5\linewidth}
        \centering
        \includegraphics[width=\textwidth]{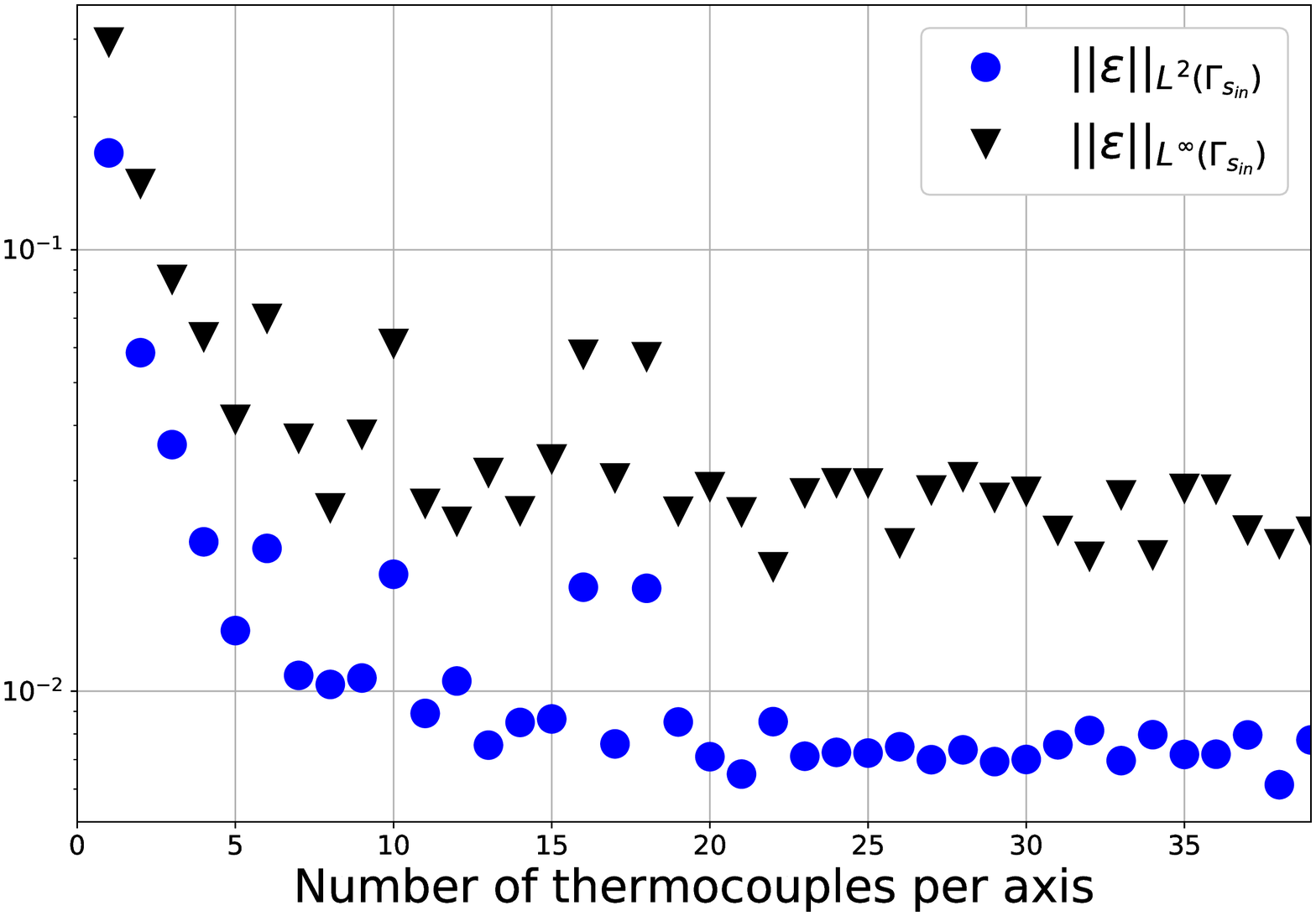}
        \caption{Measurements refinement.}
    \end{subfigure}%
    \caption{Behavior of the relative error norms (\ref{eq:analyticalBenchmark_errorNorms}) with the grid (A) and thermocouples number (B) refinement in the analytical benchmark case.}
\label{fig:analyticalBenchmark_y02_CGnoInt_inverseGridRefinement}
\end{figure}

We now investigate how the grid refinement influences the results.
Figure~\ref{fig:analyticalBenchmark_y02_CGnoInt_inverseGridRefinement} (a) shows the behavior of the relative error of the estimated heat flux with the grid refinement.
This test is performed with the stopping criterion $J_1 < J_{1_{tol}} = 1e-4 K^2$.
The error in general decreases by increasing the mesh refinement.
However, the decrease is not monotonic with a small increase for the $40^3$ elements grid.
To further investigate the convergence of the method, we tested the effect of increasing the number of thermocouples, keeping the same number of thermocouples along the x- and y-axis equal.
Figure~\ref{fig:analyticalBenchmark_y02_CGnoInt_inverseGridRefinement} (b) shows the obtained results for the $40^3$ elements grid.
Notice that the error converges non-monotonically. 

\subsubsection{Parameterization of the Boundary Condition}

We now test the performances of the parameterization method described in Section~\ref{sec:parameterizedBC_method}.
In particular, we consider the effects that the selection of the basis functions have on the results and the conditioning of the linear system (\ref{eq:linSys_parametrizedBC}).
Moreover, also in this case, we test the effect of the  mesh refinement on the estimated heat flux.

As already mentioned, we consider Gaussian RBFs as basis functions for the parameterization of the boundary heat flux.
Recalling (\ref{eq:GaussianRBF}), the basis functions are given by
\begin{equation*}
    \phi_j(\mathbf{x}) = e^{- \left(\eta \left\lVert \mathbf{x} - \pmb{\xi}_j  \right\rVert_2\right)^2},\quad \text{ for } j=1,2,\dots,M,
\end{equation*}
where we locate the centers $\pmb{\xi}_j$ at the projection of the virtual measurement points on the boundary $\Gamma_{s_{in}}$.

Both the choice of the basis functions (\ref{eq:GaussianRBF}) and of the position of their center are arbitrary.
However, they come suggested from the physics of the problem.
The Gaussian RBFs are selected because with their radial decay reduce the correlation between bases which are far away.
For a similar reason, the RBFs are centered at the projection of the measurements to have a relationship between bases and measurements.
This reasoning applies well to CC molds because we have the thermocouples located in a surface parallel and close to the boundary where we want to estimate the heat flux.
In a more general scenario, these choices lose their motivation.

\begin{figure}[!htb]
    \centering
    \begin{subfigure}[htb]{.5\linewidth}
        \centering
        \includegraphics[width=\textwidth]{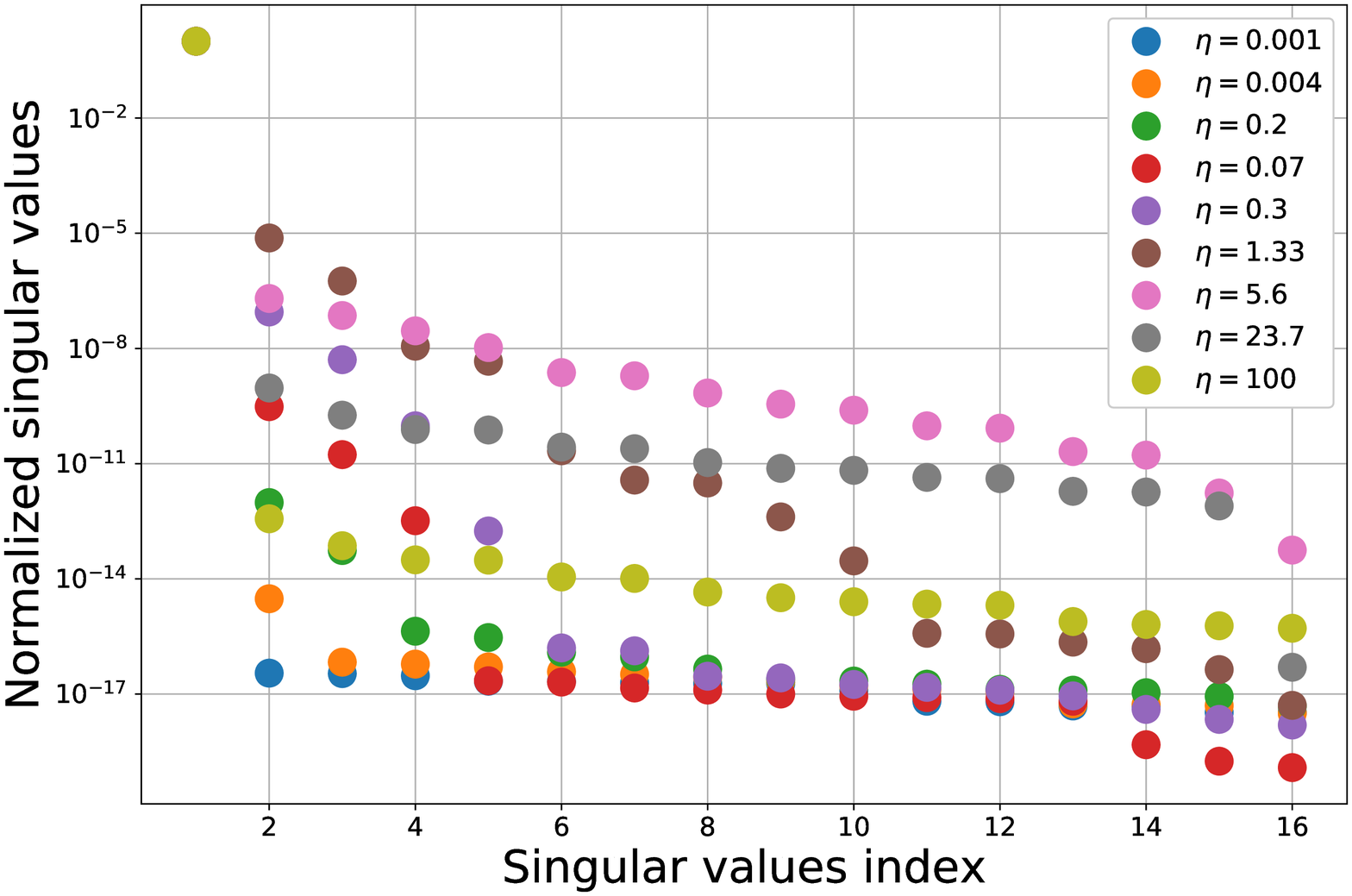}
        \caption{Decay of the singular values.}
    \end{subfigure}%
    \begin{subfigure}[htb]{.5\linewidth}
        \centering
        \includegraphics[width=\textwidth]{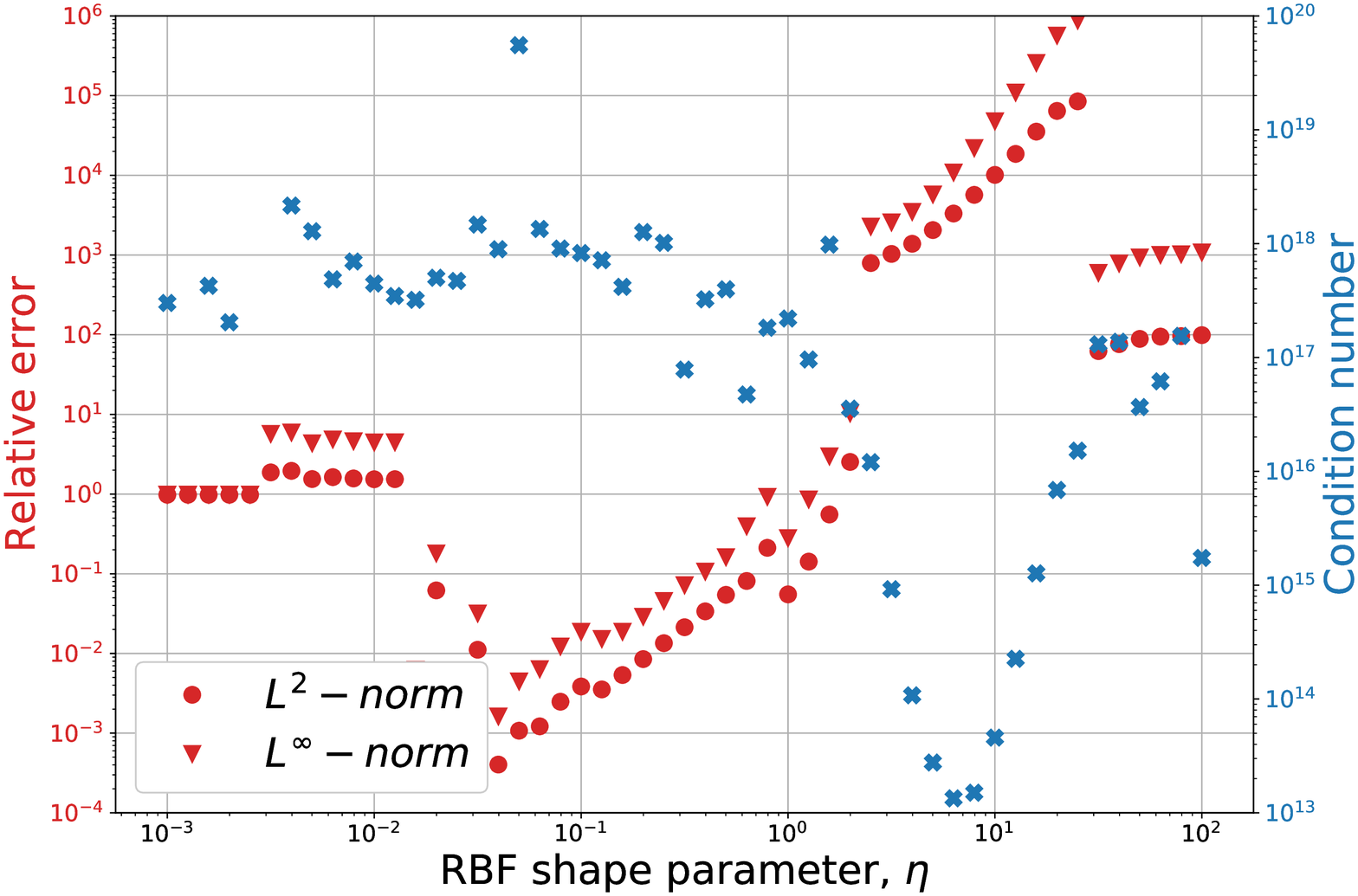}
        \caption{Relative error and condition number.}
    \end{subfigure}%
    \caption{Effect of the RBFs shape parameter, $\eta$, on (A) the normalized singular values of the matrix $\Theta^T \Theta$  and on (B) the $L^2$- and $L^\infty$-norms of the relative error  and on the condition number of the linear system.}
    \label{fig:analyticalBenchmark_y02_paramBC_RBFshapeParameterEffects}
\end{figure}

To completely define the basis functions, we still must tune the shape parameter $\eta$. 
Then, the first analysis we perform is the influence of $\eta$ on the invertibility of system (\ref{eq:linSys_parametrizedBC}) and on the boundary heat flux estimation.
This parameter controls the decay of the RBF.
For bigger (smaller) values of $\eta$ the decay is faster (slower).
Figure~\ref{fig:analyticalBenchmark_y02_paramBC_RBFshapeParameterEffects} (a) shows the decay of the normalized singular values of $\Theta^T\Theta$ for different $\eta$.
The singular values are normalized by dividing them all by the first one.
In general, we can see that to bigger values of the shape parameter, correspond a slower decay of the singular values.

\begin{figure}[!htb]
    \begin{subfigure}[htb]{\linewidth}
        \centering
        \begin{subfigure}[c]{.33\linewidth}
            \includegraphics[width=.95\textwidth]{images/analyticalBenchmark_y02_CGnoInt_gTrue.png}
            \caption{$g_{an}$}
        \end{subfigure}%
        \begin{subfigure}[c]{.33\linewidth}
            \includegraphics[width=.95\textwidth]{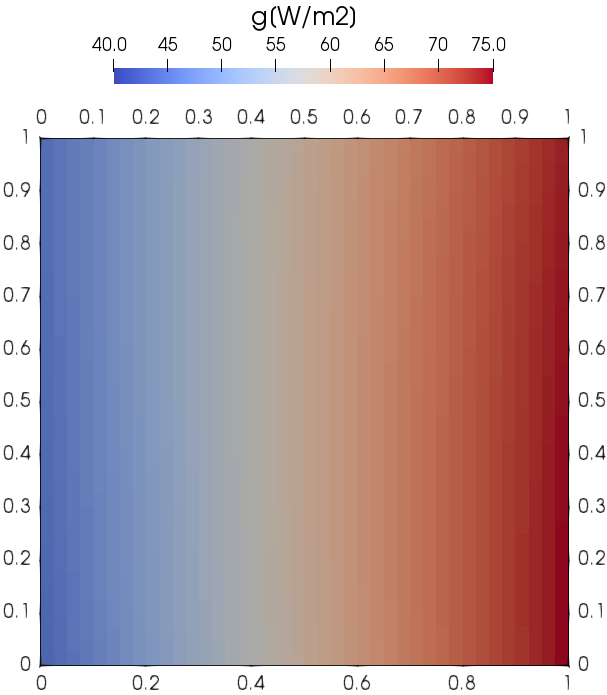}
            \caption{Estimated}
        \end{subfigure}%
        \begin{subfigure}[c]{.33\linewidth}
            \includegraphics[width=.95\textwidth]{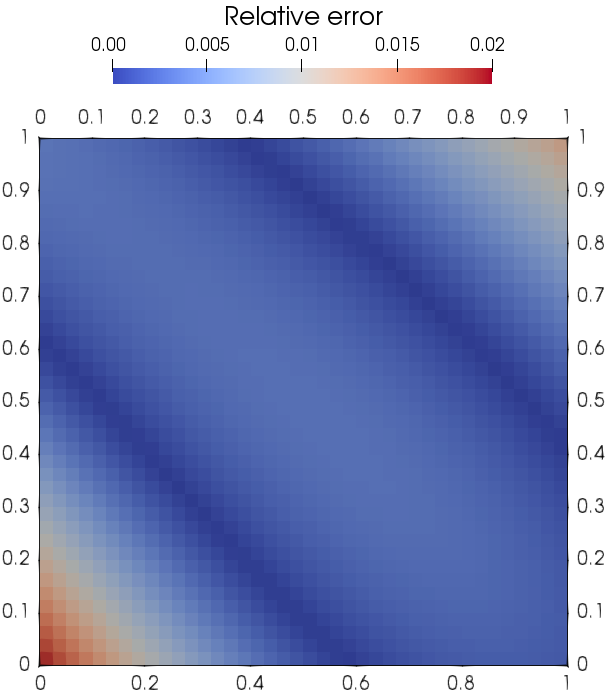}
            \caption{Relative Error}
        \end{subfigure}%
    \end{subfigure}%
    \caption{Comparison of the analytical (A) and estimated (B) boundary heat flux for the analytical benchmark case.
    This result is obtained by using the parameterization method with RBF shape parameter $\eta = 0.1$.}
\label{fig:analyticalBenchmark_y02_paramBC_heatFlux}
\end{figure}

Figure~\ref{fig:analyticalBenchmark_y02_paramBC_RBFshapeParameterEffects} (b) shows the condition number of the linear system (\ref{eq:linSys_parametrizedBC}).
The condition number is computed as the ratio between the bigger and the smaller singular value
\begin{equation}
  \kappa_{\Theta^T\Theta} = \frac{\sigma_{max}}{\sigma_{min}}.
    \label{eq:conditionNumber}
\end{equation}
The figure shows it together with the $L^2$- and $L^\infty$-norms of the relative error (\ref{eq:analyticalBenchmark_errorNorms}).
The method used for the solution of (\ref{eq:linSys_parametrizedBC}) is standard LU factorization with full pivoting.
In the figure, we see that the best results are obtained for $\eta = 0.1$ (see Figure~\ref{fig:analyticalBenchmark_y02_paramBC_heatFlux}).
Interestingly, looking at the behavior of the condition number, we can conclude that the quality of the results is not correlated to the conditioning of (\ref{eq:analyticalBenchmark_errorNorms}).

\begin{figure}[!htb]
    \centering
    \begin{subfigure}[htb]{.5\linewidth}
        \centering
        \includegraphics[width=\textwidth]{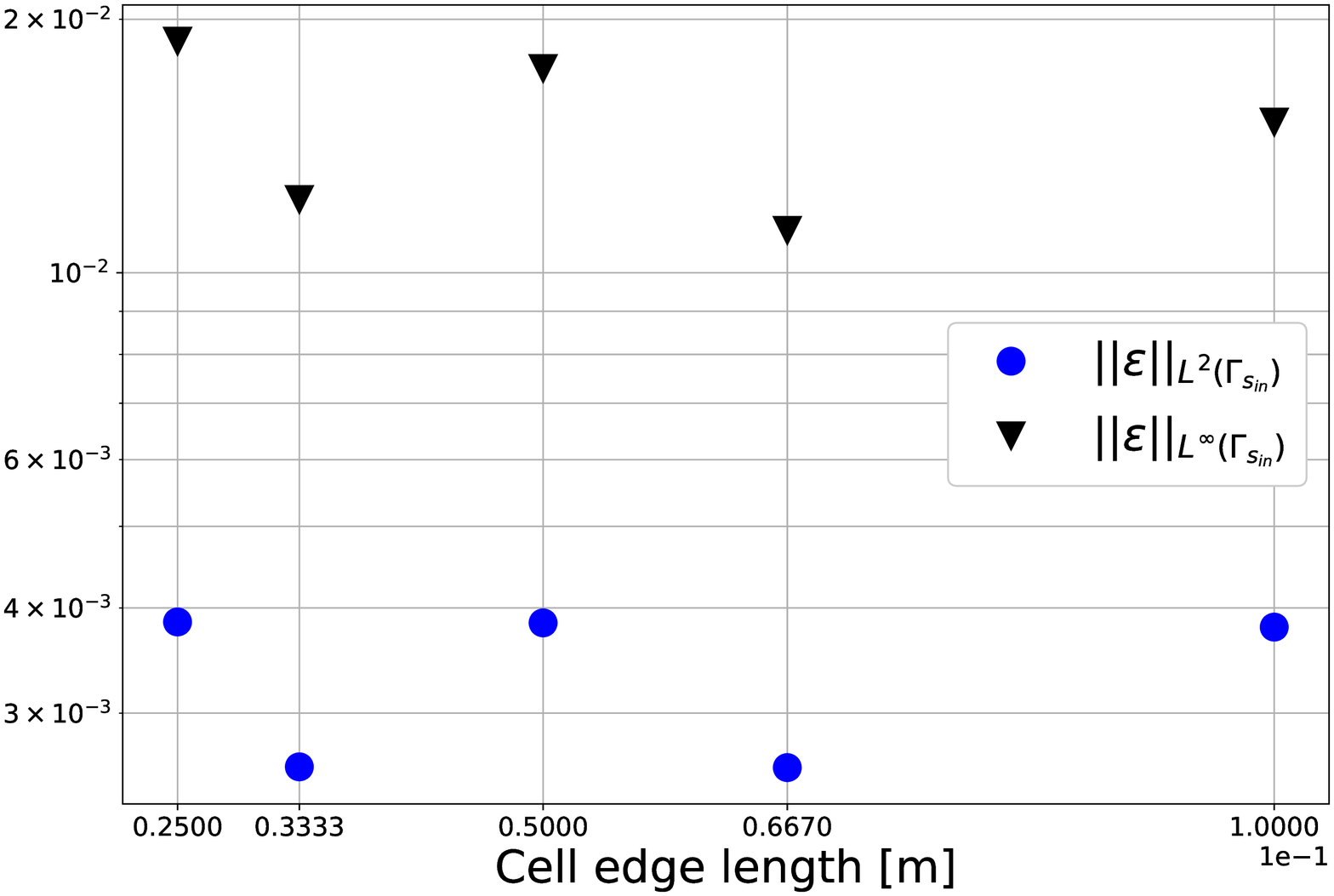}
        \caption{Grid refinement.}
    \end{subfigure}%
    \begin{subfigure}[htb]{.5\linewidth}
        \centering
        \includegraphics[width=\textwidth]{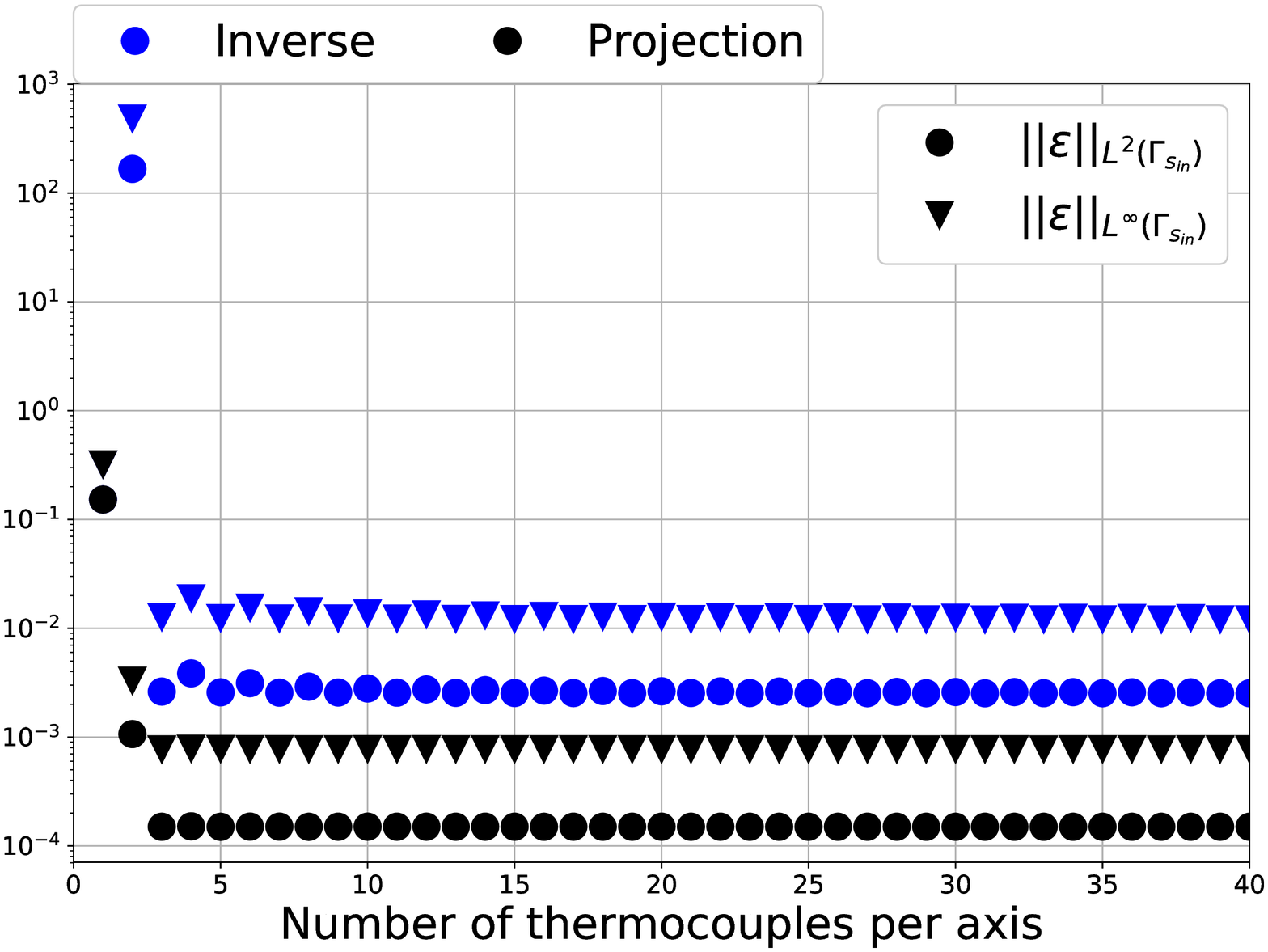}
        \caption{Measurements refinement.}
    \end{subfigure}%
    \caption{Behavior of the relative error norms (\ref{eq:analyticalBenchmark_errorNorms}) with the grid (A) and measurements (B) refinement using the parameterization method for the analytical benchmark case.
    In the figure (B), the blue results are obtained by solving the inverse problem and the black ones are the best possible approximation of the true heat flux in the parameterized space (remember that by increasing the number of thermocouples we increase the number of basis of the heat flux parameterization).}
    \label{fig:analyticalBenchmark_y02_paramBC_inverseGridRefinement}
\end{figure}

As for Alifanov's regularization, we test the effects of grid refinement on the estimated heat flux.
Figure~\ref{fig:analyticalBenchmark_y02_paramBC_inverseGridRefinement} (b) shows that we do not have a decrease of the relative error with grid refinement.
In fact, the error is oscillating between two very close values.
We obtain this result because the parameterization of the boundary heat flux is the same for all the grids and we reach the best description that the RBF parameterization can provide of the analytical heat flux as suggested by Figure~\ref{fig:analyticalBenchmark_y02_paramBC_inverseGridRefinement} (b).

\subsubsection{Noise in the Measurements}
In all previous tests, we considered the measurement to be free of noise.
This is not the real case.
In fact, thermocouples' measurements are notoriously noisy.
Thus, we analyze in this section the effects that the measurement noise have on the estimated heat flux, $g$.
From the industrial point of view, this analysis is of particular interest for our application.

We perform this analysis by adding to the measurements vector the Gaussian random noise $\mathbf{w}_n$
\begin{equation}
    \hat{\mathbf{T}}_{\mathbf{w}} = \hat{\mathbf{T}} + \mathbf{w}_n, \quad \mathbf{w}_n = \mathcal{N}(\pmb{\mu}, \Sigma),
\end{equation}
where, $\pmb{\mu}\in \R^M$ is the mean vector and $\Sigma\in \R^{M\times M}$ is the covariance matrix.
In particular, we choose $\mathbf{w}_n$ to be an independent and identically distributed (IID) random variable with zero mean, i.e. $\mathbf{w}_n = \mathcal{N}(\mathbf{0}, \omega^2 I)$.

To study the effect of noise, we perform several solutions of the inverse problem using $\hat{\mathbf{T}}_{\mathbf{w}}$ as thermocouples' measurements.
For each test, we compute 200 samples.
All these computations are done on the $40^3$ elements grid.
Then, we analyze the statistical and qualitative properties of the obtained results. 
In our first test, we analyze the behavior of the relative error (\ref{eq:analyticalBenchmark_errorNorms}) for different values of the noise standard deviation $\omega$.

Using Alifanov's regularization for the minimization of $J_1$, we must use a stopping criterion that regularize the solution.
In fact, the regularization parameter is the iteration counter $i$.
Here, we use the Discrepancy Principle (DP) as stopping criterion.\cite{Bardsley2018}
Thus, the iterations are stopped when 
\begin{equation}
    J_1[g^{i+1}] < \left(\frac{\omega^2 M}{2}\right)^2,
\end{equation}
where $M$ is the number of thermocouples.

\begin{figure}[!htb]
    \begin{subfigure}{\linewidth}
        \centering
        \begin{subfigure}[c]{.5\linewidth}
            \includegraphics[width=\textwidth]{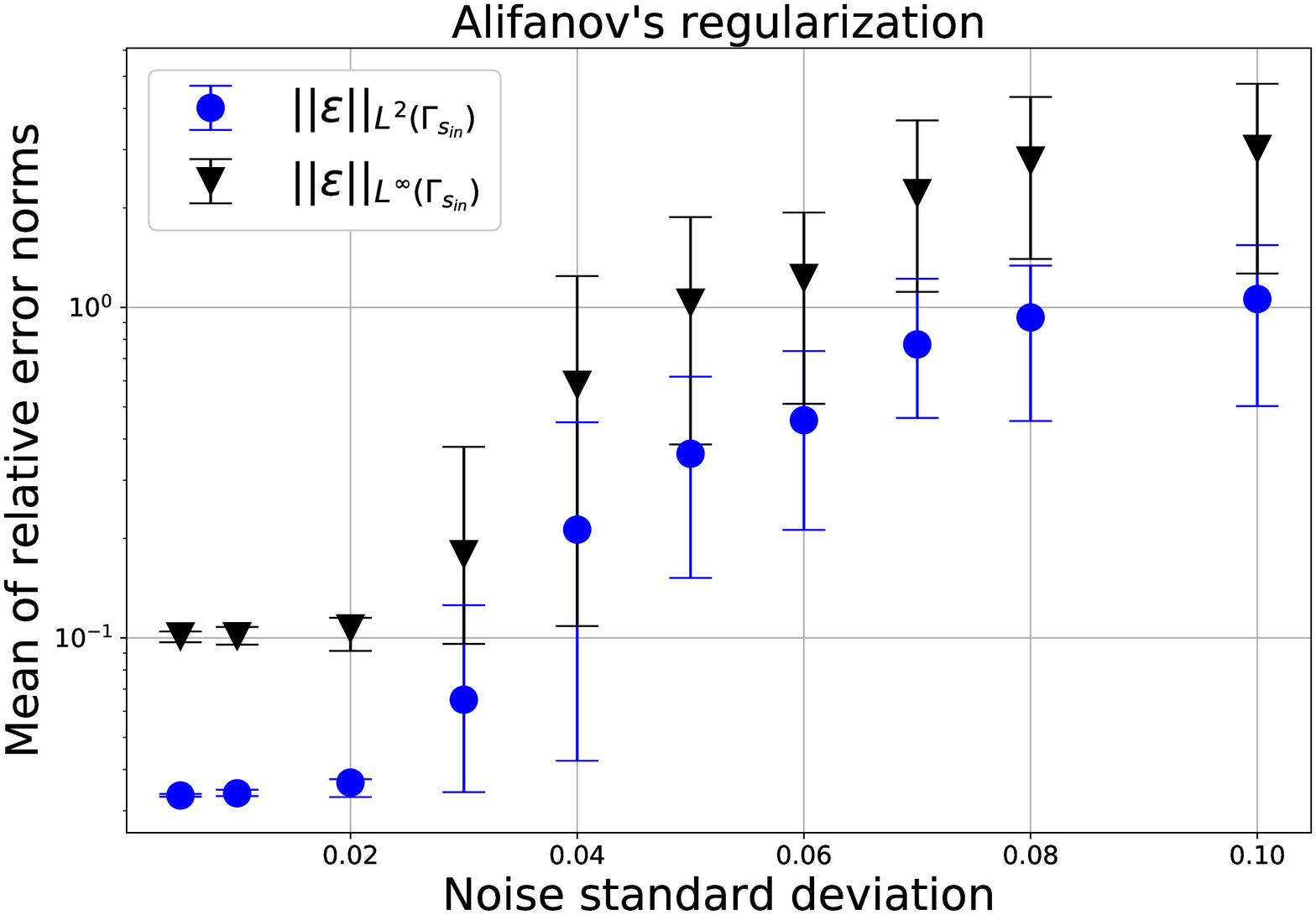}
            \caption{Alifanov's regularization}
        \end{subfigure}%
        \begin{subfigure}[c]{.5\linewidth}
            \includegraphics[width=\textwidth]{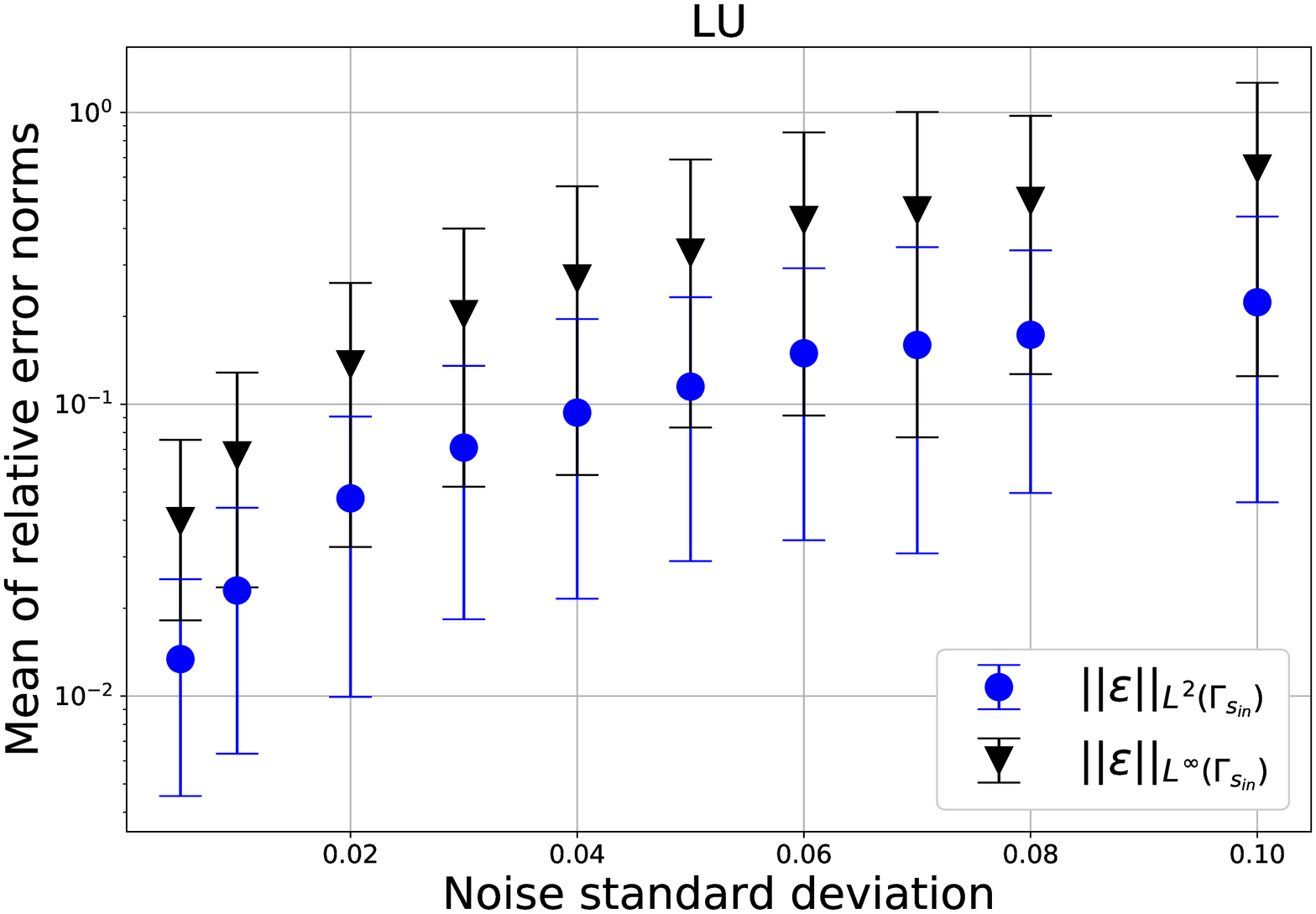}
            \caption{Parameterization method}
        \end{subfigure}%
    \end{subfigure}%
    \caption{Behavior of the relative error with respect to the standard deviation of the noise in the measurements for the Alifanov's regularization (A) and parameterization method (B) in the analytical benchmark case (90\% quantile bars shown).}
\label{fig:analyticalBenchmark_y02_noiseLevel}
\end{figure}

Figure~\ref{fig:analyticalBenchmark_y02_noiseLevel} illustrates the results of this first test.
We notice that Alifanov's regularization is able to filter the noise only for $\omega < 0.02$.
On the other hand, we see that for the parameterization method with LU factorization the results are spread around the mean value.
It suggests that the noise is propagating from the measurements into the solution.
As already mentioned, we require some regularization technique in solving (\ref{eq:linSys_parametrizedBC}) to regularize the solution.

As described in Section~\ref{sec:parameterizedBC_method}, we use TSVD regularization in the parameterization method.
We opt for this technique because it is effective when we have jumps in the singular values decay (see Figure~\ref{fig:analyticalBenchmark_y02_paramBC_RBFshapeParameterEffects}).
As already said, attention must be paid when using regularization techniques in selecting the regularization parameter.
In our case, the regularization parameter, $\alpha_{TSVD}$ is the number of singular values used in the truncation.
Different methodologies are available in the literature, e.g. unbiased predictive risk estimator, DP, L-curve, U-curve, generalized cross validation.\cite{Bardsley2018}
However, to show the dependency of the results on the regularization parameter, we performed numerical tests.

\begin{figure}[!htb]
    \centering
    \begin{subfigure}[htb]{\linewidth}
        \centering
        \begin{subfigure}[c]{.5\linewidth}
            \includegraphics[width=\textwidth]{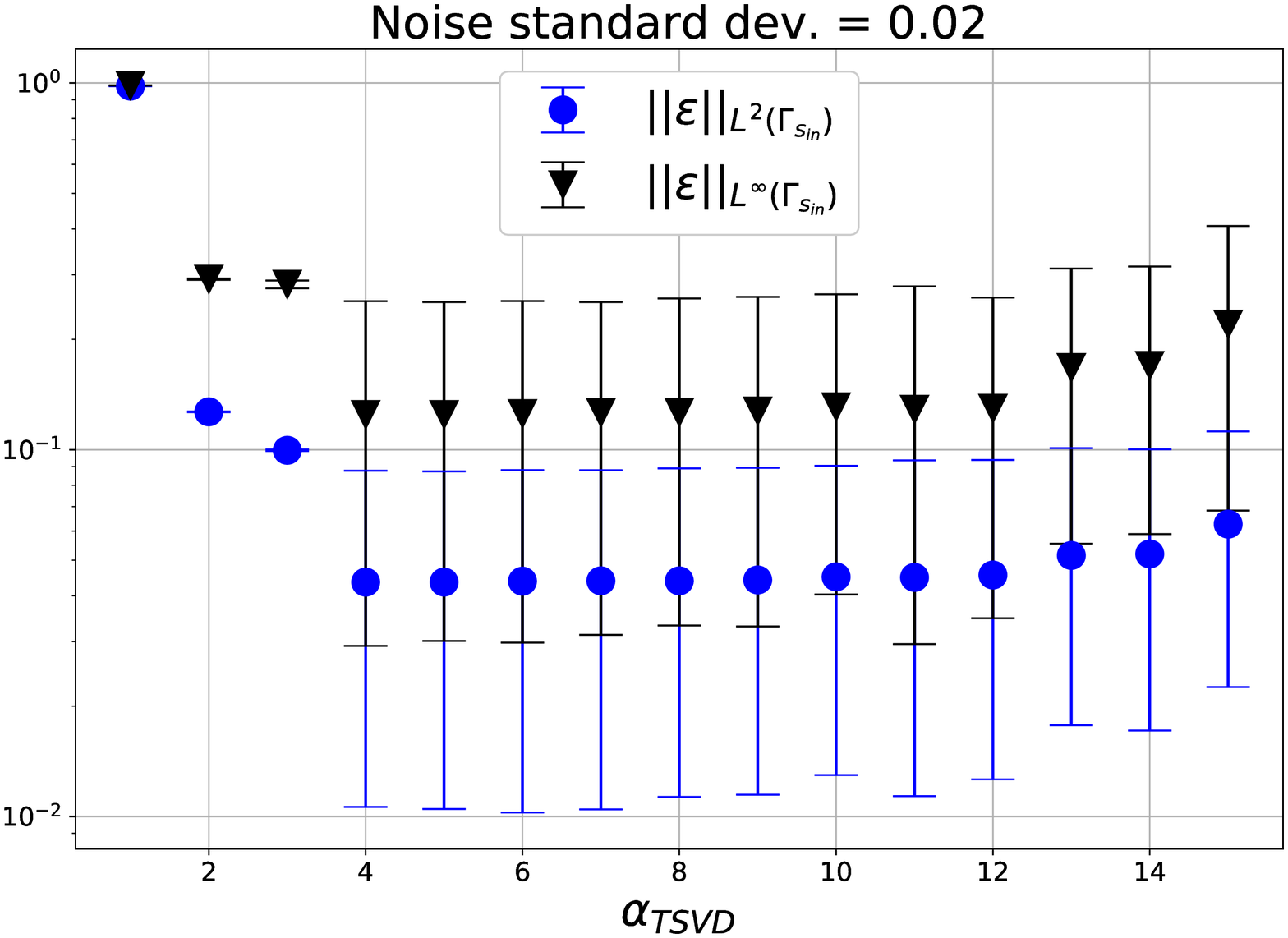}
        \end{subfigure}%
        \begin{subfigure}[c]{.5\linewidth}
            \includegraphics[width=\textwidth]{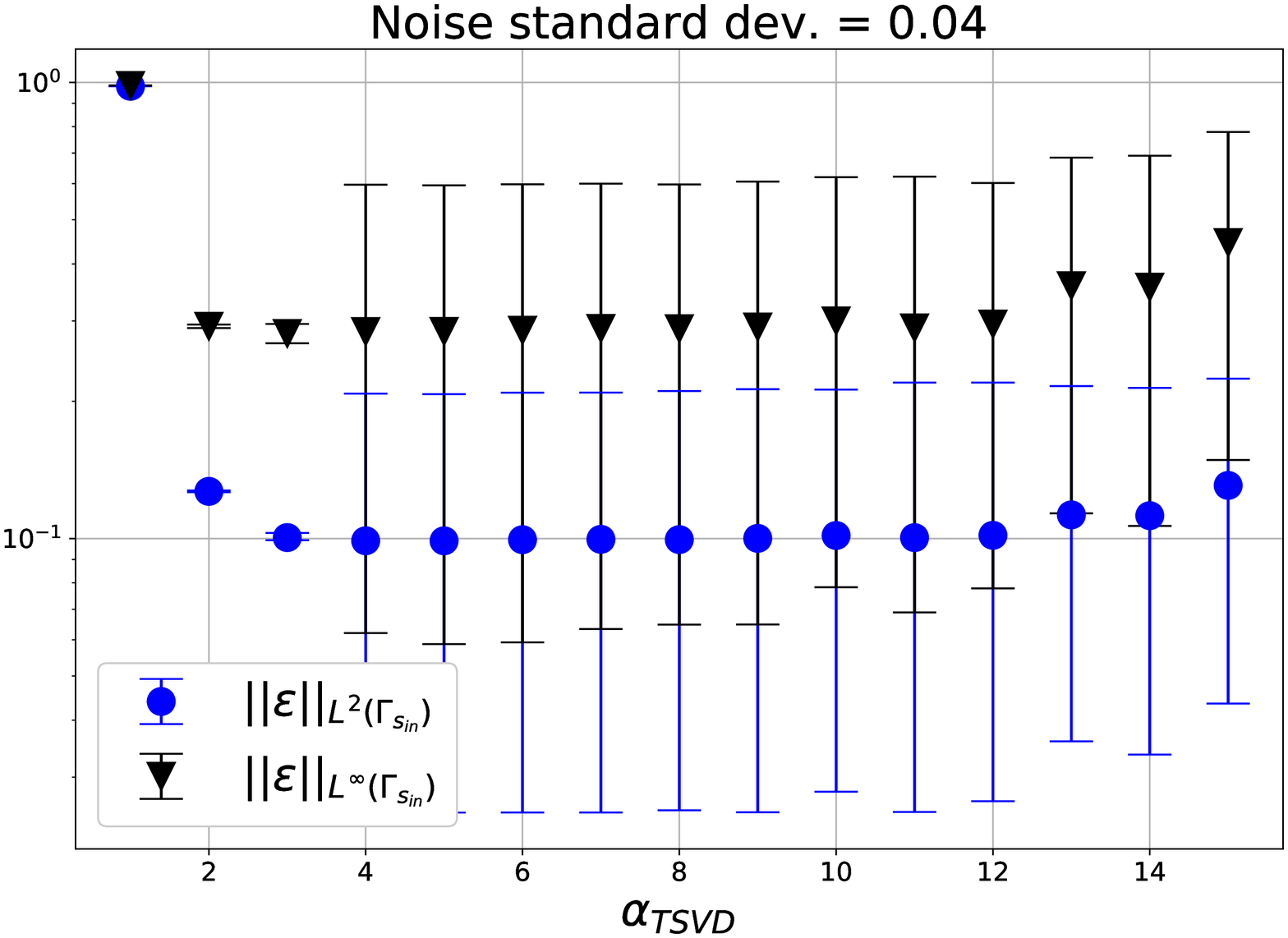}
        \end{subfigure}%
    \end{subfigure}%

    \centering
    \begin{subfigure}[htb]{\linewidth}
        \centering
        \begin{subfigure}[c]{.5\linewidth}
            \includegraphics[width=\textwidth]{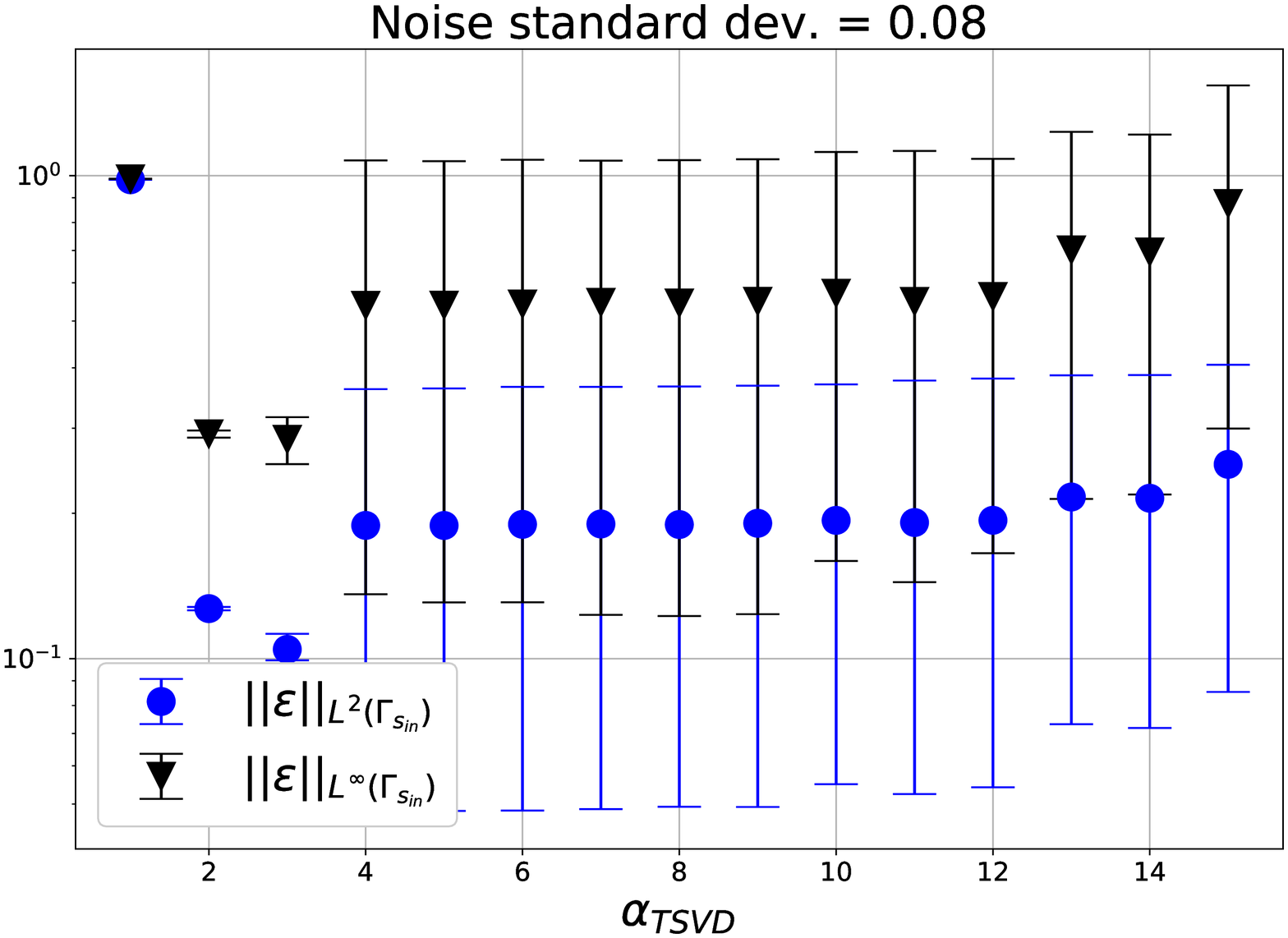}
        \end{subfigure}%
        \begin{subfigure}[c]{.5\linewidth}
            \includegraphics[width=\textwidth]{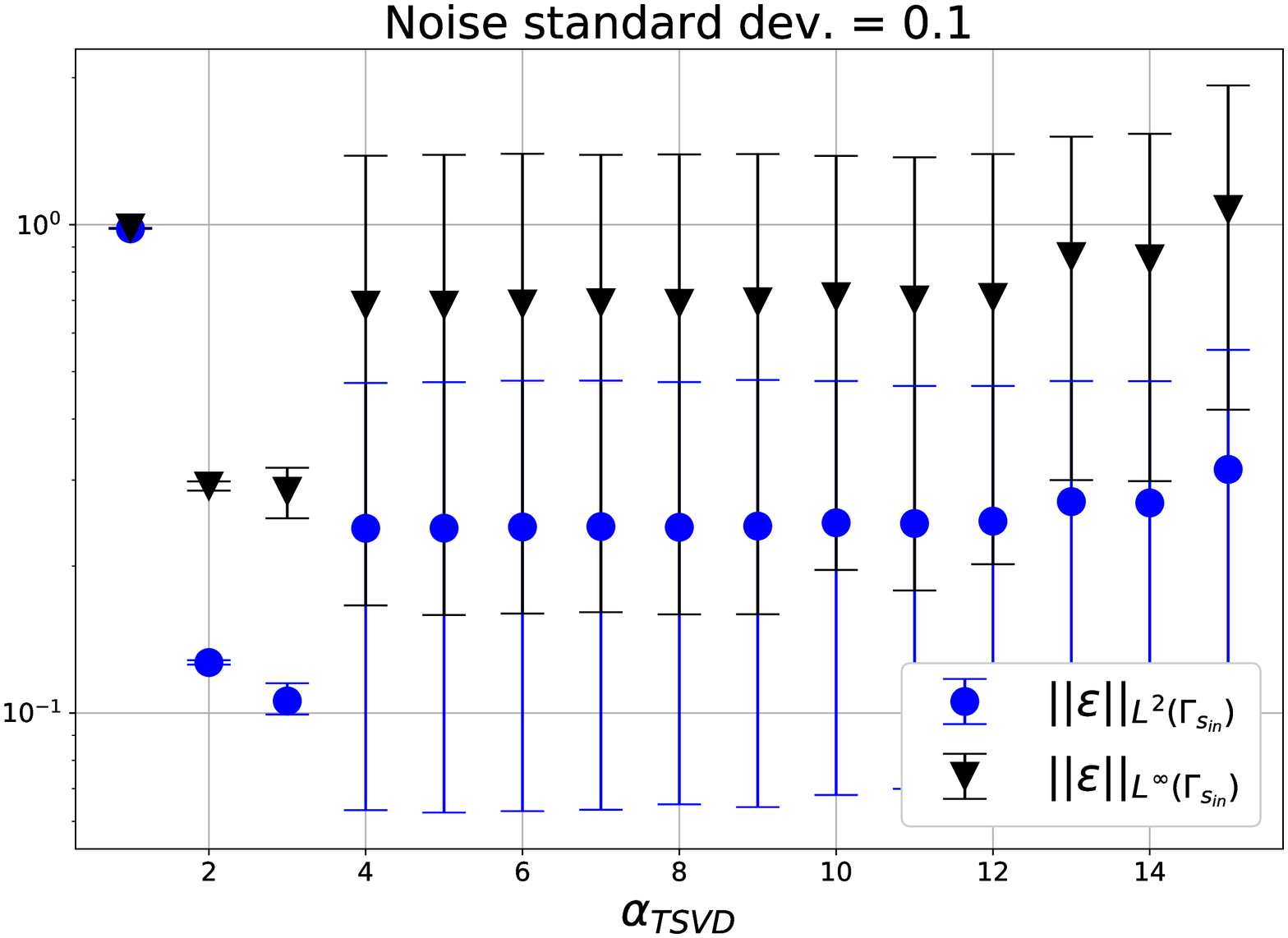}
        \end{subfigure}%
    \end{subfigure}%
    \caption{Effect of the regularization parameter $\alpha_{TSVD}$ using the TSVD in parameterization method for the analytical benchmark case (90\% quantile bars shown).}
\label{fig:analyticalBenchmark_y02_paramBC_noise_TSVDregPar}
\end{figure}

Figure~\ref{fig:analyticalBenchmark_y02_paramBC_noise_TSVDregPar} shows the behavior of the $L^2$- and $L^\infty$-norm of the relative error with respect to regularization parameter $\alpha_{TSVD}$, for different values of the noise standard deviation $\omega$.
As expected, the optimal value of the regularization parameter depends on the noise variance.
In fact, for low noise level we should use higher values of $\alpha_{TSVD}$ reducing it as the noise increases and vice versa.

\begin{figure}[!htb]
    \begin{subfigure}[htb]{\linewidth}
        \centering
        \begin{subfigure}[c]{.5\linewidth}
            \includegraphics[width=\textwidth]{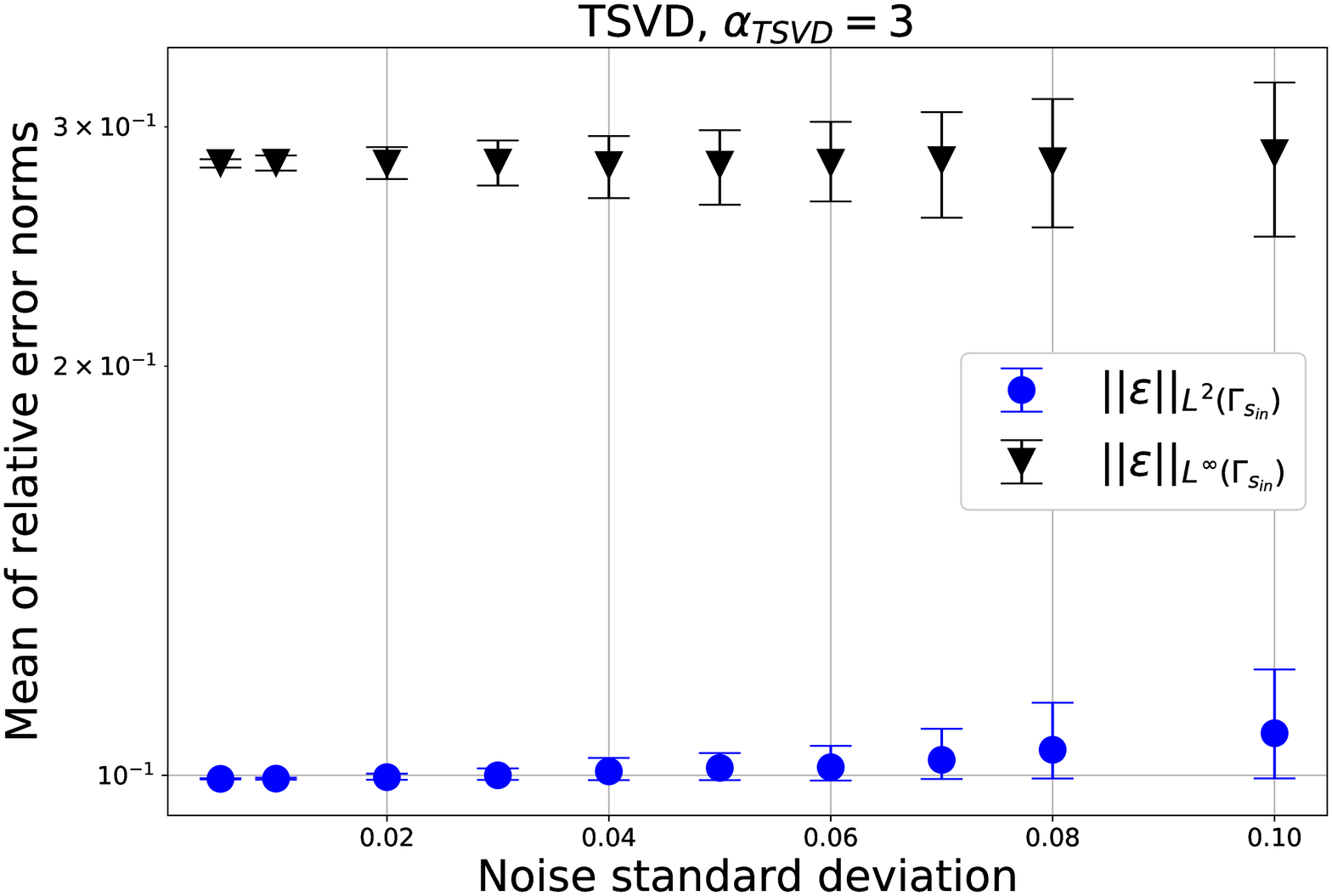}
        \end{subfigure}%
        \begin{subfigure}[c]{.5\linewidth}
            \includegraphics[width=\textwidth]{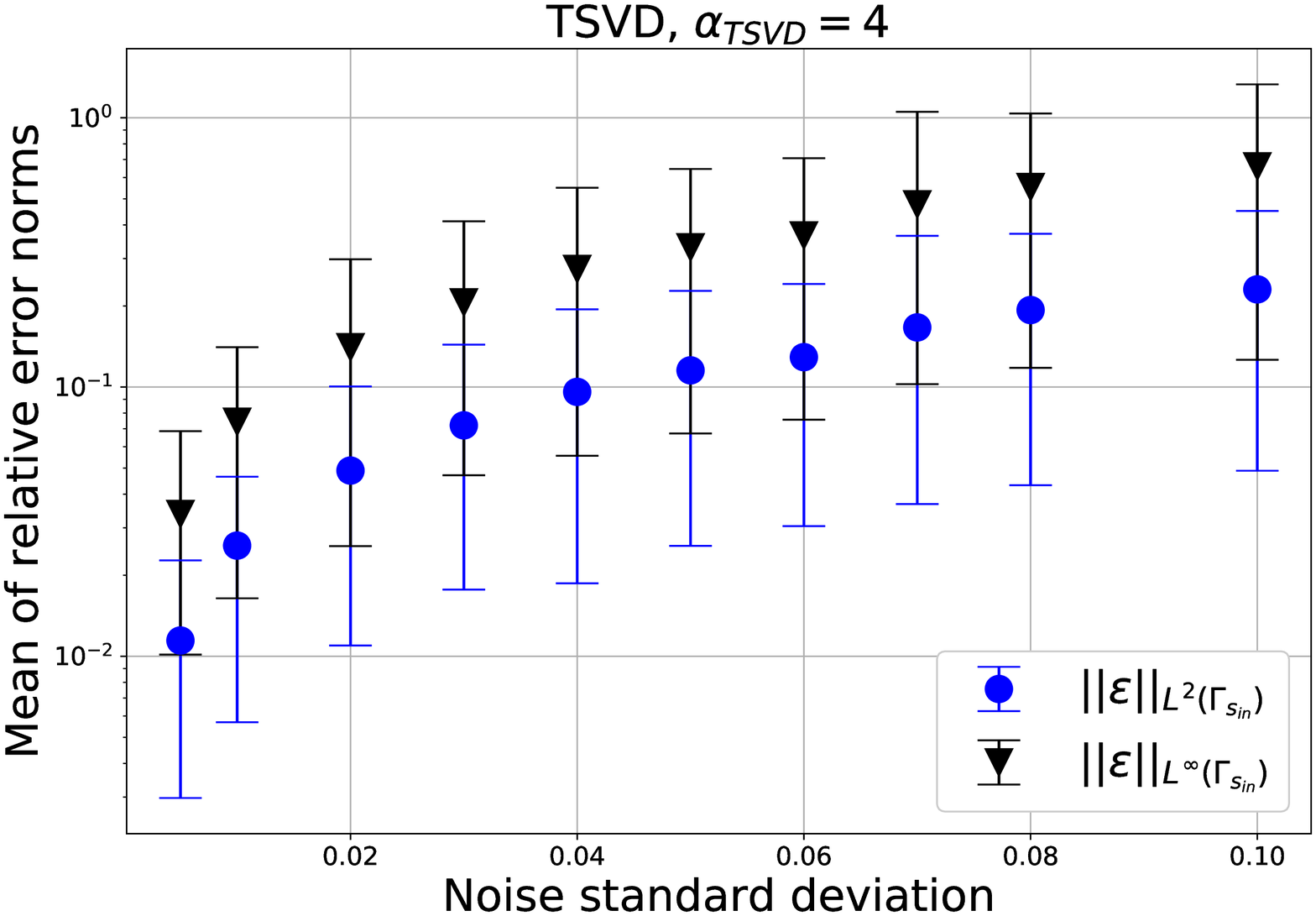}
        \end{subfigure}%
    \end{subfigure}%
    \caption{Behavior of the relative error with respect to the standard deviation of the noise in the measurements using the parameterization method with TSVD regularization in the analytical benchmark case (90\% quantile bars shown).}
\label{fig:analyticalBenchmark_y02_noiseLevel_TSVD}
\end{figure}

Testing again the TSVD regularization fixing $\alpha_{TSVD}$ and increasing the noise standard deviation, we clearly see the regularizing effect of the TSVD.
Figure~\ref{fig:analyticalBenchmark_y02_noiseLevel_TSVD} shows the obtained results.
In the figure, we appreciate the importance of a right choice of the regularizing parameter.
In fact, if we have very low noise in the measurements, we should opt for higher values of $\alpha_{TSVD}$ and vice versa.

\begin{figure}[!htb]
    \begin{subfigure}[htb]{\linewidth}
        \centering
        \begin{subfigure}[c]{.33\linewidth}
            \includegraphics[width=.95\textwidth]{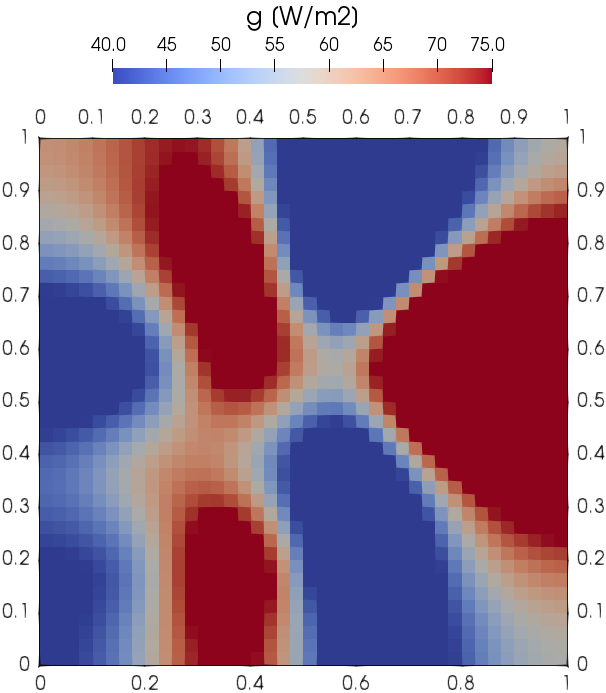}
        \end{subfigure}%
        \begin{subfigure}[c]{.33\linewidth}
            \includegraphics[width=.95\textwidth]{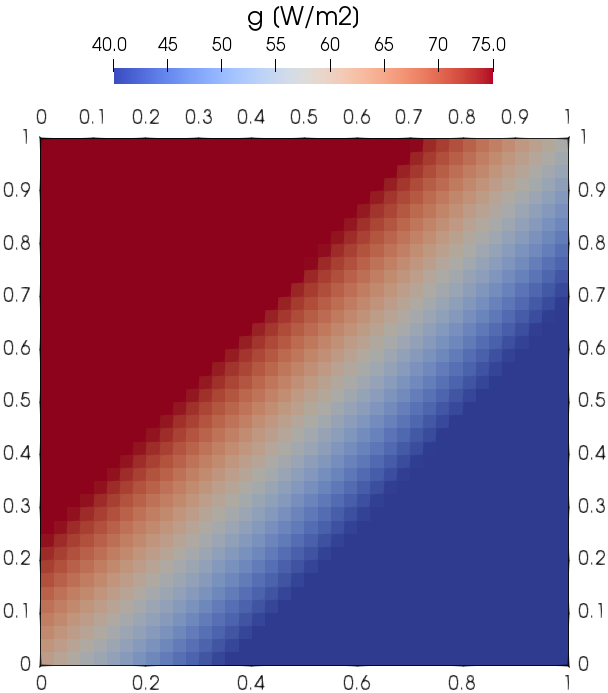}
        \end{subfigure}%
        \begin{subfigure}[c]{.33\linewidth}
            \includegraphics[width=.95\textwidth]{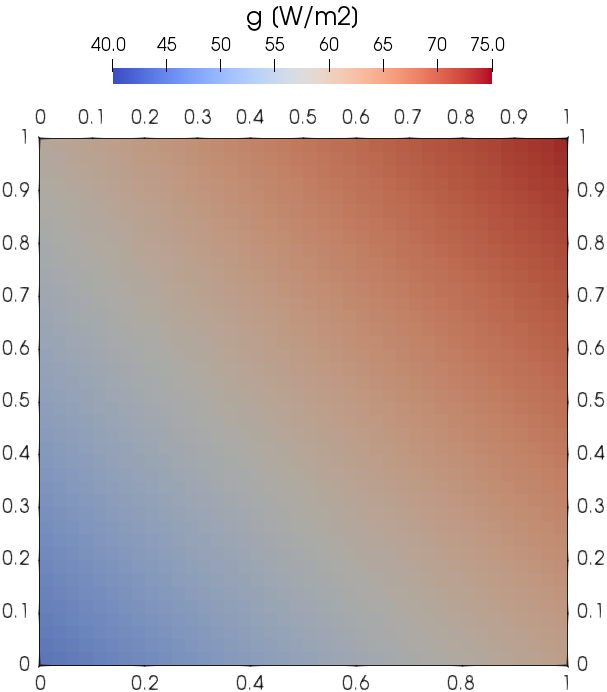}
        \end{subfigure}%
    \end{subfigure}%

    \begin{subfigure}[htb]{\linewidth}
        \centering
        \begin{subfigure}[c]{.33\linewidth}
            \includegraphics[width=.95\textwidth]{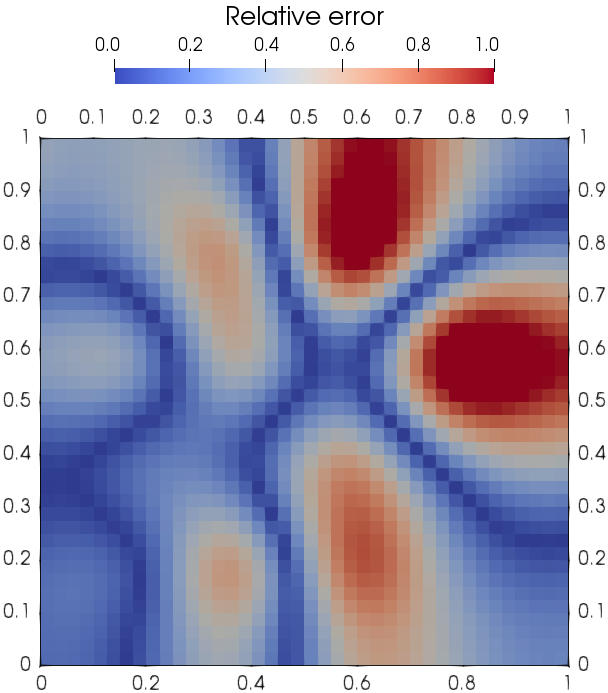}
            \caption{Alifanov's regularization}
        \end{subfigure}%
        \begin{subfigure}[c]{.33\linewidth}
            \includegraphics[width=.95\textwidth]{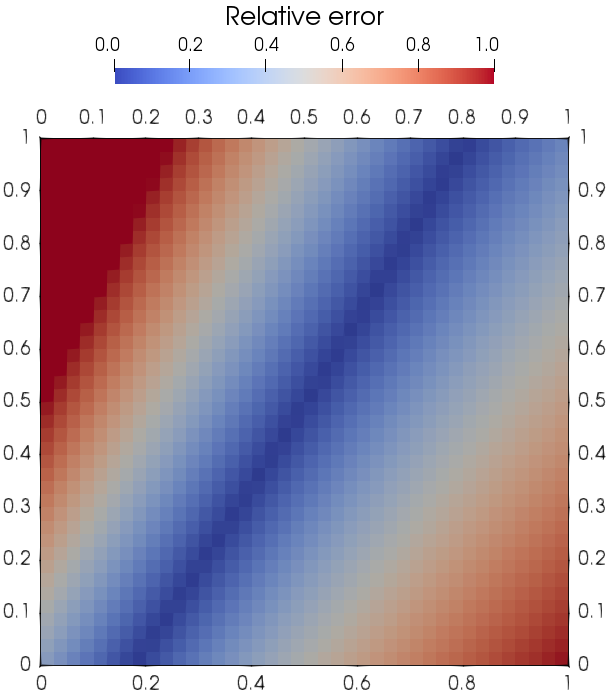}
            \caption{LU w. full pivoting}
        \end{subfigure}%
        \begin{subfigure}[c]{.33\linewidth}
            \includegraphics[width=.95\textwidth]{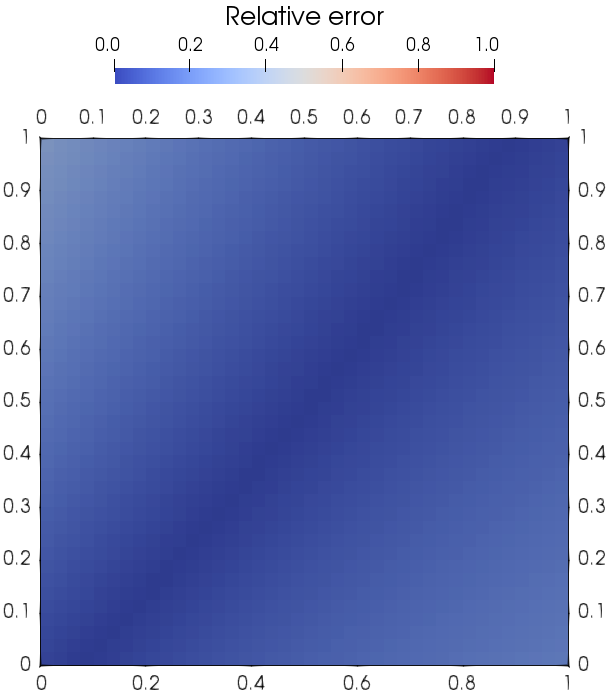}
            \caption{TSVD, $\alpha_{TSVD}=3$}
        \end{subfigure}%
    \end{subfigure}%
    \caption{Comparison of the estimated heat flux for Alifanov's regularization (A) and parameterization method with (C) and without (B) regularization for for the analytical benchmark case with noisy measurements  (noise standard deviation $\omega = 0.08$).}
\label{fig:analyticalBenchmark_y02_noiseWall}
\end{figure}

We conclude this noise analysis by looking at a realization of the computed heat flux with the different methods.
Figure~\ref{fig:analyticalBenchmark_y02_noiseWall} provides a qualitative example of the performances of the inverse solvers for $\omega = 0.08$.
As expected, the noise is not well filtered by the Alifanov's regularization while the parameterized method with TSVD provides a smooth solution in good agreement with the true value. 

\subsection{Inverse Problem with Temperature and Total Heat Flux Measurements}
In this section, we discuss the numerical solution of the inverse Problem~\ref{inverseProblemTotalHeat_3DhcModelSteady} where $T[g](\mathbf{x}_i)$ is the solution of Problem~\ref{prob:analyticalBenchmark} at points $\mathbf{x}_i$, for all $i=1,2,\dots,M$ and $\hat{G} = \int_{\Gamma_{s_{in}}} g_{an} d\Gamma$, $g_{an}$ being defined by (\ref{eq:analyticalHeatFlux}).
All computations are performed on the $40^3$ elements grid and the basis in the parameterization method are as in the previous section.

\begin{figure}[!htb]
    \begin{subfigure}[htb]{\linewidth}
        \centering
        \begin{subfigure}[c]{.5\linewidth}
            \includegraphics[width=\textwidth]{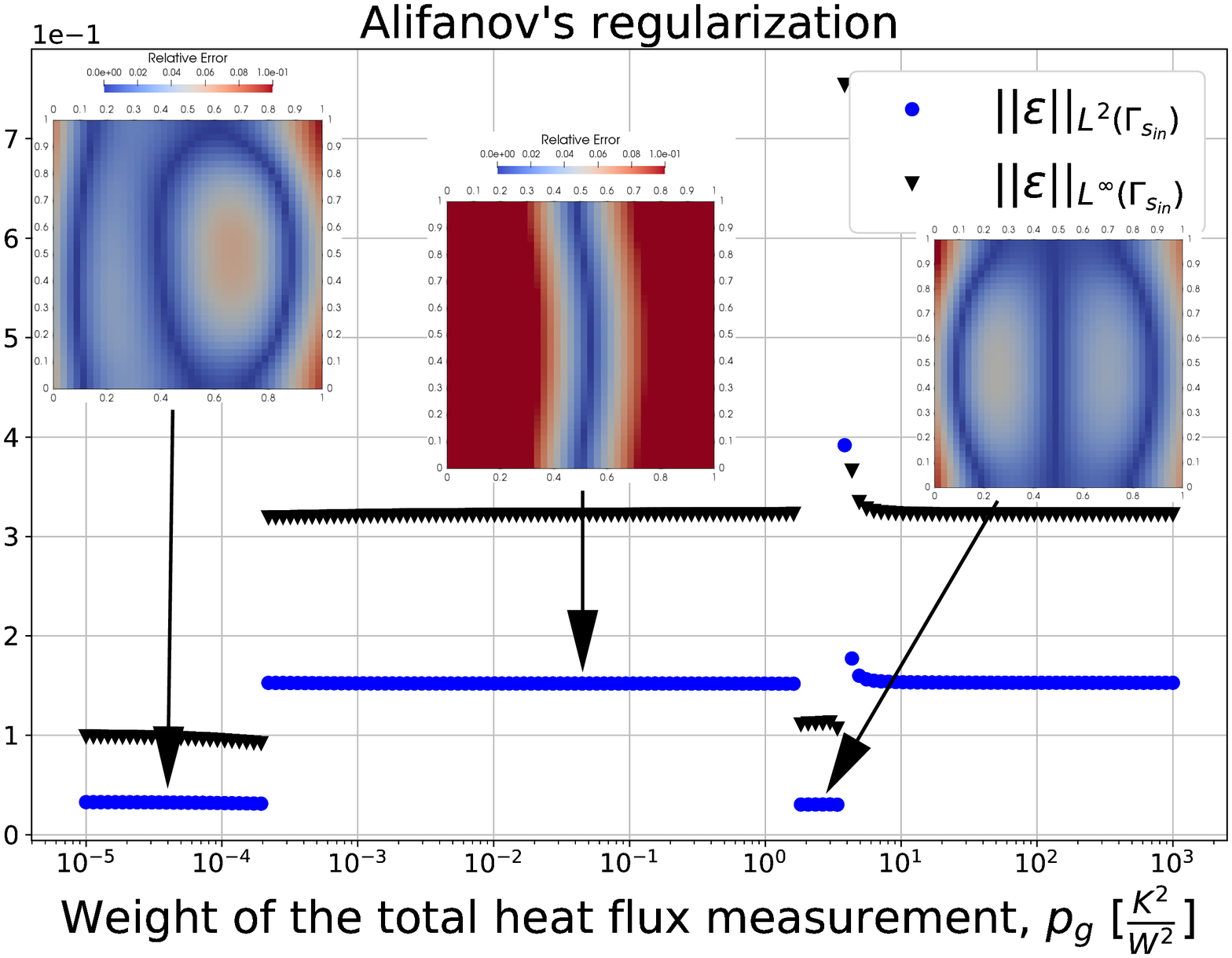}
            \caption{Alifanov's regularization.}
        \end{subfigure}%
        \begin{subfigure}[c]{.5\linewidth}
            \includegraphics[width=\textwidth]{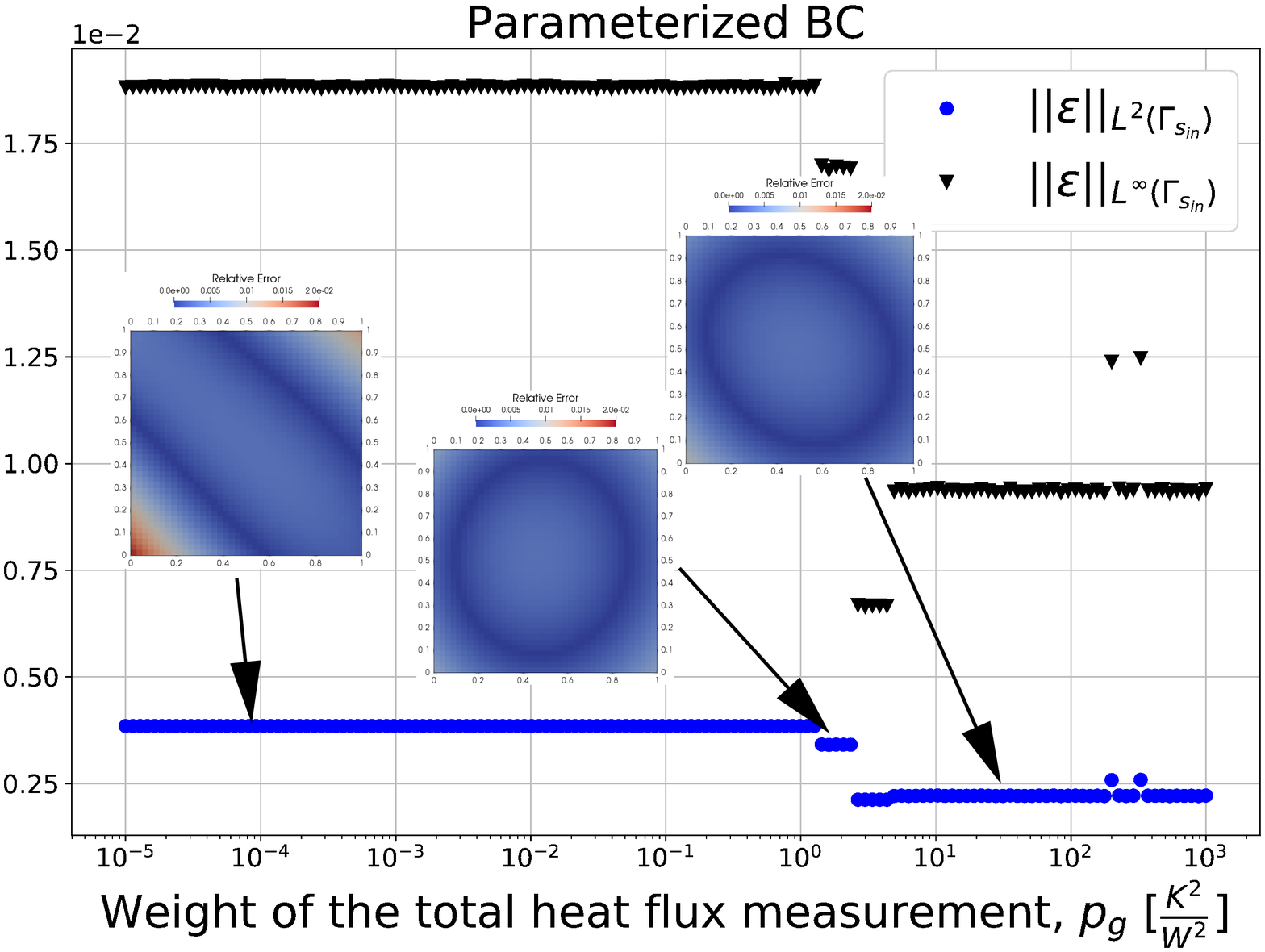}
            \caption{Parameterization method.}
        \end{subfigure}%
    \end{subfigure}%
    \caption{Behavior of the relative error with respect to the total heat measurement weight, $p_g$, using Alifanov's regularization (A) and parameterization of the heat flux with LU decomposition with full pivoting (B) in the analytical benchmark case.}
    \label{fig:analyticalBenchmark_y02_totalHeat_measureWeight}
\end{figure}

With respect to the previous section, we have one additional parameter: the total heat weight, $p_g$.
Since, it is not possible to set it a priori, we analyze its effects on the solution.
Figure~\ref{fig:analyticalBenchmark_y02_totalHeat_measureWeight} (a) shows the behavior of the $L^2$- and $L^\infty$-norm of the relative error for different values of $p_g$ using Alifanov's regularization for the solution of the inverse problem.
On the other hand, Figure~\ref{fig:analyticalBenchmark_y02_totalHeat_measureWeight} (b) shows the same graph for the parameterization method with LU decomposition with full pivoting.
These computations are performed without errors in the measurements.

Comparing the two figures (notice the different order of magnitude on the y-axis), we see that adding the total heat measurement improves the boundary heat flux estimation only for the parameterization method.
In the Alifanov's regularization, we have a very small decrease of the relative error for $p_g$ about $1e-4 \frac{K^2}{W^2}$ before having a sudden negative jump.
In Figure~\ref{fig:analyticalBenchmark_y02_totalHeat_measureWeight} (a), we appreciate an interesting the jump in the error for $1.5 \frac{K^2}{W^2} < p_g < 3 \frac{K^2}{W^2}$.
For these values, we recast similar results to those obtained for $p_g < 1e-4 \frac{K^2}{W^2}$.

\begin{figure}
    \centering
    \includegraphics[width=0.6\textwidth]{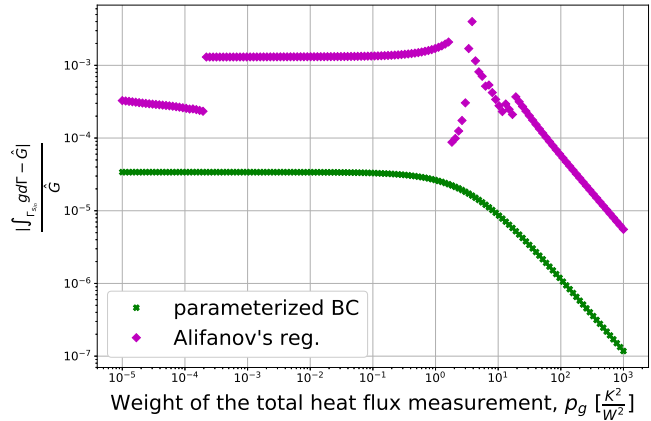}
    \caption{Behavior of the relative error in computing the total heat flux on $\Gamma_{s_{in}}$ with respect to $p_g$.}
    \label{fig:analyticalBenchmark_y02_totalHeat_integralError.eps}
\end{figure}

Figure~\ref{fig:analyticalBenchmark_y02_totalHeat_integralError.eps} provides further information on the effect of $p_g$.
In the parameterization method, $p_g$ does not have any effect on the solution for $p_g < 1$.
Then for higher values of $p_g$, the relative error decreases linearly.
On the other hand, the figure confirms the interesting behavior of Alifanov's regularization for $1.5 \frac{K^2}{W^2} < p_g < 4 \frac{K^2}{W^2}$.
However, it shows also for this method an almost linear decrease of the total heat relative error for $p_g > 4 \frac{K^2}{W^2}$.

\subsection{Conclusions}

To draw final conclusions on the performances of the tested inverse solvers, we compare their computational cost.
This is of particular interest in our research because we want to achieve real-time performances.
Table~\ref{tab:analyticalBenchmark_inverseProblem_timeComparison} illustrates the CPU time required for the computations with no error in the measurements and $J_{tol} = 10^{-4}~K^2$ in the case of only temperature measurements available.
Notice that all the computations were performed in serial on a Intel\textsuperscript{\textregistered} Core\textsuperscript{\texttrademark} i7-8550U CPU processor.

\begin{table}[htb]
  \centering
  \caption{Inverse problem CPU time comparison for the analytical benchmark case.}
  \label{tab:analyticalBenchmark_inverseProblem_timeComparison}
  \begin{tabular}{ |l|c|c|c| }
    \hline
      & Alifanov's reg. &  \multicolumn{2}{c|}{Parameterized heat flux} \\
    \hline
      & & offline & online \\
    \hline
    CPU time  &  $18.8~s$ & $7.21~s$ & $0.0056~s$\\
    \hline
  \end{tabular}
\end{table}

These results confirm that the offline-online decomposition makes the parameterized heat flux method eligible for real-time applications.
On the other hand, the Alifanov's regularization requires several solutions of direct, adjoint and sensitivity problem, so it cannot be employed in real-time as it is.

With this final remark, we conclude that the parameterization method outclasses Alifanov's regularization both in the quality of the estimation provided and, in the robustness (with TSVD regularization) with respect to errors in the measurements.

Moreover, thanks to its offline-online decomposition, the parameterization method has proved to be able to achieve real-time computation.
In fact, it requires a computationally expensive offline computation in which we solve several direct problems.
In the online phase it is very fast since we only solve a linear system with dimension equal to the number of basis used in the parameterization of the heat flux.

Finally, we considered the case of having as data for the inverse problem also a total heat flux measurement.
The parameterization method results are improved under every aspect by introducing this additional data.
On the other hand, Alifanov's regularization is only slightly affected by this additional data.

\section{Industrial Benchmark}\label{sec:numericalBenchmark}
The benchmark case presented in this section is a numerical test case.
This benchmark is designed to mimic the real industrial scenario of a CC mold.
In particular, the domain is a simplification of a mold plate and the physical quantities have typical industrial values.
Also the thermocouples' number and positioning are those of a real mold.
Table~\ref{tab:numericalBenchmark_parameters} summarizes the physical properties for this test case and the chosen heat flux, $g_{true}$.

As for the previous benchmark, the direct problem is a steady-state heat conduction problem in a homogeneous isotropic solid with a rectangular parallelepiped domain.
The domain $\Omega$ is as in Figure~\ref{fig:analyticalBenchmarkDomainSchematic} with $\Gamma_{s_{ex}} = \Gamma_{I} \cup \Gamma_{II} \cup \Gamma_{III} \cup \Gamma_{IV}$.
The mathematical formulation of the direct problem is that of Problem~\ref{prob:3DhcModelSteady}.

\begin{table}[!htb]
\centering
\caption{Physical parameters of the industrial benchmark case.}
\label{tab:numericalBenchmark_parameters}
\begin{tabular}{ ll }
\hline
\textbf{Parameter}    &   \textbf{Value}\\
\hline
    Thermal conductivity, $k$ 	    & $300.0~W/(m K)$ \\
    Heat transfer coefficient, $h$  & $5.66e4~W/(m^2 K)$\\
    Water temperature, $T_f$ 	    & $303 + 8 (1.2 - z)~K$\\
    Heat flux (Figure~\ref{fig:numericalSteadyBenchmark_trueHeatFlux} (a)), $g_{true}$           & $1e5[2(x-1)^2 - 2 z  - 5]~W/m^2$\\
    $L$ 			    & $2~m$\\
    $W$				    & $0.1~m$\\
    $H$				    & $1.2~m$\\
\hline
\end{tabular}
\end{table}

\begin{figure}[!htb]
    \begin{subfigure}[htb]{\linewidth}
        \centering
        \begin{subfigure}{.5\linewidth}
            \centering
            \includegraphics[width=0.8\textwidth]{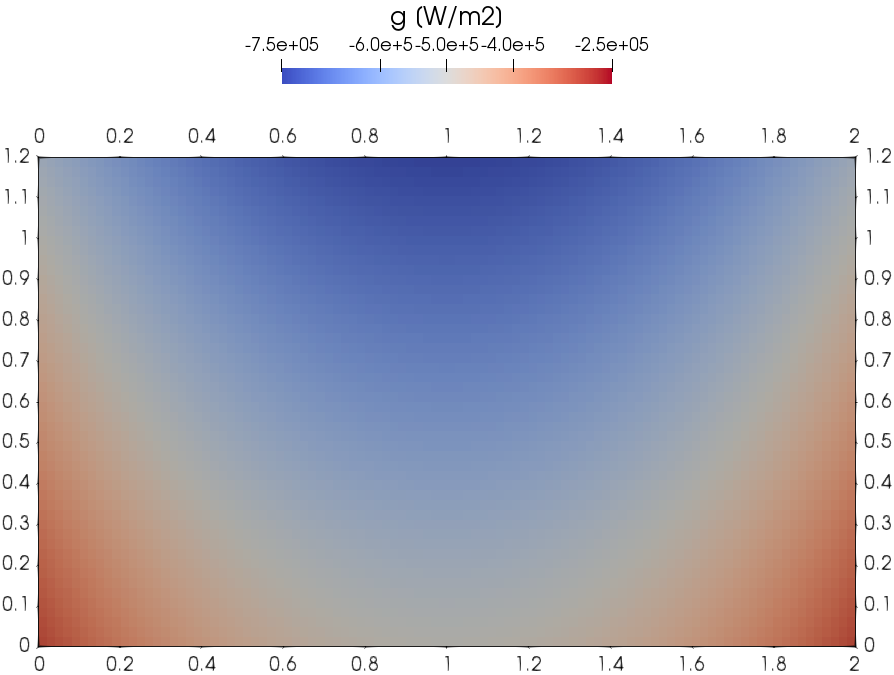}
            \caption{True heat flux.}
        \end{subfigure}%
        \begin{subfigure}{.5\linewidth}
            \centering
            \includegraphics[width=0.75\textwidth]{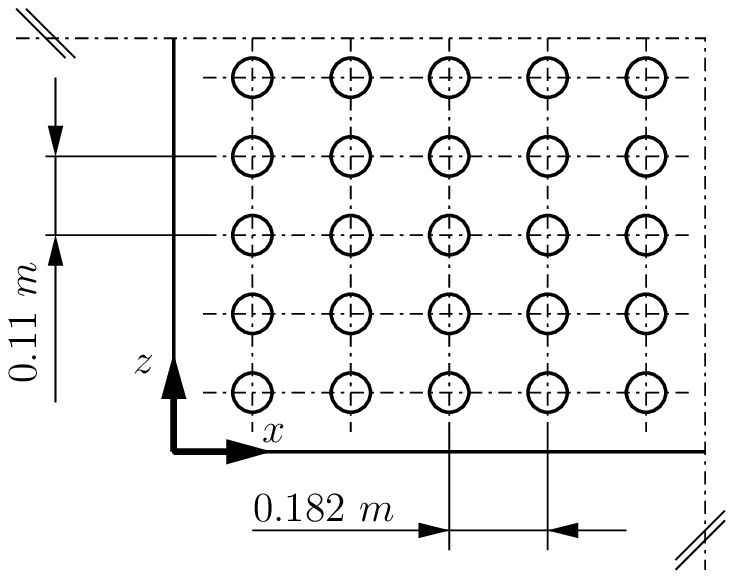}
            \caption{Thermocouples locations.}
        \end{subfigure}%
    \end{subfigure}%
    \caption{True heat flux (A) and position of the 100 thermocouples at the plane $y=0.02~m$ (B) for the industrial benchmark case.}
    \label{fig:numericalSteadyBenchmark_trueHeatFlux}
\end{figure}

For the discretization of the domain, we use a structured orthogonal grid with uniformly distributed elements along the three axes.
We use 200, 50 and 100 elements on the x-, y- and z-axis respectively.
Thus, the grid size is 1e6 elements.

The direct problem does not have an analytical solution.
Then for this benchmark, we assume that the direct problem is well solved and focus our attention on the solution of the inverse problem.

As in the real industrial case under study, we locate the virtual thermocouples in the plane $y=0.02~m$.
In this plane, they are equally distributed on the $x$- and $z$-axis as shown in Figure~\ref{fig:numericalSteadyBenchmark_trueHeatFlux} (b).

\subsection{Inverse Problem with Temperature Measurements}\label{section:numericalBenchmark_onlyTCmeas}

In the present section, we analyze the performance of the proposed methods for the solution of the inverse Problem~\ref{inverseProblem_3DhcModelSteady} for the introduced numerical test case.
First, we analyze the performances of Alifanov's regularization (see Section~\ref{sec:inverse_3DhcModelSteady_alifanov}).
Table~\ref{tab:numericalBenchmark_AlifanovParameters} shows the parameters used for the simulation.

\begin{table}[htb]
\centering
\caption{Parameters used in the Alifanov's regularization algorithm for the solution on the industrial benchmark case.}
\label{tab:numericalBenchmark_AlifanovParameters}
\begin{tabular}{ cl }
\hline
    \multicolumn{1}{l}{\textbf{Parameter}}    &   \textbf{Value}\\
\hline
    $g^0$       & $0~\frac{W}{m^2}$\\
    $J_{1_{tol}}$   & $1e2~K^2$\\
    $\frac{\norm{J_1^n - J_1^{n-1}}}{J^n}$ & $1e-2$\\
\hline
\end{tabular}
\end{table}

\begin{figure}[!htb]
    \begin{subfigure}{\linewidth}
        \centering
        \begin{subfigure}[c]{.5\linewidth}
            \includegraphics[width=.95\textwidth]{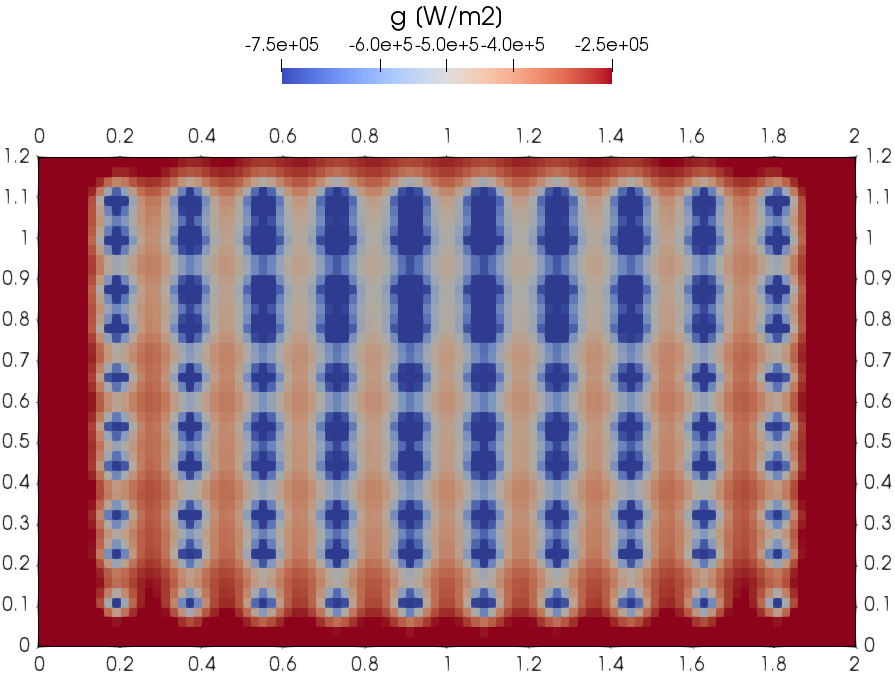}
            \caption{Iteration 1}
        \end{subfigure}%
        \begin{subfigure}[c]{.5\linewidth}
            \includegraphics[width=.95\textwidth]{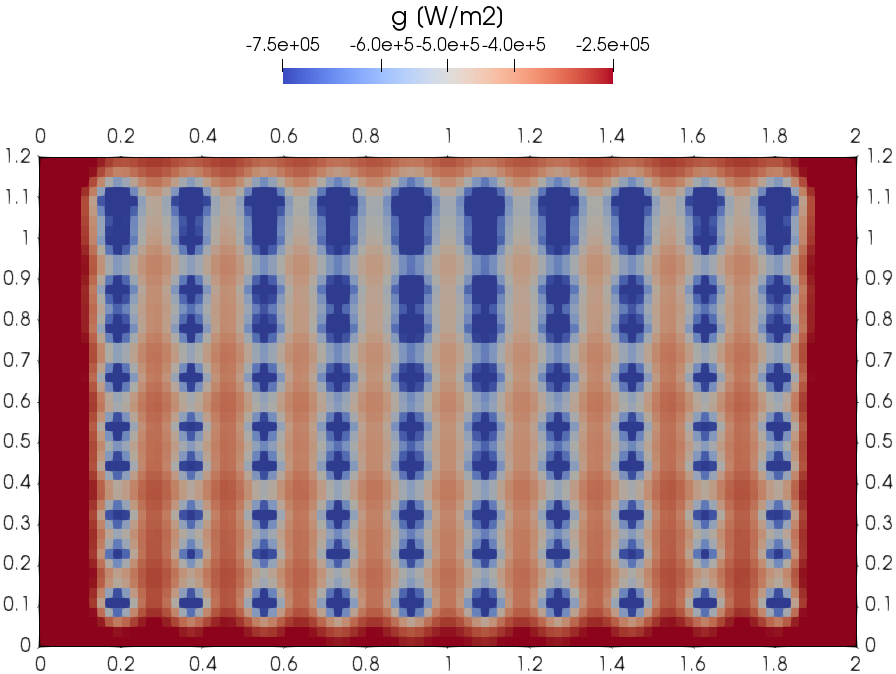}
            \caption{Iteration 80}
        \end{subfigure}%
    \end{subfigure}%
    \caption{Comparison of the computed heat flux by Alifanov's regularization at different iterations of the algorithm.}
\label{fig:numericalSteadyBenchmark_CGnoInt_heatFlux}
\end{figure}

Figure~\ref{fig:numericalSteadyBenchmark_CGnoInt_heatFlux} illustrates the estimated heat flux, $g$, at different iterations of the algorithm.
We notice that the algorithm provides a solution not in agreement with $g_{true}$.
In particular, it overestimates the heat flux close to the measurement points while far from the measurements the initial estimate is not modified. 
Moreover, increasing the number of iterations does not improve the results.
Due to the inability of estimating the heat flux also in the simplest case without measurement noise, we do not perform further tests with this method.

\begin{figure}[!htb]
    \centering
    \begin{subfigure}{0.5\linewidth}
        \centering
        \includegraphics[width=.95\textwidth]{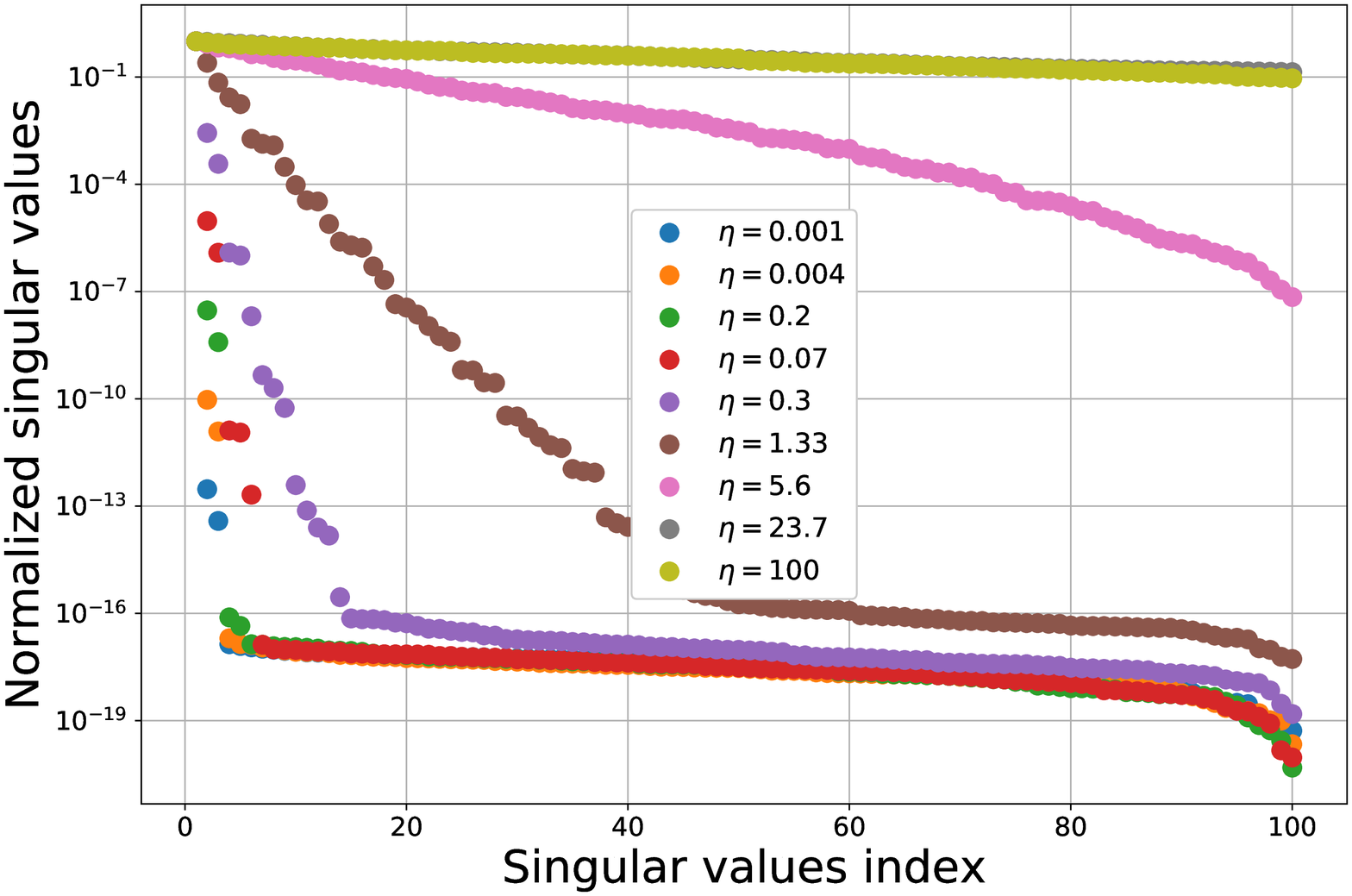}
        \caption{Decay of the singular values.}
    \end{subfigure}%
    \begin{subfigure}{0.5\linewidth}
        \centering
        \includegraphics[width=.9\textwidth]{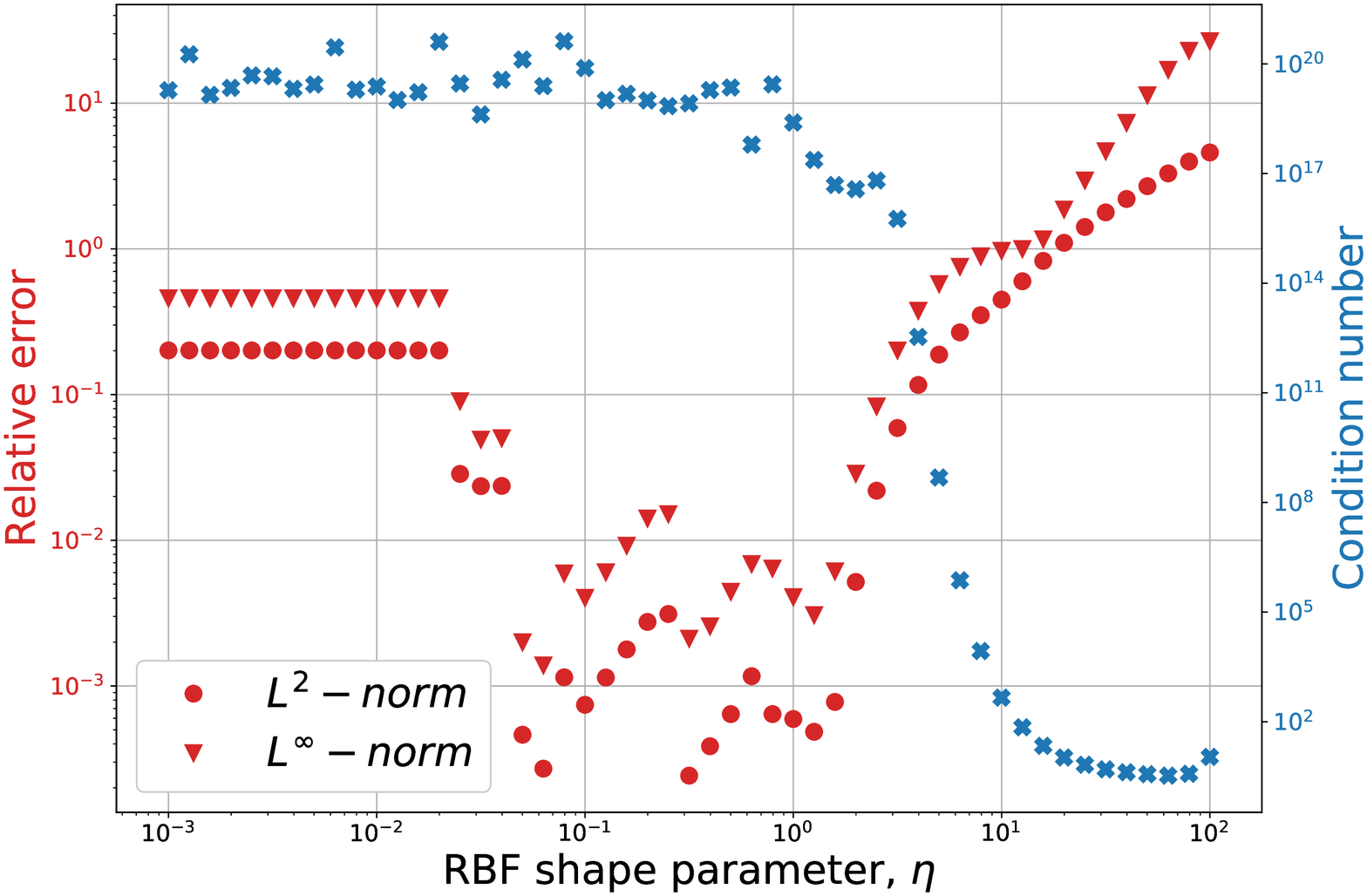}
        \caption{Relative error and condition number.}
    \end{subfigure}%
    \caption{Effect of the RBFs shape parameter on (A) the singular values of the matrix $\Theta^T \Theta$ and on (B) the relative error norms (\ref{eq:analyticalBenchmark_errorNorms}) using LU with full pivoting and on the condition number of (\ref{eq:linSys_parametrizedBC}).}
    \label{fig:numericalBenchmark_y02_paramBC_RBFshapeParameterEffects}
\end{figure}

We now consider the parameterization method of Section~\ref{sec:parameterizedBC_method}.
As for the previous benchmark, we start by performing a numerical analysis on the influence of the RBF shape parameter, $\eta$, on the invertibility of system (\ref{eq:linSys_parametrizedBC}) and on the estimated heat flux.
Figure~\ref{fig:numericalBenchmark_y02_paramBC_RBFshapeParameterEffects} (a) shows the decay of the singular values of $\Theta^T\Theta$ for different $\eta$.
As for the previous test case, to bigger values of the shape parameter correspond a slower decay of the singular values.
Moreover, we see from this singular value decay and in Figure~\ref{fig:numericalBenchmark_y02_paramBC_RBFshapeParameterEffects} (b) that for $\eta > 1$ the condition number of the system decreases.
However, the relative error of the heat flux estimation increases significantly for these values of $\eta$.

\begin{figure}[!htb]
    \centering
    \begin{subfigure}[c]{.5\linewidth}
        \includegraphics[width=.95\textwidth]{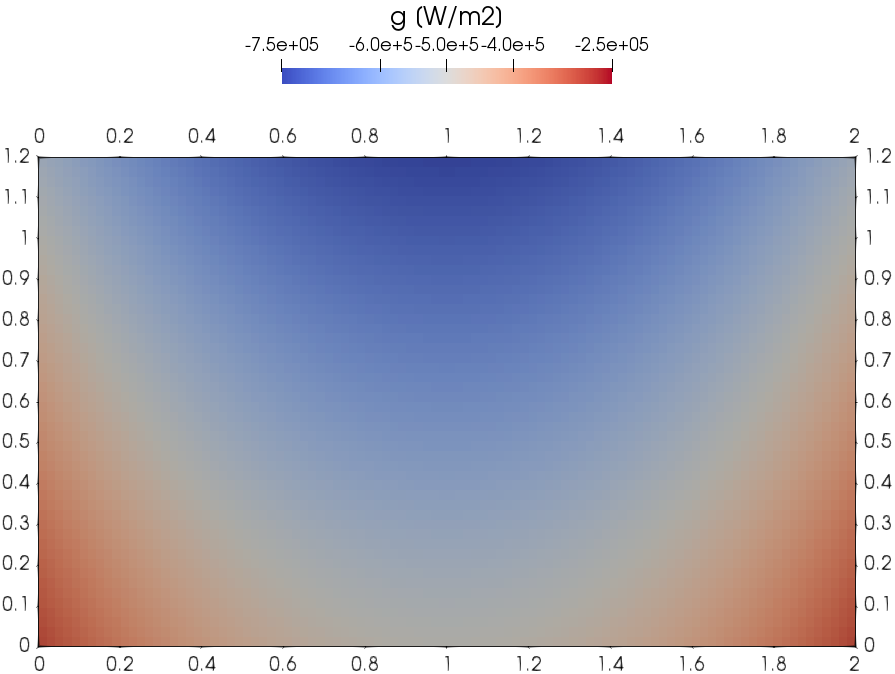}
        \caption{Heat flux}
    \end{subfigure}%
    \begin{subfigure}[c]{.5\linewidth}
        \includegraphics[width=.95\textwidth]{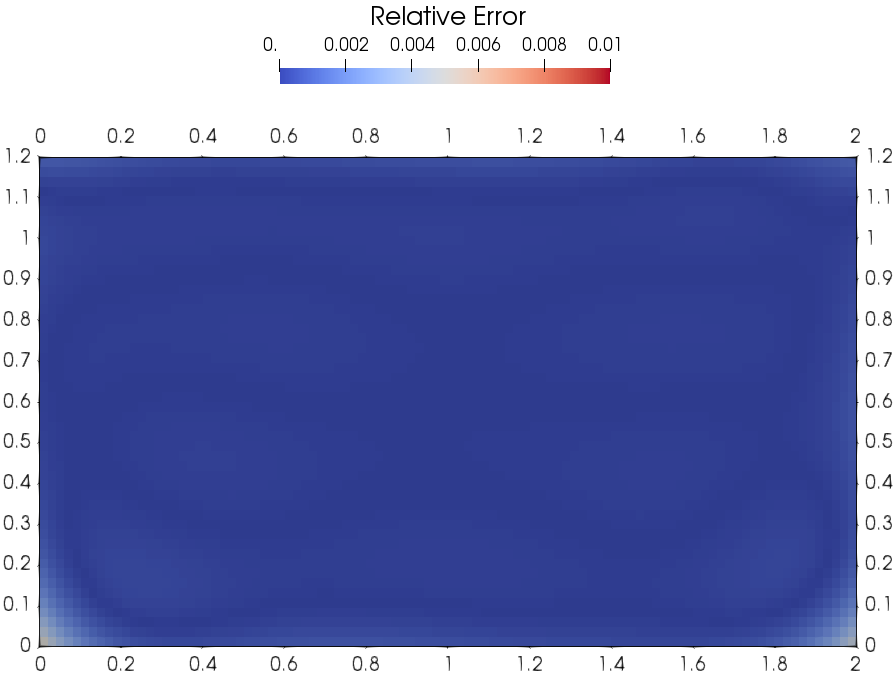}
        \caption{Relative error}
    \end{subfigure}%
    \caption{Estimated heat flux (A) and the respective relative error (B) using the parameterization method with RBF shape parameter $\eta=0.3$ in the industrial benchmark case.}
\label{fig:numericalSteadyBenchmark_paramBC_heatFlux}
\end{figure}

To conclude, there is no relationship between the condition number of the linear system and the obtained results for this industrial benchmark test.
However according to Figure~\ref{fig:numericalBenchmark_y02_paramBC_RBFshapeParameterEffects} (b), we obtain the best results for $\eta = 0.3$.
Figure~\ref{fig:numericalSteadyBenchmark_paramBC_heatFlux} shows the results obtained for this value of the RBF shape parameter.
Then, we use this value in the following tests.

\begin{figure}[!htb]
\centering
\includegraphics[width=0.5\textwidth]{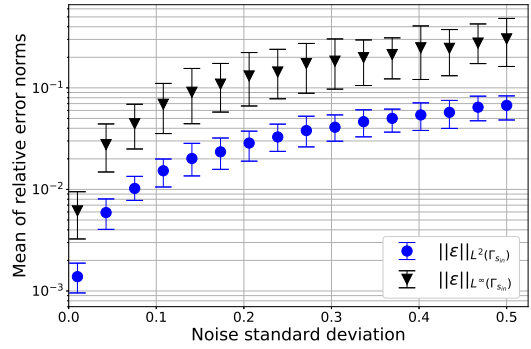}
    \caption{Effect of the measurements noise on the solution of the parameterization method with LU factorization with full pivoting in the industrial benchmark case (90\% quantile bars shown).}
\label{fig:numericalSteadyBenchmark_paramBC_noiseLevelTest_LU}
\end{figure}

We now analyze the effect of noise in the measurement.
Figure~\ref{fig:numericalSteadyBenchmark_paramBC_noiseLevelTest_LU} shows the effect of different noise levels on $L^2$- and $L^\infty$-norms of the relative error (\ref{eq:analyticalBenchmark_errorNorms}) using LU factorization with full pivoting for the solution of (\ref{eq:linSys_parametrizedBC}).
The relative error increases linearly with the noise level.

\begin{figure}[!htb]
    \centering
    \begin{subfigure}[htb]{\linewidth}
        \centering
        \begin{subfigure}[c]{.5\linewidth}
        \centering
            \includegraphics[width=0.9\textwidth]{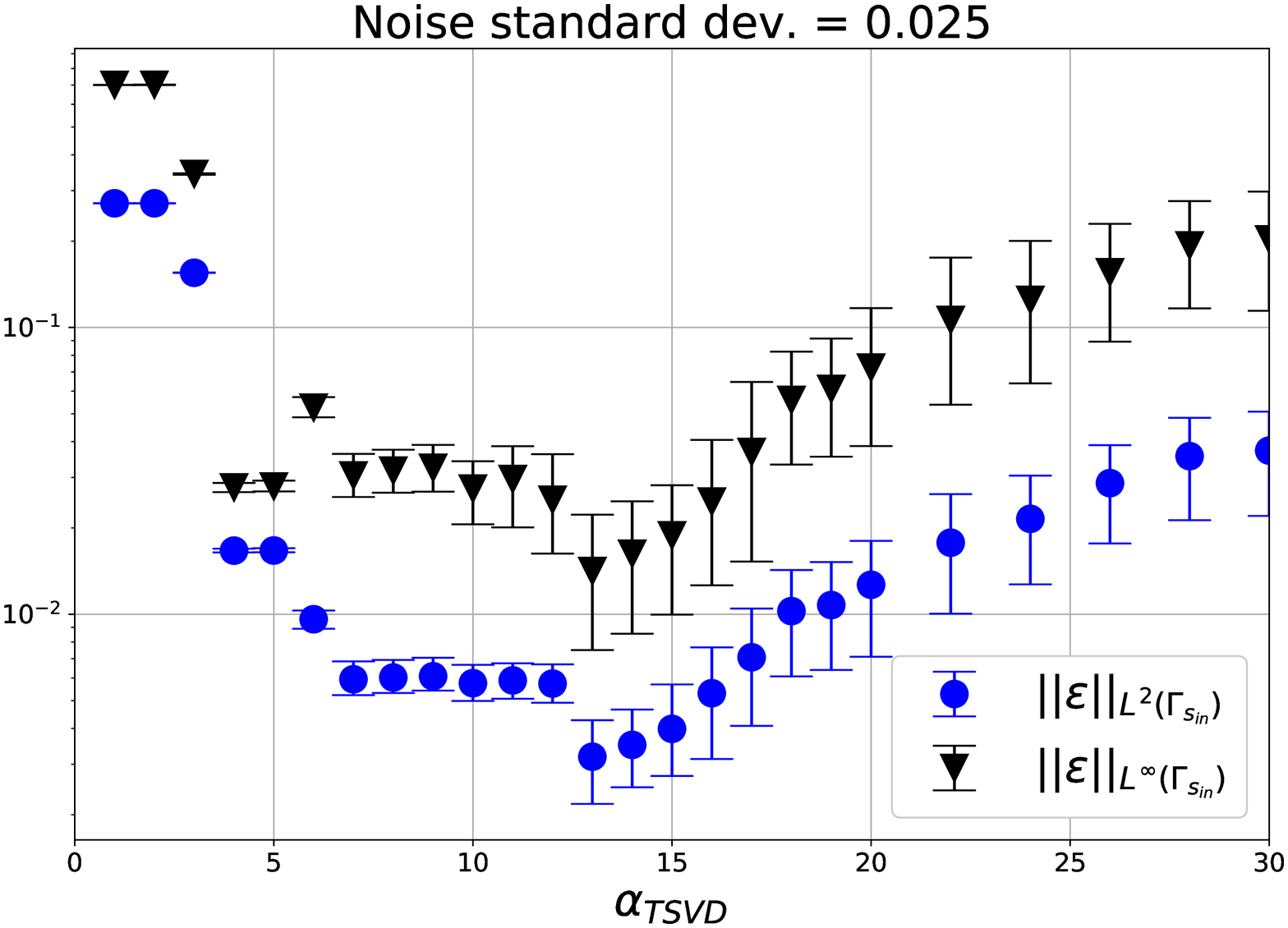}
        \end{subfigure}%
        \begin{subfigure}[c]{.5\linewidth}
        \centering
            \includegraphics[width=0.9\textwidth]{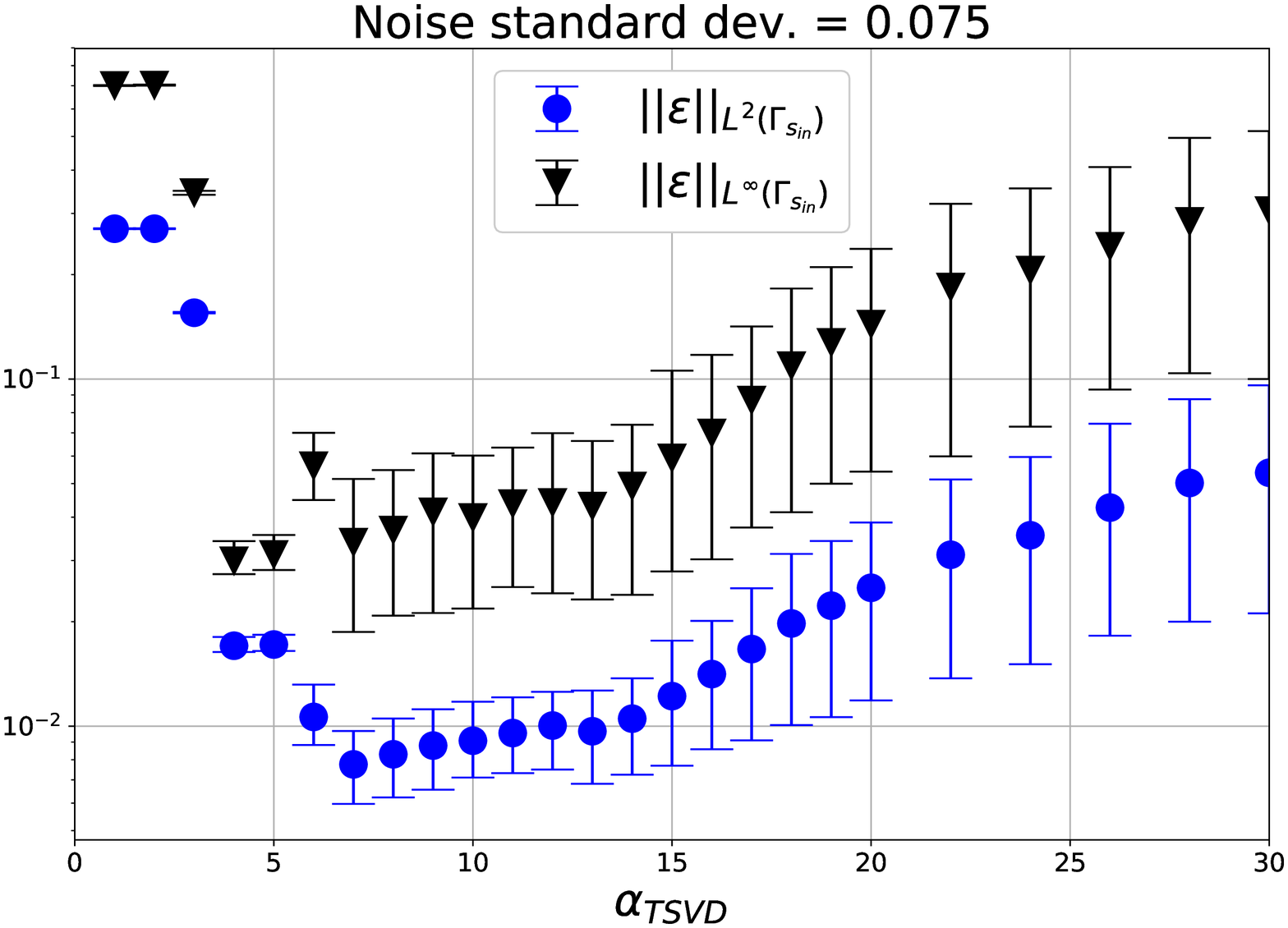}
        \end{subfigure}%
    \end{subfigure}%

    \centering
    \begin{subfigure}[htb]{\linewidth}
        \centering
        \begin{subfigure}[c]{.5\linewidth}
        \centering
            \includegraphics[width=.9\textwidth]{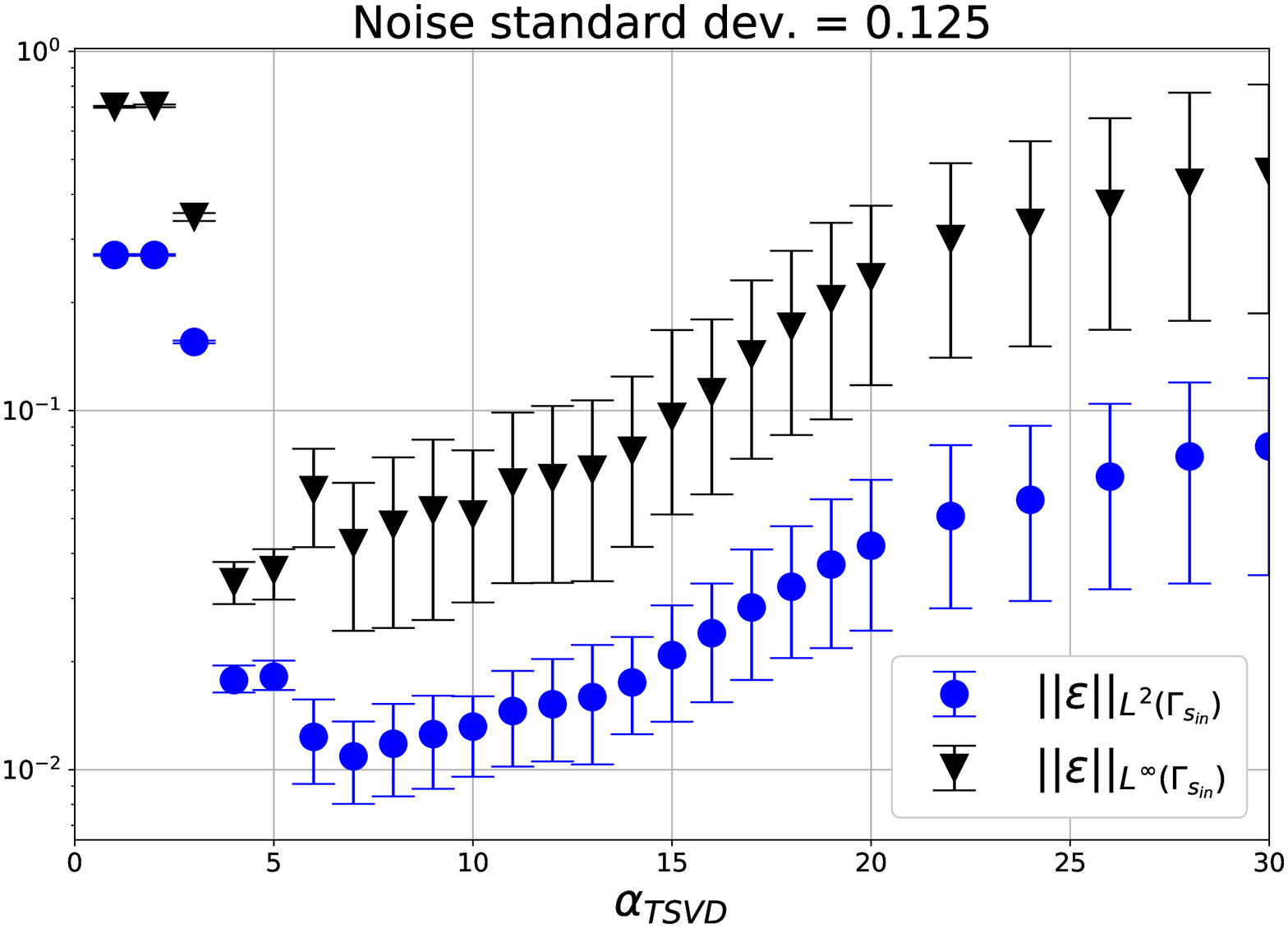}
        \end{subfigure}%
        \begin{subfigure}[c]{.5\linewidth}
        \centering
            \includegraphics[width=.9\textwidth]{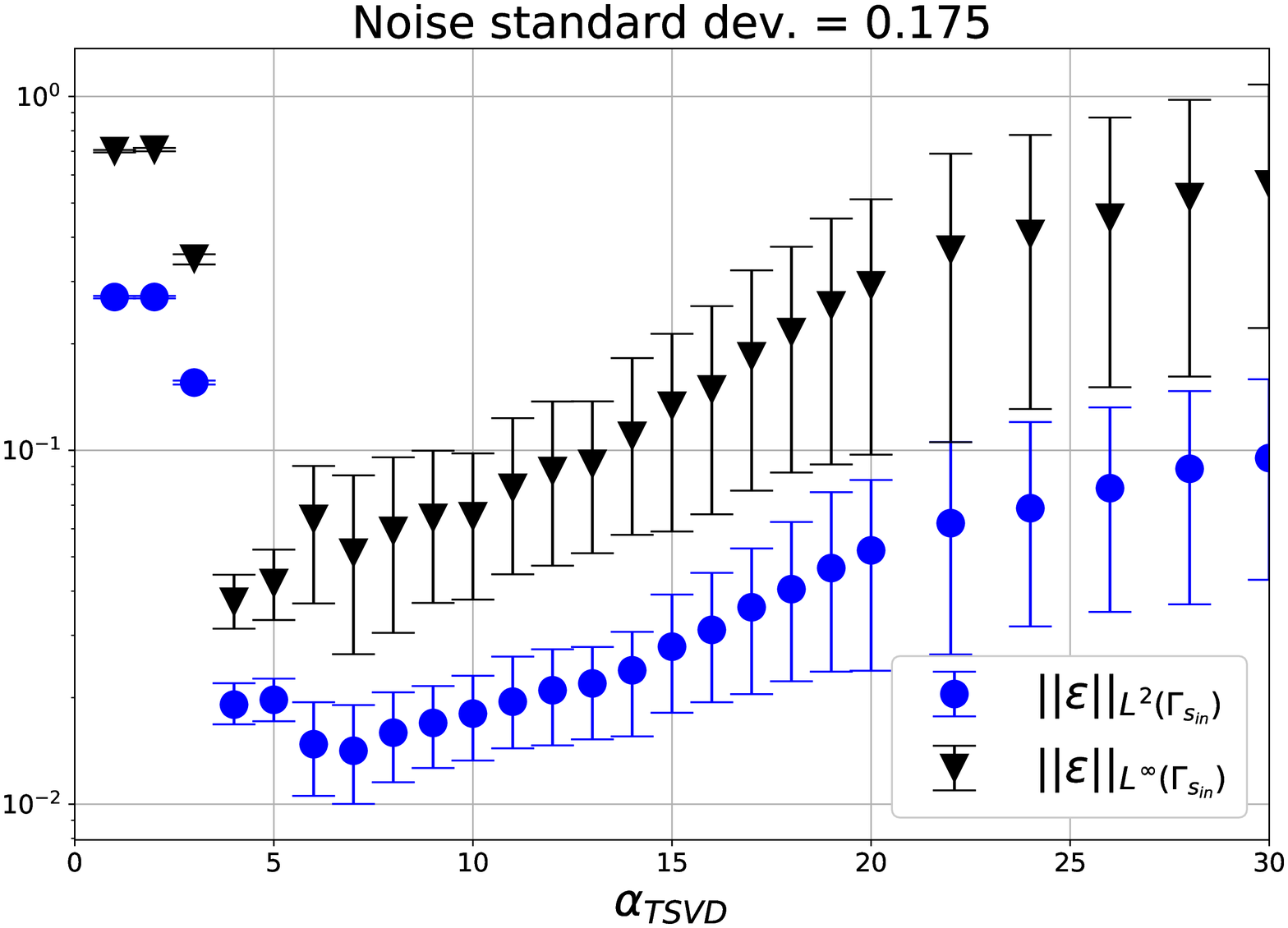}
        \end{subfigure}%
    \end{subfigure}%
    \caption{Effect of the regularization parameter $\alpha_{TSVD}$ using the TSVD in parameterization method for the industrial benchmark case (90\% quantile bars shown).}
\label{fig:numericalSteadyBenchmark_paramBC_TSVDregPar}
\end{figure}

As for the previous benchmark, we test the regularization properties of TSVD on this problem.
Figure~\ref{fig:numericalSteadyBenchmark_paramBC_TSVDregPar} shows the effect of the regularization parameter $\alpha_{TSVD}$ on the $L^2$- and $L^\infty$-norms of the relative error for different values of the noise standard deviation, $\omega$.
As expected, the optimal value of the regularizing parameter $\alpha_{TSVD}$ decreases as the noise increases.
However, for all the considered cases, we are able to achieve a relative error that in the $L^2$-norm is below $2\%$.

\begin{figure}[!htb]
    \begin{subfigure}{\linewidth}
        \centering
        \begin{subfigure}[c]{.5\linewidth}
        \centering
            \includegraphics[width=.9\textwidth]{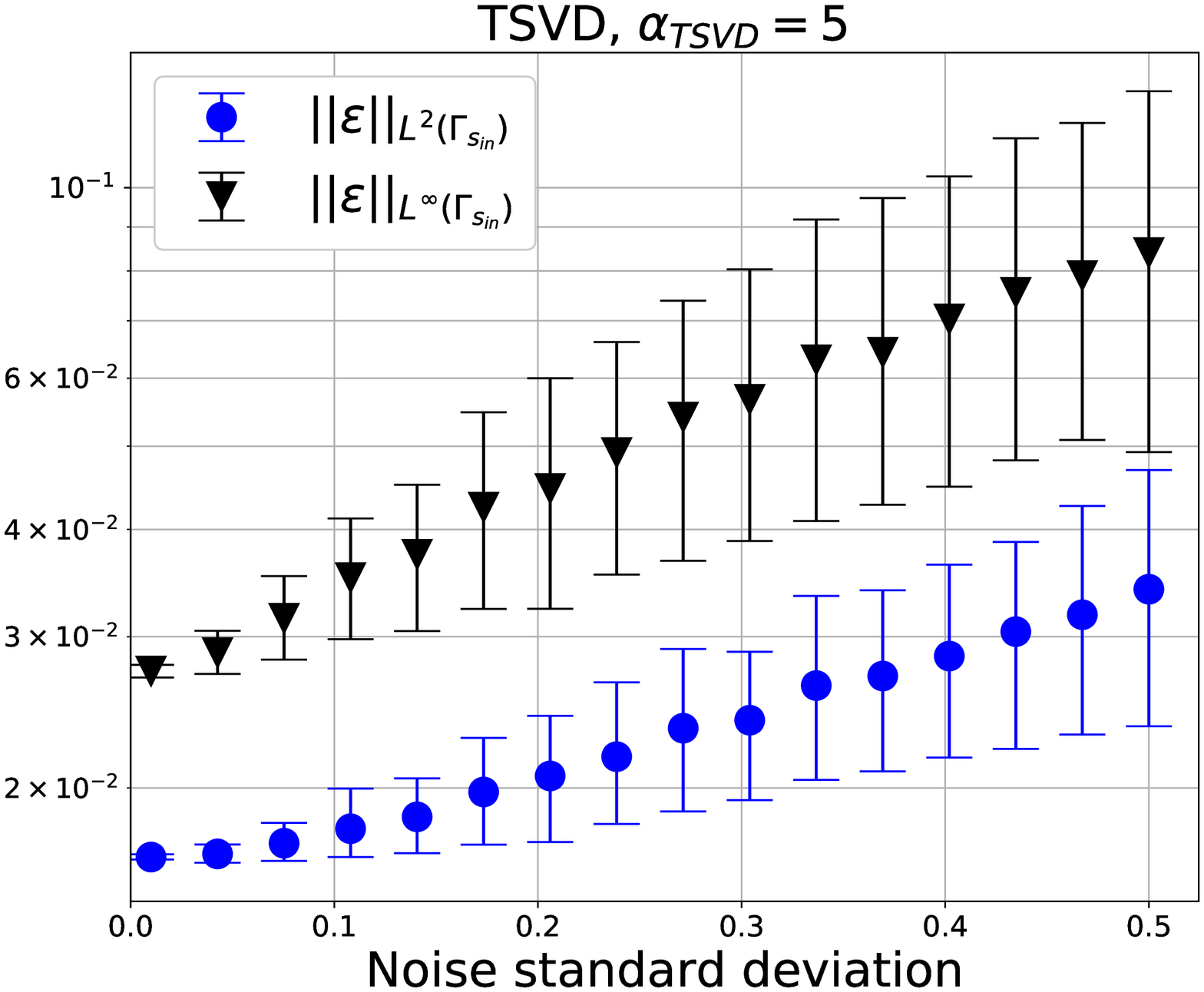}
        \end{subfigure}%
        \begin{subfigure}[c]{.5\linewidth}
        \centering
            \includegraphics[width=.9\textwidth]{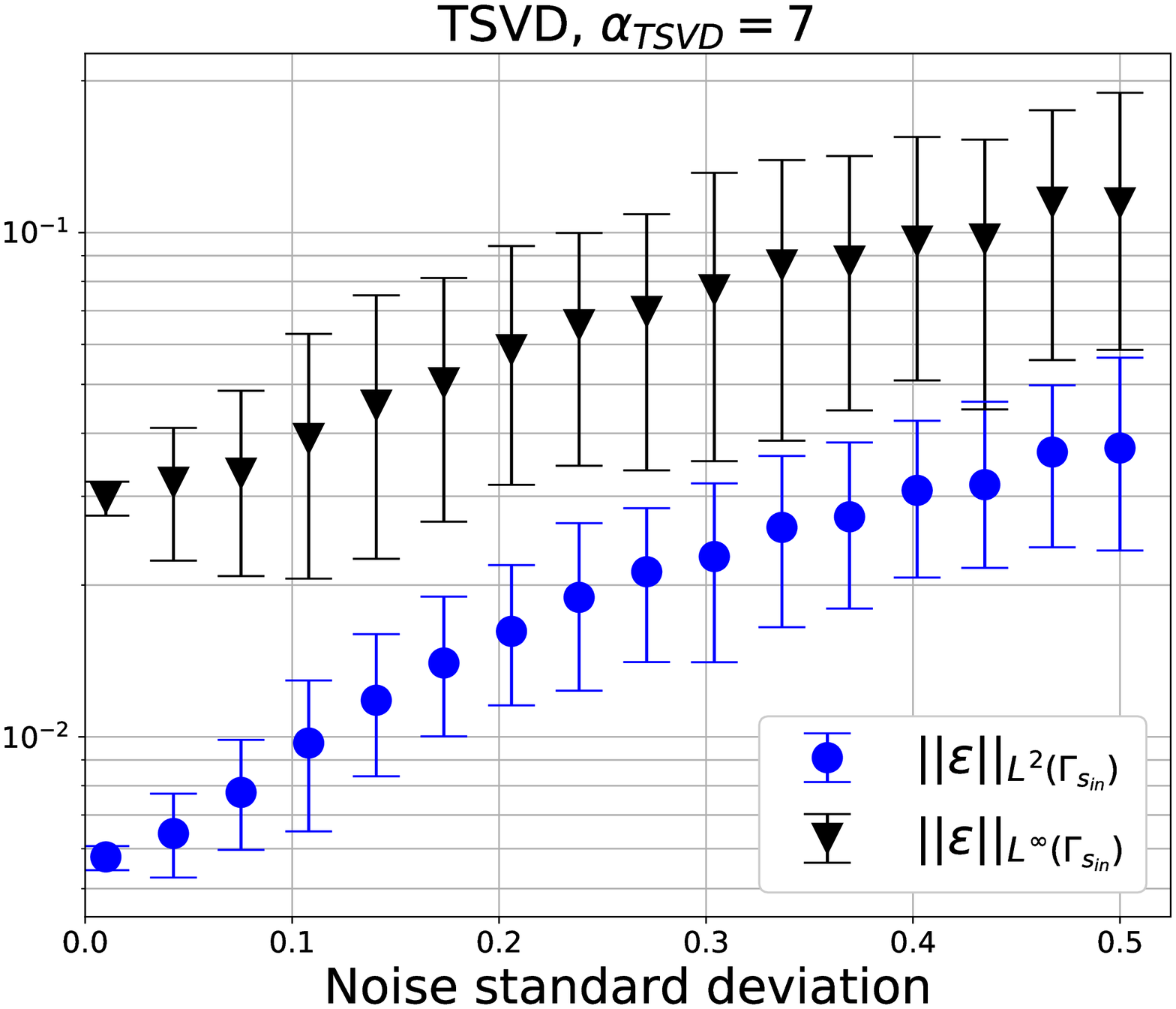}
        \end{subfigure}%
    \end{subfigure}%
    \caption{Behavior of the relative error with respect to the standard deviation of the noise in the measurements using the parameterization method with TSVD regularization in the industrial benchmark case (90\% quantile bars shown).}
\label{fig:numericalSteadyBenchmark_paramBC_noiseLevel_TSVD}
\end{figure}

To conclude, Figure~\ref{fig:numericalSteadyBenchmark_paramBC_noiseLevel_TSVD} shows the behavior of the relative error increasing the measurement noise for $\alpha_{TSVD} = 5$ and $\alpha_{TSVD} = 7$.
Notice that also for severe noise in the thermocouples' measurements, we are able to obtain a valid reconstruction of the boundary heat flux.

\subsection{Inverse Problem with Temperature and Total Heat Flux Measurements}

In this section, we discuss the numerical solution of the inverse Problem~\ref{inverseProblemTotalHeat_3DhcModelSteady} where $T[g](\mathbf{x}_i)$ is the solution of Problem~\ref{prob:3DhcModelSteady} at points $\mathbf{x}_i$, for all $i=1,2,\dots,M$, $g_{true}$ as in Table~\ref{tab:numericalBenchmark_parameters} and $\hat{G} = \int_{\Gamma_{s_{in}}} g_{true} d\Gamma$.

\begin{figure}[!htb]
    \begin{subfigure}[c]{.5\linewidth}
        \includegraphics[width=\textwidth]{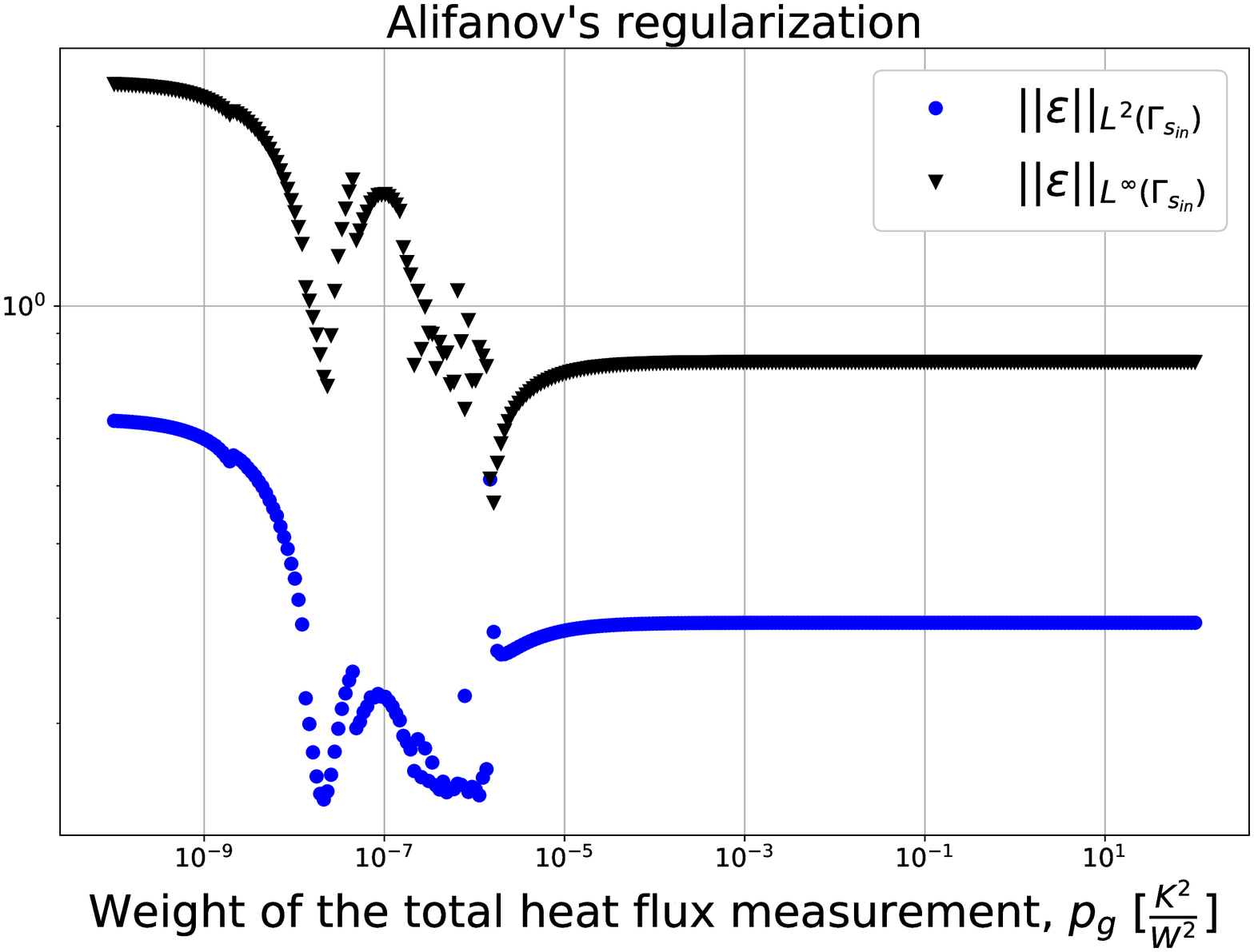}
        \caption{Alifanov's regularization.}
    \end{subfigure}%
    \begin{subfigure}[c]{.5\linewidth}
        \includegraphics[width=\textwidth]{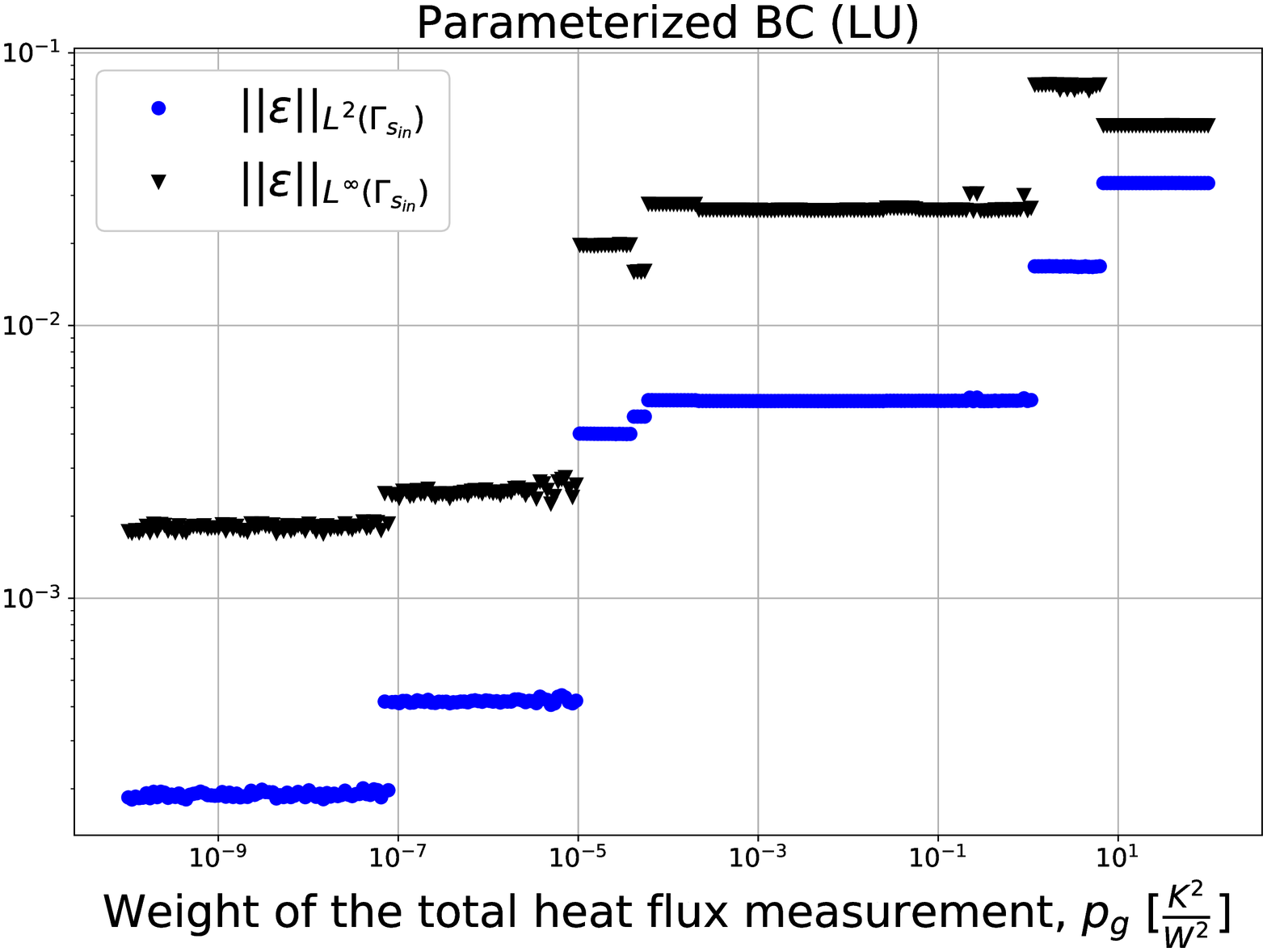}
        \caption{Parameterization method.}
    \end{subfigure}%
    \caption{Effect of the functional weight $p_g$ on the $L^2$- and $L^\infty$-norms of the relative error (\ref{eq:analyticalBenchmark_errorNorms}) for the Alifanov's regularization (A) and the parameterization method with LU factorization (B). 
    The thermocouples' measurements are free of noise.}
    \label{fig:numericalSteadyBenchmark_totalHeat_measureWeight}
\end{figure}

With respect to the previous section, we have one additional parameter: the total heat weight, $p_g$.
Since, it is not possible to set it a priori, we analyze its effects on the solution.
Figure~\ref{fig:numericalSteadyBenchmark_totalHeat_measureWeight} shows the behavior of the $L^2$- and $L^\infty$-norm of the relative error for different values of $p_g$ using Alifanov's regularization and the parameterization method with LU factorization for the solution of the inverse problem.
All these computations are performed without noise in the measurements.

Analyzing Figure~\ref{fig:numericalSteadyBenchmark_totalHeat_measureWeight}, we appreciate a different behavior for the two methods.
Alifanov's regularization improves its results, reaching a minimum of the relative error for $p_g\approx 1e-8$.
Then, the error goes quickly to a plateau in which the estimated heat flux is uniform.
On the other hand, the parameterization method error increases at jumps with $p_g$.

\begin{figure}[!htb]
\centering
    \includegraphics[width=0.5\textwidth]{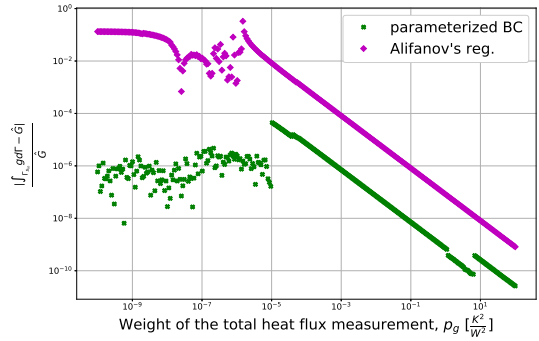}
\caption{Effect of the total heat measurement weight $p_g$ on the relative error in the total heat.}
\label{fig:numericalSteadyBenchmark_totalHeat_totalHeatError}
\end{figure}

Figure~\ref{fig:numericalSteadyBenchmark_totalHeat_totalHeatError} shows the relative error on the total heat flux.
It also provides interesting information.
While both methods linearly improve their performance for $p_g > 1e-5$, the parameterization method shows a very peculiar dependence on the weight for lower values.
However, the parameterization method has a relative error two orders of magnitude smaller with respect to Alifanov's regularization.

%
%
%

\subsection{Conclusions}
In this industrial benchmark case, we tested the methods presented in Sections~\ref{sec:inverseProblem_soloThermocouples} and~\ref{sec:inverseProblem_totalHeat} in an industrial setting.
Alifanov's regularization proved to perform very poorly.
Due to the thermocouples located very close to the boundary $\Gamma_{s_{in}}$, this regularization method overestimates the heat flux close to the measurement points, underestimating it away from the measurements.
Including the total heat flux measurement in the cost functional improves the obtained results, but not to a satisfactory level.

Also in this test case, the parameterization method proved to perform very well providing excellent estimation of the heat flux.
In this case, introducing the total heat measurement caused a degradation of the estimated heat flux.
For this method, the TSVD regularization was used to filter the measurement noise.
It allowed to obtain nice heat flux estimations also in the noisy scenario.

To conclude, Table~\ref{tab:numericalBenchmark_inverseProblem_timeComparison} illustrates the CPU time required for the computations with no error in the measurements and $J_{tol} = 10^{2}~K^2$ in the case of only temperature measurements available.
Notice that all the computations were performed in serial on a Intel\textsuperscript{\textregistered} Core\textsuperscript{\texttrademark} i7-8550U CPU processor.
Recalling that in this application the thermocouples sample at $1~Hz$, the parameterization method allows the real-time estimation of the boundary heat flux.

\begin{table}[!htb]
  \centering
  \caption{Inverse problem CPU time comparison for the industrial benchmark case.}
  \label{tab:numericalBenchmark_inverseProblem_timeComparison}
  \begin{tabular}{ |l|c|c|c| }
    \hline
      & Alifanov's reg. &  \multicolumn{2}{c|}{Parameterized heat flux} \\
    \hline
      & & offline & online \\
    \hline
    CPU time  &  $221~s$ & $121.4~s$ & $0.15~s$\\
    \hline
  \end{tabular}
\end{table}
\section{Conclusions and Future Work}

The objective of this work was to develop a methodology for the real-time estimation of the steel-mold heat flux in CC molds.

We approached this problem by first studying the mold modeling (the direct problem).
With physical considerations on the problem, we justified some simplifying assumptions for the mold model that allowed us to use a steady heat conduction model for the solid portion of the mold.
This model was equipped with convective BCs on the portion of the boundary in contact with the cooling water and a Neumann BC in the portion in contact with the cooling steel.
This latter BC is the heat flux that we want to estimate.

For the setup of the inverse problem, we considered two different measurement settings: having as measurements only the thermocouples' pointwise temperature measurements or having them together with the total boundary heat flux measurement.
For the definition of the inverse problems, we used a deterministic least square approach.

To solve the inverse problems, we used two different methodologies.
The first one is a traditional regularization method called the Alifanov's regularization.
As a second method, we developed an inverse solver that exploit the parameterization of the boundary heat flux.
The latter is very attractive for our problem because it allows for an offline-online decomposition.
It means that we have a computationally expensive offline phase in which we solve several direct problems and a fast online phase that can be computed in real-time.

We finally tested the developed methodologies in two different benchmark cases: an academic test and an industrial one.
In both cases, we tested the quality of the heat flux reconstruction and the robustness of the methods to the measurements noise.
The results shown that the parameterization method outperforms Alifanov's regularization in all the tests.
Moreover, it provided good solutions also in presence of significant noise in the measurements.
Finally, it allows the real-time estimation of the boundary heat flux while Alifanov's regularization cannot be employed in real-time as it is.

In future work, we will focus mainly on two aspects.
First, we will develop a methodology for the real-time solution of this inverse problem in the unsteady case comparing the obtained results to the steady case.
Second, we will move to a Bayesian approach to the inverse problem.\cite{Matthies2016, Cotter2010}
With this approach, we will be able to better deal with the errors not only in the measurements but also in the model.
For the industrial point of view, it is very valuable since allows to conduct uncertainty quantification on the heat flux estimation.
In both cases, to achieve real-time computations we will exploit reduced order modeling techniques.\cite{Rozza2014, Lassila2014, Hesthaven2015, Chinesta2017, Chen2017, Georgaka2019}

\subsection{Acknowledge}
We would like to acknowledge the financial support of the European Union under the Marie Sklodowska-Curie Grant Agreement No. 765374.
We also acknowledge the partial support by the Ministry of Economy, Industry and Competitiveness through the Plan Nacional de I+D+i (MTM2015-68275-R), by the Agencia Estatal de Investigacion through project [PID2019-105615RB-I00/ AEI / 10.13039/501100011033], by the European Union Funding for Research and Innovation - Horizon 2020 Program - in the framework of European Research Council Executive Agency:  Consolidator Grant H2020 ERC CoG 2015 AROMA-CFD project 681447 "Advanced Reduced Order Methods with Applications in Computational Fluid Dynamics" and  INDAM-GNCS project "Advanced intrusive and non-intrusive model order reduction techniques and applications", 2019.
Moreover, we gratefully thank Gianfranco Marconi, Federico Bianco and Riccardo Conte from Danieli \& C.Officine Meccaniche SpA for helping us in better understanding the industrial problem and for the fruitful cooperation.

\bibliographystyle{ama}
\bibliography{references}%

\end{document}